\documentclass[reqno]{amsart}

\usepackage {fullpage}						
\usepackage{csquotes}
\usepackage [utf8]{inputenc}
\usepackage {graphicx}
\usepackage [bookmarks=false, pdfstartview={FitH}, colorlinks=true, linkcolor=blue, citecolor=blue, urlcolor=blue]{hyperref} 
\usepackage {amsfonts,amsmath,amsthm,amssymb}
\usepackage {latexsym}
\usepackage {graphics, graphicx}
\usepackage {fancyhdr}
\usepackage[table]{xcolor}
\usepackage {epstopdf} 						
\usepackage {multicol}						
\usepackage {multirow}  					
\usepackage {soul} 							
\usepackage {datetime} 						
\usepackage {subfigure}						
\usepackage {parskip}						
\usepackage {setspace}						
\usepackage[font=small,textfont=it,width=0.97\textwidth]{caption}		
\usepackage{mathtools}						
\usepackage{tikz,pgfplots}
\usepackage{nicefrac}
\usepackage{bm}
\usepackage{cite}
\usepackage{morefloats}
\usepackage[section]{placeins}
\usepackage{siunitx}

\usetikzlibrary{external}
\usetikzlibrary{matrix}
\newlength\figureheight
\newlength\figurewidth

\theoremstyle{remark}
\newtheorem{remark}{Remark}[section]

\numberwithin{equation}{section}

\allowdisplaybreaks[1]

\renewcommand{\(}{\left(}
\renewcommand{\)}{\right)}

\newcommand{\w}{\bm w}
\newcommand{\x}{\bm x}

\newcommand{\F}{\bm F}

\newcommand{\D}{\mathfrak D}
\newcommand{\R}{\mathfrak R}

\def\softd{{\leavevmode\setbox1=\hbox{d}%
		\hbox to 1.05\wd1{d\kern-0.4ex{\char039}\hss}}}

\newcommand{\rhoprime}{\rho^\prime}
\newcommand{\pprime}{p^\prime}
\newcommand{\uprime}{\bm{u}^\prime}
\newcommand{\thetaprime}{\theta^\prime}
\newcommand{\rhothetaprime}{(\rho\theta)^\prime}
\newcommand{\rhobar}{\bar{\rho}}
\newcommand{\pbar}{\bar{p}}
\newcommand{\ubar}{\bar{\bm{u}}}
\newcommand{\thetabar}{\bar{\theta}}
\newcommand{\rhothetabar}{\overline{\rho\theta}}
\newcommand{\uvec}{\bm{u}}

\newcommand{\cblue}{\color{black}}
\newcommand{\cred}{\color{black}}
\definecolor{darkgreen}{rgb}{0,0.5,0.1}

\graphicspath{{figures/}}

\title{Stochastic Galerkin method for cloud simulation}

\author{A. Chertock}
\address[A. Chertock]{\newline Department of Mathematics, North Carolina State University}
\email{chertock@math.ncsu.edu}

\author{A. Kurganov}
\address[A. Kurganov]{\newline Department of Mathematics, Southern University of Science and Technology and Mathematics Department, Tulane
University}
\email{kurganov@math.tulane.edu}

\author{M. Luk\'a\v{c}ov\'a-Medvi\v{d}ov\'a}
\address[M. Luk\'a\v{c}ov\'a-Medvi{\softd}ov\'a]{\newline Institute of Mathematics, Johannes Gutenberg-University Mainz}
\email{lukacova@uni-mainz.de}

\author{P. Spichtinger}
\address[P. Spichtinger]{\newline Institute of Atmospheric Physics, Johannes Gutenberg-University Mainz}
\email{spichtin@uni-mainz.de}

\author{B. Wiebe}
\address[B. Wiebe]{\newline Institute of Mathematics, Johannes Gutenberg-University Mainz}
\email{b.wiebe@uni-mainz.de}

\begin{document}
	\maketitle
	
	\vspace{-1cm}
	
\begin{abstract}
We develop a stochastic Galerkin method for a coupled Navier-Stokes-cloud system that models dynamics of warm clouds. Our goal is to
explicitly describe the evolution of uncertainties that arise due to unknown input data, such as model parameters and initial or boundary
conditions. The developed stochastic Galerkin method combines the space-time approximation obtained by a suitable finite volume method with
a spectral-type approximation based on the generalized polynomial chaos expansion in the stochastic space. The resulting numerical scheme
yields a second-order accurate approximation in both space and time and exponential convergence in the stochastic space. Our numerical
results demonstrate the reliability and robustness of the stochastic Galerkin method. We also use the proposed method to study the
behavior of clouds in certain perturbed scenarios, for examples, the ones leading to changes in macroscopic cloud pattern as a shift from
hexagonal to rectangular structures.
\end{abstract}
	
\medskip
\section{Introduction}
Clouds constitute one of the most important component in the Earth-atmosphere system. They influence the hydrological cycle and by
interacting with radiation they control the energy budget of the system. However, clouds are one of the most uncertain components, which,
unlike the atmospheric flows, cannot be modeled using first principles of physics.

Clouds are composed by myriads of water particles in different phases (liquid and solid), and thus they need to be described by a large
ensemble in a statistical sense. A common way of obtaining such an ensemble is by using a mass or size distribution, which would lead to a
Boltzmann-type evolution equation. Although there are some approaches available in literature to formulate cloud models in such a way
\cite{beheng2010,khvorostyanov1995,khain_etal2000}, a complete and consistent description is missing. Since measurements of size
distributions of cloud particles are difficult, we are often restricted to averaged quantities such as, for example, mass of water per dry
air (mass concentrations). Therefore, models are often formulated in terms of so-called bulk quantities, that is, mass and number
concentrations of the respective water species. Many cloud processes are necessary to describe the time evolution of the cloud as a
statistical ensemble, that is, particle formation or annihilation, growth/evaporation of particles, collision processes, and sedimentation
due to gravity. For each of the processes, we have to formulate a representative mathematical term in the sense of a rate equation. Although
for some processes the physical mechanisms are quite understood, the formulation of the process rates usually contain uncertain parameters,
thus cloud models come with inherent uncertainty. On the other hand, the initial conditions for atmospheric flows and the embedded clouds
are also not perfectly constrained, leading to uncertainties in the environmental conditions. It is well-known from former studies that
uncertainties in cloud processes and in environmental conditions can lead to drastic changes in simulations, thus these uncertainties
influence predictability of moist atmospheric flows, clouds and precipitation in a crucial way; for instance, the distribution of latent
heat is changed, which in turn can influence frontogenesis \cite{igel_vandenheever2014} or convection
\cite{grabowski2015,marinescu_etal2017}.

For investigations of the impact of these uncertain cloud model parameters as well as the impact of variations in environmental conditions
on atmospheric flows, sensitivity studies are usually carried out. Since one or more parameters are (randomly) varied, the Monte Carlo
approach can be used. This, however, requires a large ensemble of simulations to be conducted, which makes Monte Carlo methods
computationally expensive and requires a very fine sampling of the parameter space and possible environmental conditions. {\cblue In most
practical studies, a much smaller set of ensembles (with about $10-100$ samples only) is used.

In order to improve both the efficiency and accuracy of a numerical method, we choose a different way of representing random variations by
using spectral expansions in the stochastic space.} This approach enables us to investigate the impact of variations in cloud model
parameters and initial conditions on the evolution of moist flows with embedded clouds.

We consider a mathematical model of cloud physics that consists of the Navier-Stokes equations  coupled with the cloud evolution equations
for the water vapor, cloud water and rain. In this  model developed in \cite{cloudpaper,porz18} and presented in
Section \ref{Deterministic_mathematical_model}, the Navier-Stokes equations describe weakly compressible flows with viscous and heat
conductivity effects, while microscale cloud physics is modeled by the system of advection-diffusion-reaction equations.

{\cblue Meteorological applications typically inherit several sources of uncertainties, such  as model parameters, initial and boundary
conditions. Consequently, purely deterministic models are insufficient in such situations and more sophisticated methods need to be applied
to analyze the influence of  uncertainties on numerical solutions.} In this paper, we study a stochastic version of the coupled
Navier-Stokes-cloud model in order to account for uncertainties in input quantities. Our main goal is to design an efficient numerical
method for quantifying uncertainties in solutions of the studied system. In recent years, a wide variety of uncertainty quantification
methods has been proposed and investigated in the context of physical and engineering applications. These  methods include stochastic
Galerkin methods based on generalized polynomial chaos (gPC) \cite{XK02,PDL,TLNE10a,WK06,DPL13}, stochastic collocation methods
\cite{MZ09,WLB09,XH05}, and multilevel Monte Carlo methods \cite{SMS13,MSS13,MS12}. Each of these groups of methods has its own pros and
cons. While results obtained by the Monte Carlo simulations are generally good, the approach is not very efficient due to a large number of
realizations required. Stochastic collocation methods are typically more efficient than the Monte Carlo ones, since they only require
solving the underlying deterministic system at the certain quadrature nodes in the stochastic space. These data are then used to reconstruct
the gPC expansion using an appropriate set of orthogonal polynomials. {\cred The Monte Carlo as well as the stochastic collocation method
fall into a class of the {\em non-intrusive methods}.} Stochastic Galerkin methods offer an alternative approach for computing the gPC
expansion. In general, they are  more rigorous and efficient than the Monte Carlo and collocation ones; see, e.g., \cite{EMPT11}. {\cred
The stochastic Galerkin method is an {\em intrusive method} since it requires changes in the underlying code. In fact,} one needs to solve
a system of PDEs for the gPC expansion coefficients.

We develop a new stochastic Galerkin method for the coupled Navier-Stokes-cloud system. {\cblue As it has already been mentioned above, the
largest source of uncertainties is} cloud physics. {\cblue Therefore,} we restrict our consideration to the case in which the uncertainties
are only in  the cloud physics representation; extension to full stochastic Navier-Stokes-cloud model is left to future studies. Thus, we need to solve the
deterministic Navier-Stokes equations coupled with the PDE system for the gPC expansion coefficients for the cloud variables. Our numerical
method is an extension of the approach proposed in \cite{cloudpaper} for the purely deterministic version of the coupled Navier-Stokes-cloud
system. This method is based on the operator splitting approach, in which the system is split into the macroscopic Navier-Stokes equations
and microscopic cloud model with random inputs. The Navier-Stokes equations are then solved by an implicit-explicit (IMEX) finite-volume
method, while for the cloud equations we develop a stochastic Galerkin method based on the gPC. The resulting gPC coefficient system is
numerically solved by a finite-volume method combined with an explicit Runge-Kutta method with an enlarged stability region \cite{Dumka3}.

The paper is organized as follows. We start in Section \ref{Deterministic_mathematical_model} with the description of the deterministic
Navier-Stokes-cloud model. The numerical method for the deterministic model is presented in Sections \ref{sec4}. In Sections \ref{sec6}, we
report on numerical experiments for well-known meteorological benchmarks---rising warm bubble and Rayleigh-B\'{e}nard convection---for the
deterministic model. We then continue in Section \ref{sec3} with the presentation of the stochastic model, which is followed by the
description of the numerical method (Section \ref{sec5}) and presentation of the numerical experiments (Section \ref{sec7}) for the
stochastic model. Our numerical results clearly demonstrate that the proposed stochastic Galerkin method is capable of quantifying the
uncertainties of complex atmospheric flows.

\section{Deterministic mathematical model}\label{Deterministic_mathematical_model}
We study a mathematical model of cloud dynamics, which is based on the compressible nonhydrostatic Navier-Stokes equations for moist
atmosphere (that is, mixture of ideal gases dry air and water vapor),
\begin{align}
\label{NS_equations}
\rho_t+\nabla\cdot\(\rho\bm{u}\)&=0,\notag\\	
(\rho\bm{u})_t+\nabla\cdot\(\rho\bm{u}\otimes\bm{u}+p\,{\rm Id}-\mu_m\rho\(\nabla\bm{u}+(\nabla\bm{u})^\top\)\)&=
-\rho g\bm{e_3},\\
(\rho\theta)_t+\nabla\cdot\({\rho\theta}\bm{u}-\mu_h\rho\nabla\theta\)&=S_\theta\notag,
\end{align}
and evolution equations for cloud variables,
\begin{align}
\label{cloud_equations}
(\rho q_v)_t+\nabla\cdot\(\rho q_v\bm{u}-\mu_q\rho\nabla q_v\)&=\rho(-C+E),\notag\\
(\rho q_c)_t+\nabla\cdot\(\rho q_c\bm{u}-\mu_q\rho\nabla q_c\)&=\rho(C-A_1-A_2),\\
(\rho q_r)_t+\nabla\cdot\(-v_q\rho q_r\bm{e_3}+\rho q_r\bm{u}-\mu_q\rho\nabla q_r\)&=\rho(A_1+A_2-E).\notag
\end{align}
Here, $\rho$ is the density, {\cblue $\uvec=(u_1,u_2,u_3)^\top$} is the velocity vector, $\theta$ is the moist potential temperature, $p$ is the
pressure, $g$ is the acceleration due to gravity, $\mu_m$ is the dynamic viscosity, $\mu_h$ is the thermal conductivity, and $\mu_q$ is the
cloud diffusivity. {\cred The cloud variables representing the mass concentration of water vapor, cloud droplets and rain drops, $q_v$,
$q_c$, $q_r$, respectively, as well as the right-hand side (RHS) terms $E$, $C$, $A_1$, $A_2$ will be defined below.} We denote by $t$ the
time variable and by $\x$ the space vector; $\x=(x_1,x_2,x_3)$ in the three-dimensional (3-D) and $\x=(x_1,x_3)$ in the two-dimensional
(2-D) cases. Furthermore, $\bm{e_3}=(0,0,1)^\top$ and $\bm{e_3}=(0,1)^\top$ in the 3-D and 2-D cases, respectively. We set $\mu_m=10^{-3}$
and $\mu_h=10^{-2}= \mu_q$. Note that the systems \eqref{NS_equations} and \eqref{cloud_equations} are coupled through the source term
$S_\theta$, which represents the impact of phase changes and will be defined below, see \eqref{S_theta}. The temperature $T$ can be obtained
from the moist adiabatic ideal gas equation
\begin{equation}
T=\frac{R}{R_m}\theta\(\frac{p}{p_0}\)^{\nicefrac{R_m}{c_p}},
\label{ideal_gas_law}
\end{equation}
where $p_0=10^5\,\mathrm{Pa}$ is the reference pressure at sea level. In addition to the usual definition of a potential temperature, we use
$R_m=(1-q_v-q_c-q_r)R+q_v R_v$ with the ideal gas constant of dry air $R=287.05\,\mathrm{\nicefrac{J}{(kg\cdot K)}}$, the gas constant of water
vapor $R_v=461.51\,\mathrm{\nicefrac{J}{(kg\cdot K)}}$ and the specific heat capacity of dry air for constant pressure
$c_p=1005\,\mathrm{\nicefrac{J}{(kg \cdot K)}}$. In order to close the system, we determine the pressure from the equation of state that includes
moisture
\begin{equation}
p=p_0\(\frac{R\rho\theta}{p_0}\)^{\gamma_m}\quad\mbox{with}\quad\gamma_m=\frac{c_p}{c_p-R_m}.
\label{state_equation}
\end{equation}
We note that in the dry case $R_m$ reduces to $R$, $S_{\theta}=0$ and the moist ideal gas equation as well as the moist equation of state
become their dry analogon.
	
In this paper, we restrict our investigations to clouds in the lower part of the troposphere, that is, to clouds consisting of liquid
droplets exclusively. All of the processes involving ice particles are left for future research. For the representation of liquid clouds in
our model we use the so-called single moment scheme, that is, equations for the bulk quantities of mass concentrations of different water
phases. For the representation of the relevant cloud processes we adapt a recently developed cloud model \cite{porz18}. Note that for bulk
models, the process rates cannot be derived completely from first principles. Consequently, some uncertain parameters show up naturally.
This underlies the need of a rigorous sensitivity study which is the goal of the present paper.

Generally, we follow the standard approach in cloud physics modeling for separating hydrometeors of different sizes, as firstly introduced
in \cite{kessler1969}. This relies on the observations that small droplets have a negligible falling velocity. In addition, measurements
indicate two different modes of droplets in the size distribution, which can be associated to small cloud droplets and large rain drops
\cite{warner1969}. Thus, we use the cloud variables $q_c$ and $q_r$ indicating mass concentration of (spatially stationary) cloud droplets
and (falling) rain drops, respectively, and the water vapor concentration $q_v$, that is,
\begin{equation*}
q_\ell=\frac{\mbox{mass of the respective phase}}{\mbox{mass of dry air}}\quad\mbox{for}\quad\ell\in\{v,c,r\}.
\end{equation*}
The rest of this section is devoted to a description of the different terms on the RHS of \eqref{cloud_equations}, which represents the
following relevant cloud processes, {\cblue see \cite{porz18}.}

{\cblue
\subsection{Single particle properties}

\begin{itemize}
\item {\em General properties of a single water particle}

As we exclusively investigate water clouds, we can assume a spherical shape of water particles. For small cloud droplets this is a very good
approximation, while for large rain drops drag effects change their shape \cite{Szakall_etal2009,Szakall_etal2010}. However, for our
investigations of ensembles of rain drops, the spherical shape approximation is appropriate. Thus, mass and radius of droplets are related
by the usual equation
\begin{equation*}
m=\frac{4}{3}\pi\rho_\ell r^3\quad\iff\quad r=\left(\frac{3}{4\pi\rho_\ell}\right)^\frac{1}{3}m^\frac{1}{3}
\end{equation*}
with the liquid water density $\rho_\ell=10^3\,\mathrm{kg\,m^{-3}}$. We make a general assumption that small cloud droplets are stationary,
while large rain drops are accelerated by gravity. After balancing gravity by frictional forces, spherical rain drops fall with a terminal
velocity, depending only on the drop mass and the density of air. According to \cite{porz18}, the terminal velocity for a droplet of mass
$m$ is given by
\begin{equation*}
\hspace*{1cm}v_\tau(m)=\alpha \,m^\beta\left(\frac{m_\tau}{m_\tau+m}\right)^\beta\left(\frac{\rho_*}{\rho}\right)^\frac{1}{2},\quad
\alpha=190.3\,\mathrm{m\,s^{-1}\,kg^{-\beta}},\quad\beta=\frac{4}{15},\quad m_\tau=1.21\cdot10^{-5}\mathrm{kg},
\end{equation*}
with the reference density $\rho_*=1.255\,\mathrm{kg\,m^{-3}}$ at $T_*=288\,\mathrm{K}$ and $p_*=101\,325\,\mathrm{Pa}$. For masses
$m\ll m_\tau$, we can approximate the terminal velocity by $v_\tau=\alpha \,m^\beta\sqrt{\frac{\rho_*}{\rho}}$; this approximation will be
used in the description of the process accretion (collection of cloud droplets by rain drops).

\item {\em Diffusion processes: Growth and evaporation}

Diffusion processes (transfer of water molecules to and from the liquid particle) can be described by the following growth equation:
$$
\frac{dm}{dt}=-4\pi rD_vG\rho(q_*-q_v)\,f_v=-\underbrace{4\pi D_vG\left(\frac{3}{4\pi\rho_l}\right)^\frac{1}{3}}_{=:d}
\rho(q_*-q_v)m^\frac{1}{3}f_v,
$$
where $D_v$ denotes the diffusion constant, $G$ determines corrections due to the latent heat release for phase changes, and $f_v$ is the
ventilation correction for large particles taking into account the effect of flows around the falling spheres. A thermodynamics equilibrium
is determined by the saturation mixing ratio $q_*=q_*(p,T)=\nicefrac{\varepsilon p_s(T)}{p}$ with the saturation water vapor pressure over
a liquid surface $p_s(T)$ given in \cite{murphy_koop2005}. By neglecting curvature effects, water particles grow for $q_v>q_*$ and evaporate
for $q_v<q_*$, respectively. The diffusion constant is given according to \cite{pruppacher_klett2010}:
\begin{equation*}
D_v=D_{v0}\left(\frac{T}{T_0}\right)^{1.94}\frac{p_0}{p},\quad D_{v0}=2.11\cdot 10^{-5}\mathrm{m^2\,s^{-1}},\quad T_0=273.15\,\mathrm{K},
\quad p_0=p_*=101\,325\,\mathrm{Pa}
\end{equation*}
and the impact of latent heat release is described by
\begin{equation*}
G=\left[ \left(\frac{L}{R_vT}-1\right)\frac{Lp_s(T)}{R_vT^2}\frac{D_v}{K_T}+1\right]^{-1},
\end{equation*}
where the latent heat of vaporisation $L=2.53\cdot10^6\mathrm{J\,kg^{-1}}$ and the heat conduction of dry air is (see \cite{dixon2007})
\begin{equation*}
K_T=\frac{a_KT^\frac{3}{2}}{T+b_K10^{\frac{c_K}{T}}},\quad a_K=0.002646\,\mathrm{W\,m^{-1}\,K^{-\frac{5}{2}}},\quad b_K=245.4\,\mathrm{K},
\quad c_K=-12\,\mathrm{K}.
\end{equation*}
Ventilation of large spherical particles of radius $r$ can be taken into account using an empirical ventilation coefficient
\begin{equation*}
f_v=a_v+b_vN_{\rm Sc}^\frac{1}{3}N_{\rm Re}^\frac{1}{2},\quad a_v=0.78,\quad b_v=0.308,
\end{equation*}
where the Schmidt and Reynolds numbers are defined as
\begin{equation}
N_{\rm Sc}=\frac{\mu}{\rho D_v}\quad\mbox{and}\quad N_{\rm Re}=\frac{\rho}{\mu}v_\tau(2r),
\label{2.4a}
\end{equation}
respectively. In \eqref{2.4a}, $\mu$ is the dynamic viscosity of air, which is expressed according to \cite{dixon2007} by
\begin{equation*}
\mu=\frac{\mu_0 T^\frac{3}{2}}{T+T_\mu},\quad\mu_0=1.458\cdot10^{-6}\mathrm{s\,Pa\,K^{-\frac{1}{2}}},\quad T_\mu=110.4\,\mathrm{K}.
\end{equation*}

For cloud droplets, we neglect the ventilation correction, thus the mass rate of diffusion for a cloud droplet of mass $m_c$ can be
expressed as
\begin{equation*}
\frac{dm_c}{dt}=d\rho(q_v-q_*)m_c^\frac{1}{3}.
\end{equation*}
For rain drops, growth due to the diffusion is negligible, and thus we obtain the mass rate for rain drops of mass $m_r$ as
\begin{equation*}
\hspace*{1.3cm}\frac{dm_r}{dt}=-d\rho(q_*-q_v)_+\left[a_Em_r^\frac{1}{3}+b_Ev_\tau(m_r)^\frac{1}{2}m_r^\frac{1}{2}\right],\quad
a_E=a_v,\quad b_E=b_v\left( \frac{\mu}{\rho D_v}\right)^\frac{1}{3}\sqrt{\frac{2\rho}{\mu}}\left(\frac{3}{4\pi\rho_\ell}\right)^\frac{1}{6}.
\end{equation*}
Here, $(\cdot)_+:=\max(\cdot,0)$ denotes the positive part.

\item {\em Collision of rain drops with cloud droplets: Accretion}

A spherical rain drop of mass $m_r$ (radius $r$) falls with terminal velocity $v_\tau(m_r)$ through a volume $V=\pi r^2v_\tau(m_r)\Delta t$
(during a time interval $\Delta t$) and collects cloud droplets of total mass $M_c=V\rho q_c$. Thus, the corresponding growth rate of the
rain drop is given by
\begin{equation*}
\frac{dm_r}{dt}=k_2'\rho\pi q_cv_\tau(m_r)\left(\frac{3}{4\pi\rho_\ell}\right)^\frac{2}{3}m_r^\frac{2}{3}
\end{equation*}
with an efficiency $k_2'>0$.
\end{itemize}

\subsection{Ensemble/collective properties}\label{section_variables}
For the description of clouds as an ensemble of water particles, we would have to introduce such averaged quantities as mass concentrations
(as described above, that is, $q_c$ and $q_r$) as well as number concentrations of cloud droplets, $n_c$, and rain drops, $n_r$. Since we do
not extend the systems of equations for these two quantities, we introduce relations between mass and number concentrations in order to keep
the main effects in a simplified way.

\begin{itemize}
\item {\em Formation of cloud droplets: Activation}

Cloud droplets can be formed by the activation of so-called cloud condensation nuclei (CCN). Liquid aerosol particles can grow by water
vapor uptake to larger sizes; this effect can be described by the K\"ohler theory; see, e.g., \cite{koehler1936,petters_kreidenweis2007}. As
described in detail in \cite{porz18}, we represent the cloud droplet number concentration $n_c$ by a nonlinear relation
\begin{equation*}
n_c=q_c\frac{N_\infty}{q_c+N_\infty m_0}\coth\(\frac{q_c}{N_0m_0}\).
\end{equation*}
Here, $N_\infty$ denotes the maximum number of CCN (depending on environmental conditions, e.g. clean or polluted air), $m_0$ can be
interpreted as the activation mass of cloud droplets, and $N_0$ is the approximated number of activated droplets at $q_v=q_*$. In our
investigations, we set these three parameters to the following values:
\begin{equation*}
N_\infty=8\cdot 10^8\,\mathrm{kg^{-1}},\quad m_0=5.236\cdot 10^{-16}\mathrm{kg},\quad N_0=10^3\,\mathrm{kg^{-1}}.
\end{equation*}
For the initialization of the cloud droplet production, we introduce an additional factor in case of supersaturation
\begin{equation*}
C_{\rm act}=N_0\,d\rho(q_v-q_*)_+\,m_0^\frac{1}{3}.
\end{equation*}

\item {\em Relation between number and mass concentration for rain drops}

In contrast to the formulation in \cite{porz18}, we do not include another equation for the number concentration of rain drops. In a similar
way as for cloud activation, we use a relation between $n_r$ and $q_r$, that is, a closure of the form $n_r=f(q_r,c)$. Since we implicitly
assume that the rain drops are distributed according to their size, this approach should be used for mimicking the shape of the distribution
in a proper way. We propose the (non)linear relation
\begin{equation*}
n_r=c_rq_r^\gamma,\quad0<\gamma\le1.
\end{equation*}
Assuming a constant mean mass of rain drops $\overline{m}_r$, we can determine the constants as $c_r=\overline{m}_r^{-1}$ and $\gamma=1$.
This approach would be meaningful for the case of a symmetric size distribution of rain droplets, centered around the mean mass. However, it
is well-known that size distributions of rain are usually skew to larger sizes, thus a linear relation is not appropriate. For sizes of rain
drops often an exponential distribution is assumed, this leads to an exponent $\gamma=\frac{1}{4}$ and a coefficient
$c_r=c_{r0}\rho^{-\frac{3}{4}}$ ($c_{r0}=23752.6753\,\mathrm{kg}^{-\frac{1}{4}}\mathrm{m}^\frac{3}{4}$, $[c_r]=\mathrm{kg}^{-1}$) as derived
in Appendix A.

\item {\em Rates for diffusion processes}

For the mass mixing ratios $q_c$ and $q_r$, the rates for the diffusion processes are given by multiplication of the single particle rates
by the number concentration of the respective particles, namely:
\begin{equation*}
\frac{dq_c}{dt}=n_c\frac{dm_c}{dt},\quad\frac{dq_r}{dt}=n_r\frac{dm_r}{dt} .
\end{equation*}
Applying the relations between the mass and number concentrations as stated above, we obtain
\begin{equation*}
C_1=\frac{dq_c}{dt}=n_c d\rho(q_v-q_*)m_c^\frac{1}{3}\stackrel{m_c=\frac{q_c}{n_c}}{=}d\rho(q_v-q_*)
\left(\frac{N_\infty}{q_c+N_\infty m_0}\coth\(\frac{q_c}{N_0m_0}\)\right)^\frac{2}{3}q_c
\end{equation*}
and
\begin{align*}
E=-\frac{dq_r}{dt}=n_r d\rho&(q_*-q_v)_+\left[a_Em_r^\frac{1}{3}+b_Ev_\tau(m_r)m_r^\frac{1}{2}\right]\\
\stackrel{m_r=\frac{q_r}{n_r}}{=}d\rho&(q_*-q_v)_+\left[a_Ec_r^\frac{2}{3}q_r^{\frac{1}{3}+\gamma\frac{2}{3}}+
b_Ec_r^\frac{1}{2}v_\tau(q_r)^\frac{1}{2}q_r^{\frac{1}{2}+\gamma\frac{1}{2}}\right]
\end{align*}
using a reformulated terminal velocity
\begin{equation*}
v_\tau(m_r)=v_\tau(q_r)=\alpha q_r^\beta\left(\frac{m_\tau}{q_r+m_\tau c_rq_r^\gamma}\right)^\beta
\left(\frac{\rho_*}{\rho}\right)^\frac{1}{2}.
\end{equation*}
Note, that the rates for activation and diffusion growth of cloud droplets are combined in the model formulation, that is,
$C=C_1+C_{\rm act}$.

\item {\em Rate for accretion}

For the rate of accretion of rain water, we obtain
\begin{equation*}
A_2=\frac{dq_r}{dt}=n_r\frac{dm_r}{dt}=k_2\rho\pi c_r^{\frac{1}{3}}\left(\frac{3}{4\pi\rho_\ell}\right)^\frac{2}{3}q_cv_\tau(q_r)
q_r^{\frac{2+\gamma}{3}}.
\end{equation*}
Note that for compensating effects of the averaging the parameter $k_2$ can be adjusted (such that $k_2=0.8\,\mathrm{kg}\neq k_2'$) and the
impact of the uncertainty of this parameter is of high interest, since it cannot be measured or derived from the first principles.

\item {\em Collision of cloud droplets, forming a rain drop: Autoconversion}

Beside the growth of an existing rain drop by collecting cloud droplets, a rain drop can be formed by the collision of two cloud droplets.
According to \cite{porz18}, we formulate the growth rate similarly to population models, namely:
\begin{equation*}
A_1=\frac{dq_c}{dt}=k_1\frac{\rho q_c^2}{\rho_\ell}.
\end{equation*}
Note that the coefficient $k_1$ cannot be measured or derived from the first principles. It is a free parameter, which must be fixed using
parameter estimations. Thus, the impact of the uncertainty of this parameter is of high interest. In our deterministic experiments, we
choose $k_1=0.0041\,\mathrm{kg\,s}^{-1}$, as indicated in \cite{porz18}.

\item {\em Sedimentation of rain mass mixing ratio}

We have introduced an additional convection term into the equation for the evolution of $q_r$ in \eqref{cloud_equations}, that is, the term
$\nabla\cdot(-v_q\rho q_r\bm{e_3})$, where
\begin{equation*}
v_q=v_q(q_r)=\alpha q_r^\beta \left(\frac{m_\tau}{q_r+m_\tau c_rq_r^\gamma}\right)^\beta\left(\frac{\rho_*}{\rho}\right)^\frac{1}{2}.
\end{equation*}
The parameter $\alpha$ can be derived empirically, but the influence of uncertainty in $\alpha$ is of high interest.
\end{itemize}
}

Note that {\cblue activation and diffusion} processes are formulated explicitly, in contrast to the usual approach of saturation adjustment
(see, e.g., \cite{LV11}), which is less accurate, but commonly used in operational weather forecast models. This explicit formulation
introduces stiffness caused by modeling cloud processes on the RHS of the cloud equations with fractional exponents between $-1$ and $1$.
{\cblue In order to avoid the root evaluation for an argument that is close to zero, we introduce a cut-off function} and replace
$\zeta^\xi$, $\xi\in(-1,1)$, with
\begin{align*}
\begin{cases}
\zeta^\xi,\quad&\mbox{if}~\zeta>10^{-16},\\
0,\quad&\mbox{otherwise}.
\end{cases}
\end{align*}

{\cblue Due to phase changes (activation and growth/evaporation of water particles) latent heat is released or absorbed. These processes are
modeled by the source term in \eqref{NS_equations}:
\begin{equation}
S_\theta=\rho\frac{L\theta}{c_p T}\(C-E\).
\label{S_theta}
\end{equation}
}

Solving the Navier-Stokes equations \eqref{NS_equations} in a weakly compressible regime is known to cause numerical instabilities due to
the multiscale effects. We follow the approach typically used in meteorological models, where
the dynamics of interest is described by a
perturbation of a background state, which is the hydrostatic equilibrium. The latter expresses a balance between the gravity and pressure
forces. Denoting by $\pbar$, $\rhobar$, $\ubar=0$, $\thetabar$ and $\rhothetabar$ the respective background state, the hydrostatic
equilibrium satisfies
\begin{equation*}
\frac{\partial \pbar}{\partial x_3}=-\rhobar g,\quad S_\theta=0,
\end{equation*}
where $\pbar$ is obtained from the equation of state \eqref{state_equation}
\begin{equation}
\pbar=p(\rhothetabar)=p_0\(\frac{R\rhothetabar}{p_0}\)^{\gamma_m}.
\label{p_equ}
\end{equation}
Let $\pprime$, $\rhoprime$, $\uprime$, $\thetaprime$ and $\rhothetaprime$ stand for the corresponding perturbations of the equilibrium
state, then
\begin{equation*}
p=\pbar+\pprime,~~\rho=\rhobar+\rhoprime,~~\theta=\thetabar+\thetaprime,~~\uvec=\uprime,~~
\rho\theta=\rhobar\thetabar+\rhobar\thetaprime+\rhoprime\thetabar+\rhoprime\thetaprime=\rhothetabar+\rhothetaprime.
\end{equation*}
The pressure perturbation $\pprime$ is derived from \eqref{state_equation} and \eqref{p_equ} using the following Taylor expansion
\begin{equation*}
p(\rho\theta)\approx p(\rhothetabar)+\frac{\partial p}{\partial(\rho\theta)}\(\rho\theta-\rhothetabar\)=
\pbar+\gamma_mp_0\(\frac{R\rhothetabar}{p_0}\)^{\gamma_m}\frac{\rhothetaprime}{\rhothetabar},
\end{equation*}	
which results in
\begin{equation*}
\pprime\approx\gamma_mp_0\(\frac{R\rhothetabar}{p_0}\)^{\gamma_m}\frac{\rhothetaprime}{\rhothetabar}.
\end{equation*}
	
{\cblue Consequently, a physical motivation from the hydrostatic balance state leads to a numerically preferable scaling and the
perturbation formulation of the Navier-Stokes equations \eqref{NS_equations}} then reads as
\begin{align}
\label{NS_equations_pert}
\rhoprime_t+\nabla\cdot\(\rho\bm{u}\)&=0,\notag\\	
(\rho\bm{u})_t+\nabla\cdot\(\rho\bm{u}\otimes\bm{u}+\pprime\,{\rm Id}-\mu_m\rho\(\nabla\bm{u}+(\nabla\bm{u})^\top\)\)&=-\rhoprime g\bm{e_3},
\\
\rhothetaprime_t+\nabla\cdot\({\rho\theta}\bm{u}-\mu_h\rho\nabla\theta\)&=S_\theta.\notag
\end{align}	
For alternative representations of cloud dynamics and
their numerical investigations, we refer the reader to \cite{kroner,COSMO} and references therein.

\section{Numerical scheme for the deterministic model}\label{sec4}
The numerical approximation of the coupled model \eqref{NS_equations_pert}, \eqref{cloud_equations} is based on the second-order Strang
operator splitting. Therefore, we split the whole system into the macroscopic Navier-Stokes flow equations and the microscopic cloud
equations. The Navier-Stokes equations \eqref{NS_equations_pert} are approximated by an IMEX finite-volume method and the cloud equations
\eqref{cloud_equations} are approximated by a finite-volume method in space and an explicit Runge-Kutta method with an enlarged stability
region in time.

\subsection{Operator form}
Let $\w:=(\rhoprime,\rho \uvec,\rhothetaprime)^\top$ and $\w_q:=(\rho q_v,\rho q_c,\rho q_r)^\top$ denote the solution vectors of
\eqref{NS_equations_pert} and \eqref{cloud_equations}, respectively. Then, the coupled system can be written as
\begin{align*}
\w_t&=-\nabla\cdot\F(\w)+\D(\w)+\R(\w),\\
(\w_q)_t&=-\nabla\cdot\F_q(\w_q)+\D_q(\w_q)+\R_q(\w_q),
\end{align*}
where $\F$ and $\F_q$ are advection fluxes and $\D$, $\R$ and $\D_q$, $\R_q$ denote the diffusion and reaction operators of the respective
systems. They are given by
\begin{align}
\label{NS_operator}
\F(\w) &:=\(\rho\bm{u},\rho\bm{u}\otimes\bm{u}+p'\,{\rm Id},\rho\theta\bm{u}\)^\top,\notag\\
\D(\w)&:=\(0,\nabla\cdot(\mu_m\rho(\nabla\bm{u}+(\nabla\bm{u})^\top)),\nabla\cdot(\mu_h\rho\nabla\theta)\)^\top,\\
\R(\w)&:=\(0,-\rho^\prime g\bm{e_3},S_\theta\)^\top,\notag
\end{align}
\begin{align*}	
\phantom{mi}\F_q(\w_q)&:=\(\rho q_v\bm{u},\rho q_c\bm{u},\rho q_r\bm{u}-v_q\rho q_r\bm{e_3}\)^\top,\\
\D_q(\w_q)&:=\(\nabla\cdot(\mu_q\rho\nabla q_v),\nabla\cdot(\mu_q\rho\nabla q_c),\nabla\cdot(\mu_q\rho\nabla q_r)\)^\top,\\
\R_q(\w_q)&:=\(-C+E,C-A_1-A_2,A_1+A_2-E\)^\top.
\end{align*}
In order to derive an asymptotically stable, accurate and computational efficient scheme for the Navier-Stokes equations, we first split the
equations into linear and nonlinear parts; see \cite{bispen1,cloudpaper} and references therein. Consequently, we introduce
\begin{itemize}
\item $\F(\w)=\F_L(\w)+\F_N(\w)$ with $\F_L(\w):=\(\rho\bm{u},p'\,{\rm Id},\bar\theta\rho\bm{u}\)^\top$ and
$\,\F_N(\w):=\(0,\rho\bm{u}\otimes\bm{u},\theta'\rho\bm{u}\)^\top$;
\item $\D(\w)=\D_L(\w)+\D_N(\w)$ with
\begin{align*}
\D_L(\w)&:=\(0,\mu_m(\Delta(\rho\bm{u})+\nabla(\nabla\cdot(\rho\bm{u}))),\mu_h\Delta(\rho\theta)'\)^\top~\mbox{and}\\
\D_N(\w)&:=\(0,-\mu_m((\Delta\rho)\bm{u}+(D^2\rho)\bm{u}+\nabla\bm{u}\nabla\rho+\nabla\rho\nabla\cdot\bm{u}),
\mu_h(\Delta(\overline{\rho\theta})-\theta\Delta\rho-\nabla\rho\cdot\nabla\theta)\)^\top;
\end{align*}
\item $\R(\w)=\R_L(\w)+\R_N(\w)$ with $\R_L(\w):=\(0,-\rho'g\bm{e_3},0\)^\top$ and $\,\R_N(\w):=\(0,0,S_\theta\)^\top$.
\end{itemize}
We would like to point out that the choice of the linear and nonlinear operators is crucial. We choose the linear part to model linear
acoustic and gravitational waves as well as linear viscous fluxes. The nonlinear part describes nonlinear advective effects together with
the remaining nonlinear viscous fluxes and the influence of the latent heat. We will use the following notation:
$$
\mathcal{L}:=-\nabla\cdot\F_L(\w)+\D_L(\w)+\R_L(\w)\quad\mbox{and}\quad\mathcal{N}:=-\nabla\cdot\F_N(\w)+\D_N(\w)+\R_N(\w).
$$

\subsection{Discretization in space}\label{discr_space}
The spatial discretization is realized by a finite-volume method. We take a cuboid computational domain $\Omega\subset\mathbb{R}^d$, which
is divided into $N$ uniform Cartesian cells. The cells are labeled in a certain order using a single-index notation. For simplicity of
notation, we assume that the cells are cubes with the sides of size $h$ so that $|C_i|=h^d$. We also introduce the notation $S(i)$ for the
set of all neighboring cells of cell $C_i$, $i=1,\ldots,N$.

We assume that at a certain time $t$ the approximate solution is realized in terms of its cell averages
\begin{align*}
\w_i(t)\approx\frac{1}{h^d}\int\limits_{C_i}\w(\x,t)\,{\rm d}\x\quad\mbox{and}\quad
(\w_q)_i(t)\approx\frac{1}{h^d}\int\limits_{C_i}\w_q(\x,t)\,{\rm d}\x,\quad i =1,\ldots,N.
\end{align*}
In order to simplify the notation, we will now omit the time dependence of $\w_i(t)$ and $(\w_q)_i(t)$. Next, we introduce the notation
$\w_h:=\{\w_i\}_{i=1}^N$ and $(\w_q)_h:=\{(\w_q)_i\}_{i=1}^N$ and consider the following approximation of the advection, diffusion and
reaction operators:
\begin{align*}
\mathcal{A}_i(\w_h)=(\mathcal{A}_L)_i(\w_h)+(\mathcal{A}_{\cblue N})_i(\w_h)&\approx
\frac{1}{h^d}\int\limits_{C_i}\nabla\cdot\F_L(\w(\x,t))\,{\rm d}\x+\frac{1}{h^d}\int\limits_{C_i}\nabla\cdot\F_N(\w(\x,t))\,{\rm d}\x,\\
\mathcal{D}_i(\w_h)=(\mathcal{D}_L)_i(\w_h)+(\mathcal{D}_N)_i(\w_h)&\approx\frac{1}{h^d}\int\limits_{C_i}\D_L(\w(\x,t))\,{\rm d}\x+
\frac{1}{h^d}\int\limits_{C_i}\D_N(\w(\x,t))\,{\rm d}\x,\\
\mathcal{R}_i(\w_h)=(\mathcal{R}_L)_i(\w_h)+(\mathcal{R}_N)_i(\w_h)&\approx\frac{1}{h^d}\int\limits_{C_i}\R_L(\w(\x,t))\,{\rm d}\x+
\frac{1}{h^d}\int\limits_{C_i}\R_N(\w(\x,t))\,{\rm d}\x.
\end{align*}
Analogous notation will be used for the approximations $(\mathcal{A}_q)_i(\w_h)$, $(\mathcal{D}_q)_i(\w_h)$ and $(\mathcal{R}_q)_i(\w_h)$ of
the cloud operators.

\subsubsection{Advection}
The advection terms are discretized using flux functions as follows:
\begin{align*}
(\mathcal{A}_L)_i(\w_h)&=\frac{1}{h}\sum_{j\in S(i)}H_{ij}^L(\w_h)\sum_{k=1}^dn_{ij}^{(k)},\\
(\mathcal{A}_N)_i(\w_h)&=\frac{1}{h}\sum_{j\in S(i)}H_{ij}^N(\w_h)\sum_{k=1}^dn_{ij}^{(k)},\\
(\mathcal{A}_q)_i((\w_q)_h)&=\frac{1}{h}\sum_{j\in S(i)}(H_q)_{ij}((\w_q)_h)\sum_{k=1}^dn_{ij}^{(k)},
\end{align*}
where the numerical fluxes $H_{ij}^L$, $H_{ij}^N$ and $(H_q)_{ij}$ approximate the corresponding fluxes between the computational cells
$C_i$ and $C_j$, and $n_{ij}^{(k)}$ denotes the $k$-th component of the outer normal unit vector of cell $C_i$ in the direction of cell
$C_j$. We use the Rusanov numerical flux for $H_{ij}^N$ and $(H_q)_{ij}$ and the central flux for $H_{ij}^{L}$. For
$(\mathcal{A}_N)_i(\w_h)$ and $(\mathcal{A}_q)_i((\w_q)_h)$ a discretization is obtained via a MUSCL-type approach using
piecewise linear reconstructions with the minmod limiter. {\cblue It is well-known this approach yields an approximation, which is almost
second-order accurate as its accuracy deteriorates at extrema.} The numerical fluxes are given by
\begin{align}
\label{fluxes}
H_{ij}^L(\w_h)&=\frac{1}{2}\(\F_L(\w_{j})+\F_L(\w_{i})\),\notag\\
H_{ij}^N(\w_h)&=\frac{1}{2}\(\F_N(\w_{ij}^+)+\F_N(\w_{ij}^-)\)-\frac{\lambda_{ij}}{2}\(\w_{ij}^+-\w_{ij}^-\),\\
(H_q)_{ij}((\w_q)_h)&=\frac{1}{2}\(\F_q((\w_q)_{ij}^+)+\F_q((\w_q)_{ij}^-)\)-\frac{(\lambda_q)_{ij}}{2}\((\w_q)_{ij}^+-(\w_q)_{ij}^-\).
\notag
\end{align}
Here, $\w_{ij}^-$, $\w_{ij}^+$ and $(\w_q)_{ij}^-$, $(\w_q)_{ij}^+$ denote the corresponding interface values, which are computed using a
piecewise linear reconstruction so that
\begin{equation*}
\w_{ij}^-=\w_i+\bm{s}_{ij}\frac{h}{2}\sum_{k=1}^dn_{ij}^{(k)},\quad\w_{ij}^+=\w_j-\bm{s}_{ji}\frac{h}{2}\sum_{k=1}^dn_{ij}^{(k)},
\end{equation*}
where the slopes $\bm{s}_{ij}$ are computed by the minmod limiter,
\begin{equation*}
\bm{s}_{ij}=\frac{1}{h}{\rm minmod}\(\w_j-\w_i,\w_i-\w_{j^*}\)\sum_{k=1}^dn_{ij}^{(k)},
\end{equation*}
applied in a component-wise manner. Here,
\begin{equation*}
{\rm minmod}(a,b)=\begin{cases}
a,~\mbox{if}~|a|<|b|~\mbox{and}~ab>0,\\
b,~\mbox{if}~|b|<|a|~\mbox{and}~ab>0,\\
0,~\mbox{if}~ab\le0,
\end{cases}
\end{equation*}
and $(\w_q)_{ij}^-$ and $(\w_q)_{ij}^+$ are obtained similarly. Thereby $C_{j*}$ is the other neighboring cell of $C_i$ in the opposite
direction from $C_j$. Finally, the values $\lambda_{ij}$ and $(\lambda_q)_{ij}$ are given by
\begin{equation*}
\lambda_{ij}=\max\left\{\sigma\(\frac{\partial\F_N(\w_{ij}^-)}{\partial\w}\),\,\sigma\(\frac{\partial\F_N(\w_{ij}^+)}{\partial\w}\)\right\},
\quad(\lambda_q)_{ij}=\max\left\{\sigma\(\frac{\partial\F_q((\w_q)_{ij}^-)}{\partial\w_q}\),\,
\sigma\(\frac{\partial\F_q((\w_q)_{ij}^+)}{\partial\w_q}\)\right\},
\end{equation*}
where $\sigma$ denotes the spectral radius of the corresponding Jacobians.
\begin{remark}
Note that in the computation of $H_{ij}^{L}$ in \eqref{fluxes}, we use the cell averages rather than the point values at the cell interfaces
for the following two reasons. First, the flux is second-order accurate. Second, in Section \ref{discr_time}, we will treat the linear part
of the flux implicitly and this is much easier to do when the numerical flux is linear as well.
\end{remark}

\subsubsection{Diffusion}
The components of the discrete diffusion operators are discretized in a straightforward manner using second-order central differences.

\subsubsection{Reaction}
The reaction terms are discretized by a direct evaluation of the reaction operators at the cell centers:
\begin{align*}
\mathcal{R}_i(\w_h)=\R_L(\w_i)+\R_N(\w_i),\quad(\mathcal{R}_q)_i((\w_q)_h)=\R_q((\w_q)_i).
\end{align*}

After the spatial discretization, we obtain the following system of time-dependent ODEs:
\begin{align}
\frac{{\rm d}}{{\rm d}t}\w_i&=-\mathcal{A}_i(\w_h)+\mathcal{D}_i(\w_h)+\mathcal{R}_i(\w_h),\label{ODE_NS}\\
\frac{{\rm d}}{{\rm d}t}(\w_q)_i&=-(\mathcal{A}_q)_i((\w_q)_h)+(\mathcal{D}_q)_i((\w_q)_h)+(\mathcal{R}_q)_i((\w_q)_h).\label{ODE_cloud}
\end{align}
This system has to be solved using an appropriate ODE solver as discussed in Section \ref{discr_time}.

\subsection{Discretization in time}\label{discr_time}
Let $\w_h^n$ and $(\w_q)_h^n$ denote the numerical approximation of the solutions $\w_h(t)$ and $(\w_q)_h(t)$ at the discrete time level
$t^n$. We evolve the solution to the next time level $t^{n+1}=t^n+\Delta t^n$, where $\Delta t^n$ is the size of the Strang operator
splitting time step. In the operator splitting approach, we first numerically solve the ODE system \eqref{ODE_NS} with
$\Delta t^n_{\rm NS}=\nicefrac{\Delta t^n}{2}$, we then numerically integrate the ODE system \eqref{ODE_cloud} with $\Delta t^n$ and
finally we solve the system \eqref{ODE_NS} again with $\Delta t^n_{\rm NS}$.

Notice that the system \eqref{ODE_NS} may be very stiff as the Navier-Stokes equations are in the weakly compressible regime. We therefore
follow the approach in \cite{bispen1} (see also \cite{diss}), and employ the second-order ARS(2,2,2) IMEX method from \cite{ARS222}:
\begin{equation}
\begin{aligned}
\w_h^{n+\frac{1}{4}}&=\w_h^n+\beta\Delta t^n_{\rm NS}\(\mathcal{L}\(\w_h^{n+\frac{1}{4}}\)+\mathcal{N}\(\w_h^n\)\),\\
\w_h^{n+\frac{1}{2}}&=\w_h^n+\Delta t^n_{\rm NS}\(\delta\mathcal{N}\(\w_h^n\)+(1-\delta)\mathcal{N}\(\w_h^{n+\frac{1}{4}}\)\)+
\Delta t^n_{\rm NS}\(\beta\mathcal{L}\(\w_h^{n+\frac{1}{2}}\)+(1-\beta)\mathcal{L}\(\w_h^{n+\frac{1}{4}}\)\),
\end{aligned}
\label{4.5}
\end{equation}
where $\alpha=1-\nicefrac{1}{\sqrt{2}}$, $\delta=1-\nicefrac{1}{2\beta}$, $t^{n+\frac{1}{2}}=t^n+\Delta t^n_{\rm NS}$,
$t^{n+\frac{1}{4}}=t^n+\nicefrac{\Delta t^n_{\rm NS}}{2}$, and $\Delta t^n_{\rm NS}$ satisfies the following CFL condition:
$$
\max_{s=1,2,3}\,\max_{i=1,\ldots,N}(|(u_s)_i|)\,\frac{\Delta t^n_{\rm NS}}{h}<0.5.
$$
For solving the linear systems arising in \eqref{4.5}, we use the generalized minimal residual (GMRES) method combined with a
preconditioner, the incomplete LU factorization (ILU). As it was shown in \cite{bispen1} (see also \cite{diss}), the resulting method is
both accurate and efficient in the weakly compressible regime.
	
The ODE system \eqref{ODE_cloud} is also stiff, but its stiffness only comes from the diffusion and power-law-type source terms. We
therefore efficiently solve it using the large stability domain third-order Runge-Kutta method from \cite{Dumka3}. We have utilized the ODE
solver DUMKA3, which is a free software that can be found in \cite{Dumka3_code}. We note that DUMKA3 selects time steps automatically, but
in order to improve its efficiency, one needs to provide the code with a time step stability restriction for the forward Euler method; see
\cite{Dumka3_code,Dumka3}. This bound is obtained by $\min\{\Delta t^n,\,\Delta t^n_{\rm cloud}\}$, where $\Delta t^n_{\rm cloud}$ satisfies
the following CFL condition for the cloud system:
$$
\max_{s=1,2}\,\max_{i=1,\ldots,N}(|(u_s)_i|,|(u_3)_i + v_q|)\,\frac{\Delta t^n_{\rm cloud}}{h}<0.5.
$$

\section{Deterministic numerical experiments}\label{sec6}
In this section, we test the numerical method described in Section \ref{sec4}. The experimental order of convergence is computed for the
so-called free convection of a moist warm air bubble and the structure formation in cloud dynamics is shown in the Rayleigh-B\'enard
convection. The latter will be simulated in both the 2-D and 3-D cases.

\subsection{Free convection of a moist warm air bubble in 2-D}\label{SWAB_det}
We start with the well-known meteorological benchmark describing the free convection of a smooth warm air bubble; see, e.g.,
\cite{BF02,DT50}.

\subsubsection*{Example 1} In this experiment, the warm bubble rises and deforms axisymmetrically due to the shear friction with the
surrounding air at the warm/cold air interface, gradually forming a mushroom-like shape. The warm air bubble is placed at $(3500\,\mathrm{m},
2000\,\mathrm{m})$ with the initial perturbation:
\begin{eqnarray*}
&&\rhoprime(\x,0)=-\rhobar(\x)\frac{\thetaprime(\x,0)}{\thetabar(\x)+\thetaprime(\x,0)},\quad\rhobar(\x)=
\frac{p_0}{R\thetabar(\x)}\pi_e(\x)^\frac{1}{\gamma-1},\quad\pi_e(\x)=1-\frac{g x_3}{c_p\thetabar},\\
&&\bm{u}(\x,0)=0,\\
&&\thetaprime(\x,0)=\begin{cases}
2\cos^2\left(\frac{\pi r}{2}\right),&r:=\sqrt{(x_1-3500)^2+(x_3-2000)^2}\le2000,\\
0,&\text{otherwise},
\end{cases}
\end{eqnarray*}
where $\thetabar=300\,\mathrm{K}$ and $p_0=\pbar=10^5\,\mathrm{Pa}$. The experiment was simulated in a domain $\Omega=[0,7000]\times[0,5000]\,\mathrm{m^2}$. For the cloud variables we choose the following
initial conditions:
{\cblue
$$
q_v(\x,0)=0.08\,\thetaprime(\x,0),\quad q_c=10^{-3}\,\thetaprime(\x,0),\quad q_r=10^{-5}\,\thetaprime(\x,0).
$$
}
Furthermore, we apply the no-flux boundary conditions $\uvec\cdot\bm{n}=0$, $\nabla\rhoprime\cdot\bm{n}=0$,
$\nabla\rhothetaprime\cdot\bm{n}=0$, $\nabla(\rho q_\ell)\cdot\bm{n}=0$, $\ell\in\{v,c,r\}$.

\begin{figure}[ht!]
\includegraphics[scale=0.95]{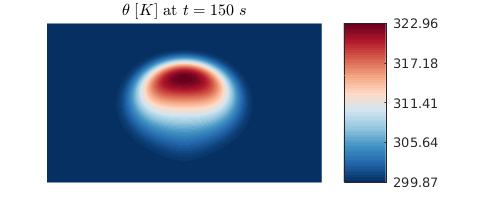}\includegraphics[scale=0.95]{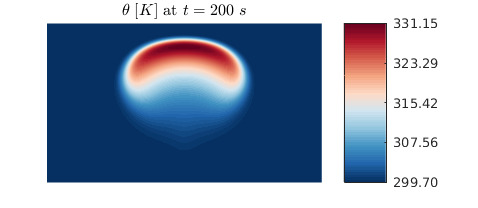}
\includegraphics[scale=0.95]{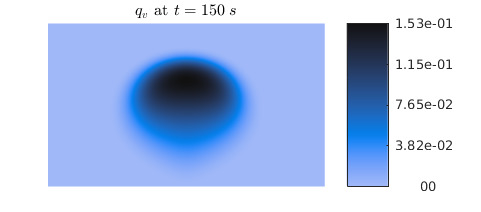}\includegraphics[scale=0.95]{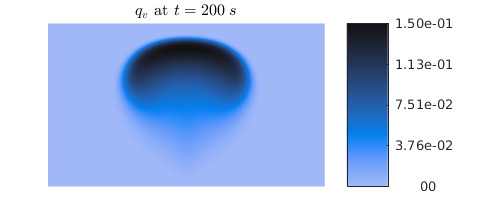}
\includegraphics[scale=0.95]{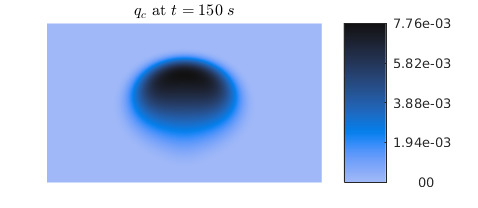}\includegraphics[scale=0.95]{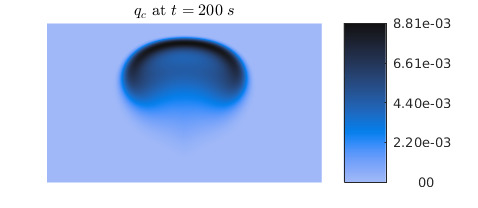}
\includegraphics[scale=0.95]{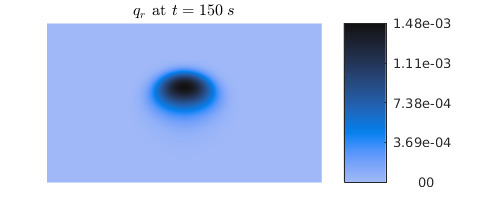}\includegraphics[scale=0.95]{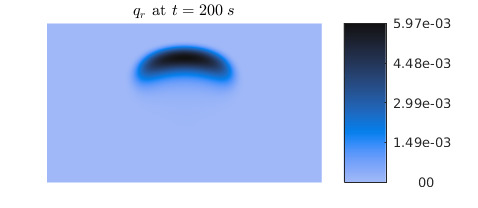}
\caption {{\cblue Example 1: Potential temperature $\theta$, water vapor concentration $q_v$, cloud drops concentration $q_c$ and rain
concentration $q_r$ at times $t=150$ (left column) and $200s$ (right column) simulated on a $160\times160$ mesh.\label{SWAB_T200}}}
\end{figure}

\begin{figure}[ht!]
\centerline{\hspace*{0.5cm}\includegraphics[scale=1.05]{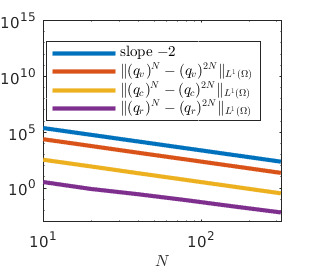}\includegraphics[scale=1.05]{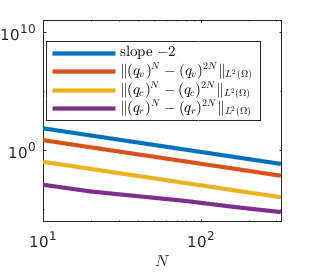}
\includegraphics[scale=1.05]{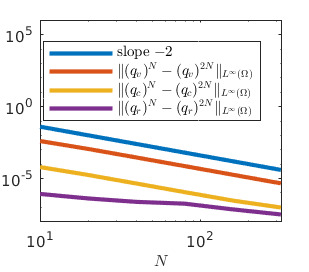}}
\caption{{\cblue Example~1: Comparison of different error norms for the cloud variables $q_v$, $q_c$ and $q_r$ computed at time $t=10s$
using $\Delta t=\nicefrac{20}{N}$.\label{cloud_L}}}
\end{figure}

\begin{figure}[ht!]
\centerline{\hspace*{0.5cm}\includegraphics[scale=1.05]{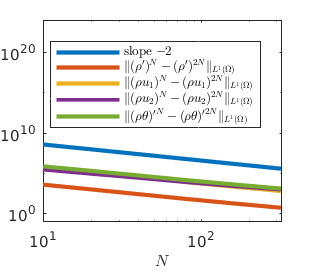}\includegraphics[scale=1.05]{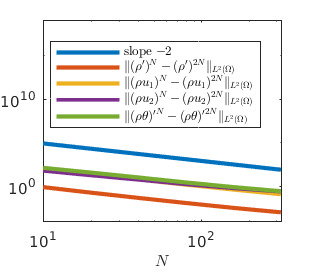}
\includegraphics[scale=1.05]{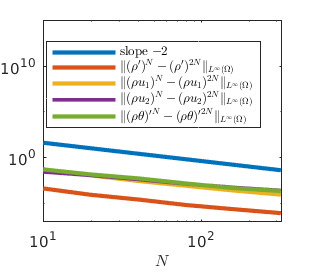}}
\caption{{\cblue Example 1: Comparison of different error norms for the flow variables $\rhoprime$, $\rho u_1$, $\rho u_2$ and
$\rhothetaprime$ at time $t=10s$ using $\Delta t=\nicefrac{20}{N}$.\label{NS_L}}}
\end{figure}

{\cblue In Figure \ref{SWAB_T200}, we show the potential temperature $\theta$ and cloud variables $q_v$, $q_c$ and $q_r$, computed on a
$160\times160$ mesh at $t=150$ and $200s$. One can clearly observe condensation taking place on the interface between cold and warm air}
{\cblue and leading to cloud formation in this region. In consequence, rain is formed in the clouds and falls towards the surface. Note
that the order of magnitude of the different water mass concentrations is very different, that is, $q_v\gg q_c\gg q_r$, as expected.} The
experimental convergence study for the cloud and flow variables is presented in Tables \ref{SWAB2D_cloudNS_T10} and \ref{SWAB2D_NS_T10},
respectively. The experimental order of convergence (EOC) has been computed in the following way:
{\cblue
$$
EOC=\log_2\left(\frac{\|v_{N,\Delta t}-v_{2N,\nicefrac{\Delta t}{2}}\|_{L^2(\Omega)}}
{\|v_{2N,\nicefrac{\Delta t}{2}}-v_{4N,\nicefrac{\Delta t}{4}}\|_{L^2(\Omega)}}\right),
$$}
where $v_{N,\Delta t}$ is the numerical solution computed on a grid with $N\times N$ grid cells and using a fixed time step $\Delta t$. As
one can clearly see, the expected second order of accuracy has been achieved. {\cblue For comparison, we present in Figures \ref{cloud_L}
and \ref{NS_L} the errors measured in the $L^1$-, $L^2$- and $L^\infty$-norms. They all give similar results.}

\begin{table}[ht!]
\begin{center}
\begin{tabular}{|c|cc|cc|cc|}
\hline
$N$&$L^2$-error in $\rho q_v$&EOC&$L^2$-error in $\rho q_c$&EOC&$L^2$-error in $\rho q_r$&EOC\\
\hline
10&6.870e+00&--~&1.019e-01&--~&1.159e-03&--~\\
20&1.711e+00&2.01&2.544e-02&2.00&3.152e-04&1.88~\\
40&4.271e-01&2.00&6.380e-03&2.00&1.240e-04&1.35~\\
80&1.080e-01&1.98&1.611e-03&1.99&4.952e-05&1.32~\\
160&2.703e-02&2.00&4.042e-04&1.99&1.571e-05&1.67~\\
320&6.765e-03&2.00&1.016e-04&1.99&5.666e-06&1.47~\\
\hline
\end{tabular}
\end{center}
\vspace{0.2cm}
\caption{{\cblue Example 1: $L^2$-errors and EOC for the cloud variables computed at time $t=10s$ using $\Delta t=\nicefrac{20}{N}$}.
\label{SWAB2D_cloudNS_T10}}
\end{table}
\begin{table}[ht!]
\begin{center}
\begin{tabular}{|c|cc|cc|cc|cc|}
\hline
$N$&$L^2$-error in $\rhoprime$&EOC&$L^2$-error in $\rho u_1$&EOC&$L^2$-error in $\rho u_2$&EOC&$L^2$-error in $\rhothetaprime$&EOC\\
\hline
10&7.522e-01&--~&1.134e+02&--~&5.805e+01&--~&1.213e+02&--~\\
20&1.757e-01&2.10&2.607e+01&2.12&1.744e+01&1.74&3.494e+01&1.80~\\
40&4.418e-02&2.00&5.604e+00&2.22&5.462e+00&1.67&9.641e+00&1.86~\\
80&1.147e-02&1.95&1.436e+00&1.96&1.658e+00&1.72&2.584e+00&1.90~\\
160&3.170e-03&1.85&3.972e-01&1.85&5.875e-01&1.50&7.577e-01&1.77~\\
320&9.810e-04&1.70&1.159e-01&1.78&2.420e-01&1.28&2.556e-01&1.57~\\
\hline
\end{tabular}
\end{center}
\vspace{0.2cm}
\caption{{\cblue Example 1: $L^2$-errors and EOC for the flow variables computed at time $t=10s$ using $\Delta t=\nicefrac{20}{N}$.
\label{SWAB2D_NS_T10}}}
\end{table}

\subsection{Rayleigh-B\'{e}nard convection}\label{RB_det}
In this experiment, we study a natural convection that is used to model structure formation. It occurs in a planar flow between two
horizontal plates, where the lower one is heated from below and the upper one is cooled from above. Due to the presence of buoyancy, and
hence gravity, the fluid develops a regular pattern of convection roles, known as the B\'enard cells. In 3-D, these convection roles form
additionally hexagonal structures; see, e.g., \cite{RB,pauluis_schumacher_2010,AR17}.

For our numerical simulations, we prescribe the following initial conditions:
\begin{align*}
&\rhoprime(\x,0)=-\rhobar(\x)\frac{\thetaprime(\x,0)}{\thetabar(\x)+\thetaprime(\x,0)},\quad
\rhobar(\x)=\frac{p_0}{R\thetabar(\x)}\pi_e(\x)^\frac{1}{\gamma-1},\quad\pi_e(\x)=1-\frac{g x_3}{c_p\thetabar},\\
&\bm{u}(\x,0)=0,\quad\thetaprime(\x,0)=\eta(\x),\quad\thetabar(\x)=284-\frac{1}{1000}x_3,
\end{align*}
where $p_0=\pbar=10^5\,\mathrm{Pa}$ and $\eta(\x)$ is a random perturbation uniformly distributed in $[-0.0021, 0.0021]$. For the cloud equations,
the following initial data are used:
\begin{eqnarray*}
q_v(\x,0)=2\cdot10^{-5}\thetabar(\x),\quad q_c=0,\quad q_r=0.
\end{eqnarray*}

We apply periodic boundary conditions in horizontal direction and the following conditions vertically: $\uvec\cdot\bm{n}=0$,
$\nabla\rhoprime\cdot\bm{n}=0$, $\nabla(\rho q_\ell)\cdot\bm{n}=0$, $\ell\in\{v,c,r\}$ with the Dirichlet boundary conditions for the
potential temperature,
\begin{equation*}
\theta(x_3=0)=284\,\mathrm{K}\quad\mbox{and}\quad\theta(x_3=1000)=283\,\mathrm{K}.
\end{equation*}

\subsubsection*{Example 2: 2-D case} 
\begin{figure}[ht!]
\centerline{\includegraphics[scale=1]{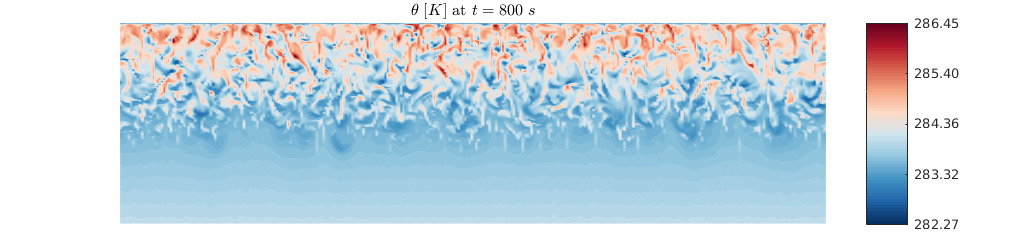}}
\centerline{\includegraphics[scale=1]{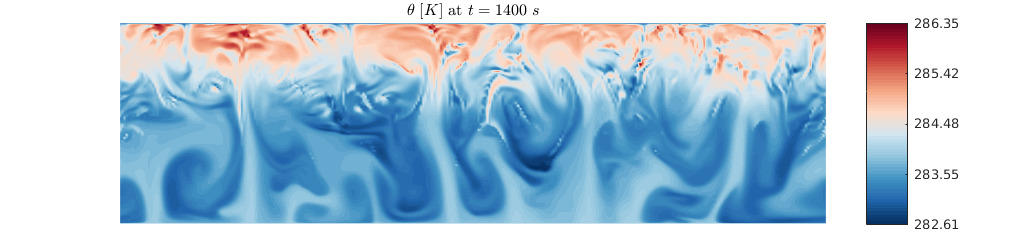}}
\caption{Example 2: Time evolution of the potential temperature $\theta$.\label{fig1}}
\end{figure}

\begin{figure}[ht!]
\centerline{\includegraphics[scale=1]{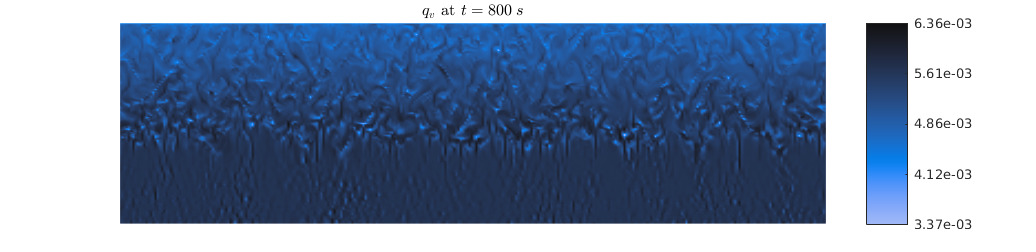}}
\centerline{\includegraphics[scale=1]{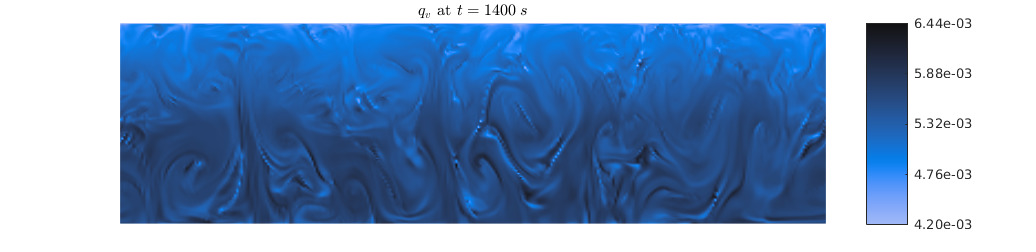}}
\caption{Example 2: Time evolution of the water vapor concentration $q_v$.\label{fig2}}
\end{figure}

\begin{figure}[ht!]
\centerline{\includegraphics[scale=1]{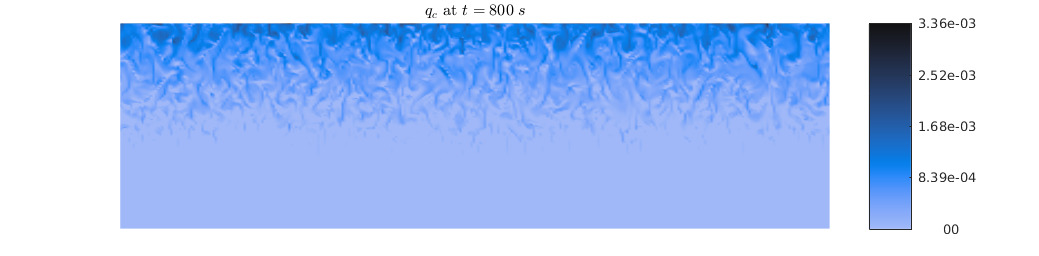}}
\centerline{\includegraphics[scale=1]{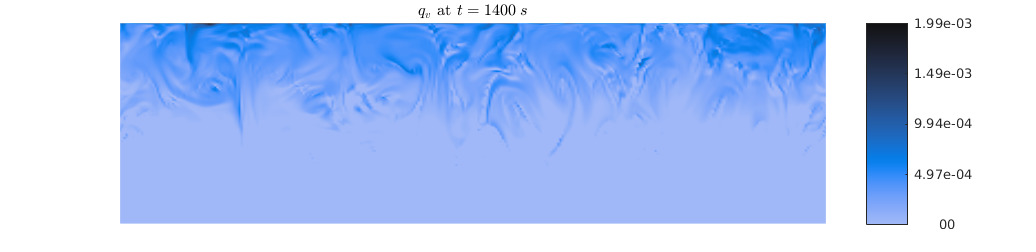}}
\caption{Example 2: Time evolution of the cloud drops concentration $q_c$.\label{fig3}}
\end{figure}

\begin{figure}[ht!]
\centerline{\includegraphics[scale=1]{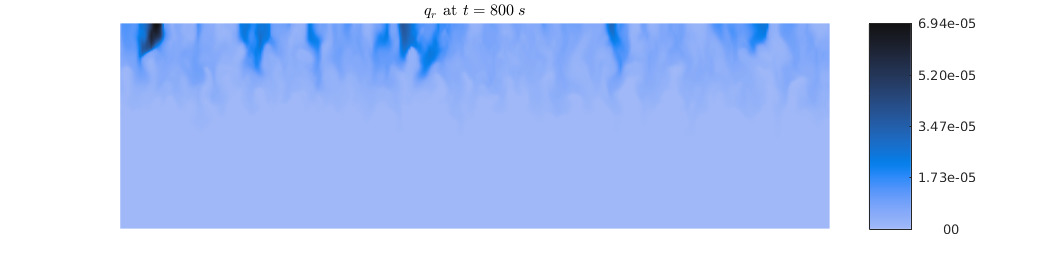}}
\centerline{\includegraphics[scale=1]{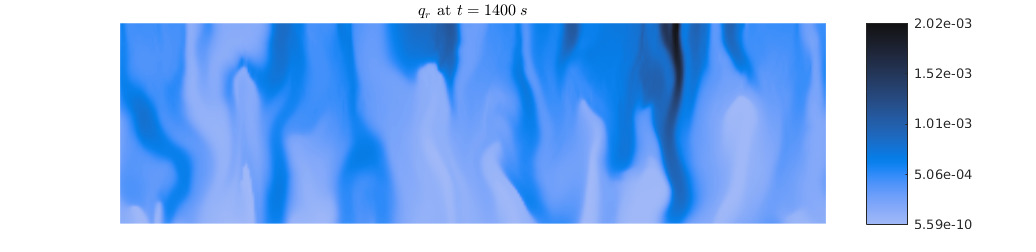}}
\caption{Example 2: Time evolution of the rain concentration $q_r$.\label{fig4}}
\end{figure}

In Figures \ref{fig1}--\ref{fig4}, we present time snapshots of the numerical solution computed in a
domain $\Omega=[0,5000]\times[0,1000]\,\mathrm{m^2}$ that has been discretized using $320\times320$ mesh cells. {\cblue The potential temperature,
water vapor mixing ratio, cloud mass and rain mass concentration are presented at two distinct times ($t=800$ and $1400s$) in Figures
\ref{fig1}--\ref{fig4}, where we can clearly see the formation and evolution of a pattern. At an earlier time $t=800\,s$, one can observe
the formation of small convective cells, visible as narrow finger-like structures reaching towards the top of the domain. Inside these
cells, the potential temperature is enhanced, partly due to the upward transport of higher values from below and partly due to phase changes
and thus latent heat release. Also the mass concentrations of water vapor and cloud water follow the small scale structure and show enhanced
values inside the fingers. Even at this early stage, rain can be formed at the top of the domain, since there the cloud water concentration
is high enough for autoconversion and accretion. Nevertheless, the structure of rain water is very different since after the formation of
rain, it is vertically transported due to sedimentation leading to vertically smeared structures. By a later time $t=1400\,s$, much larger
structures, which are similar to classical structures for dry thermal convection, have been formed. In the variables $\theta$, $q_v$ and
$q_c$, the spatially extended convective cells can be clearly seen. In contrast, rain water is not following the convective structure
although some larger features can be seen. In general, smearing due to sedimentation is again a major feature of the rain mass
concentration.}


\subsubsection*{Example 3: 3-D case}
\begin{figure}[ht!]
\includegraphics[scale=0.18]{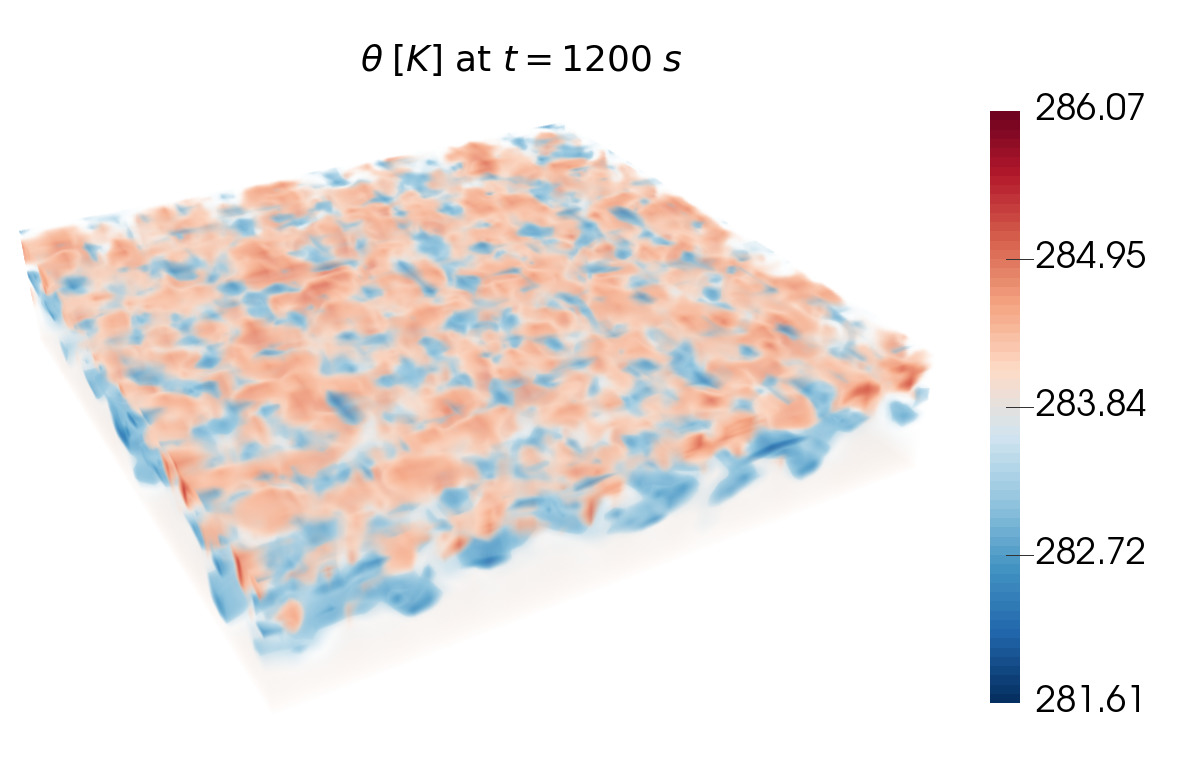}\hspace*{0.5cm}\includegraphics[scale=0.18]{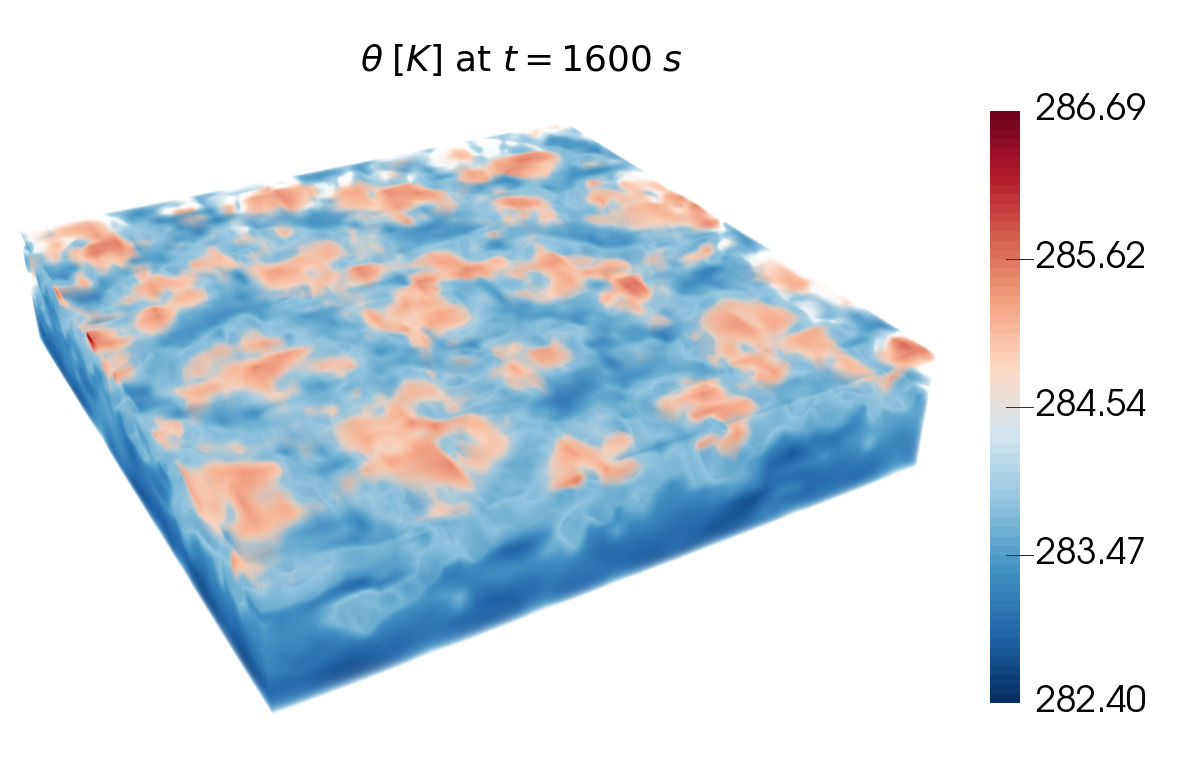}
\includegraphics[scale=0.18]{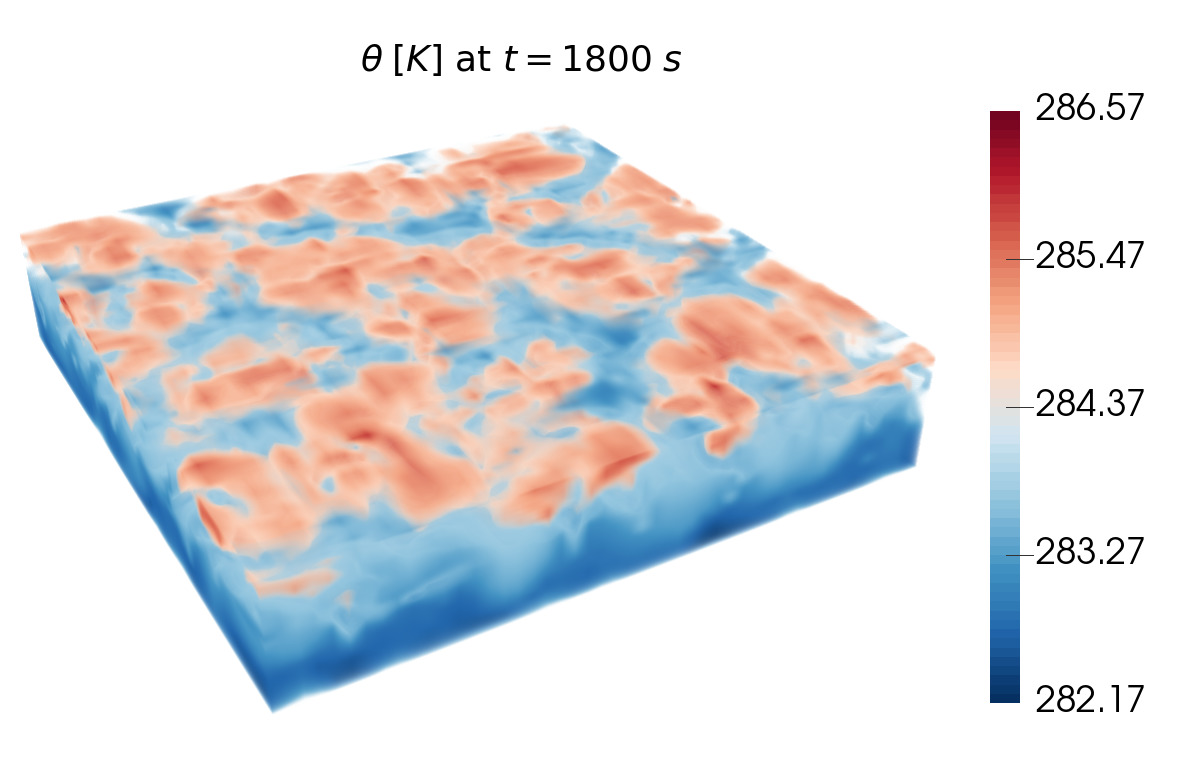}\hspace*{0.5cm}\includegraphics[scale=0.18]{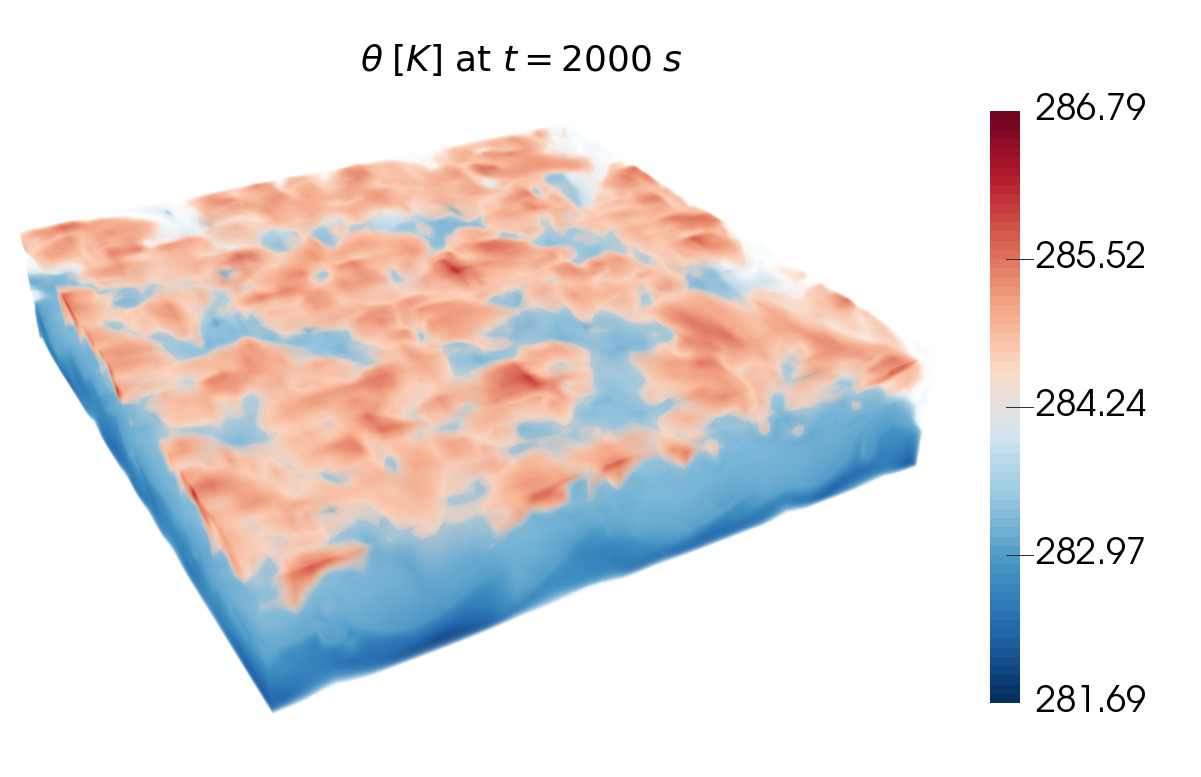}
\caption {Example 3: Time evolution of the potential temperature $\theta$.\label{fig5}}
\end{figure}

\begin{figure}[ht!]
\includegraphics[scale=0.18]{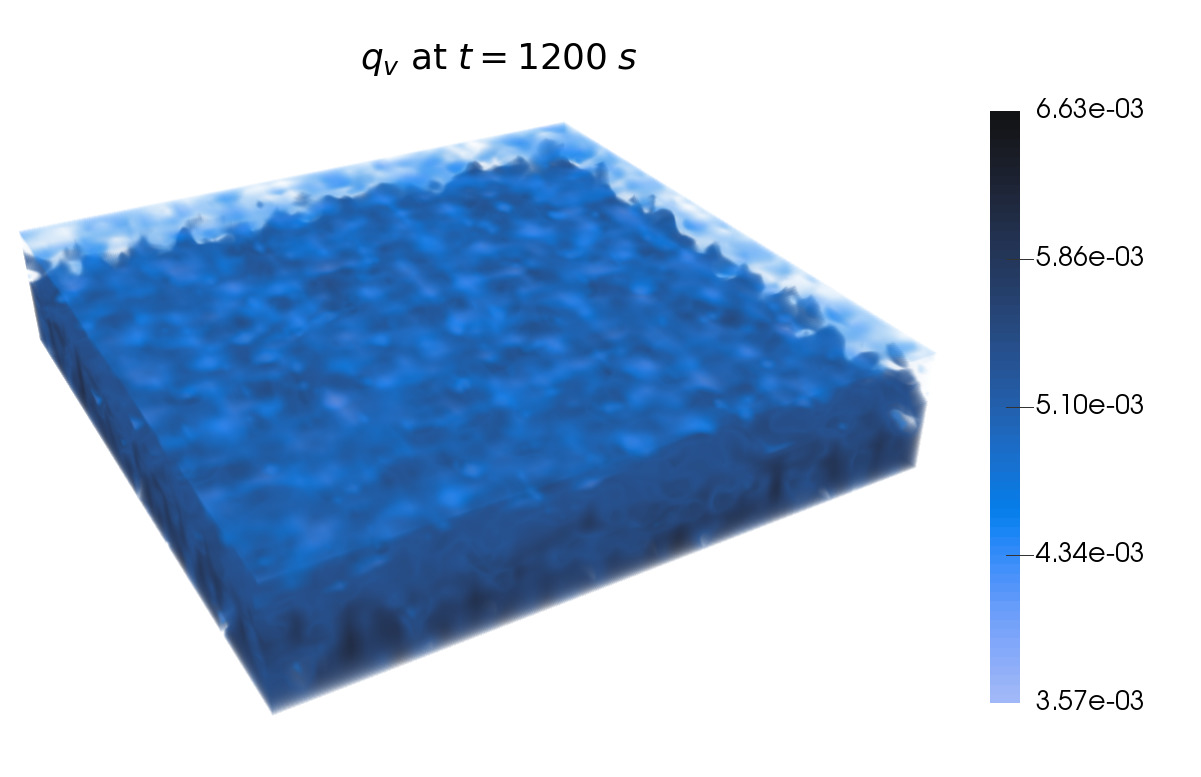}\hspace*{0.5cm}\includegraphics[scale=0.18]{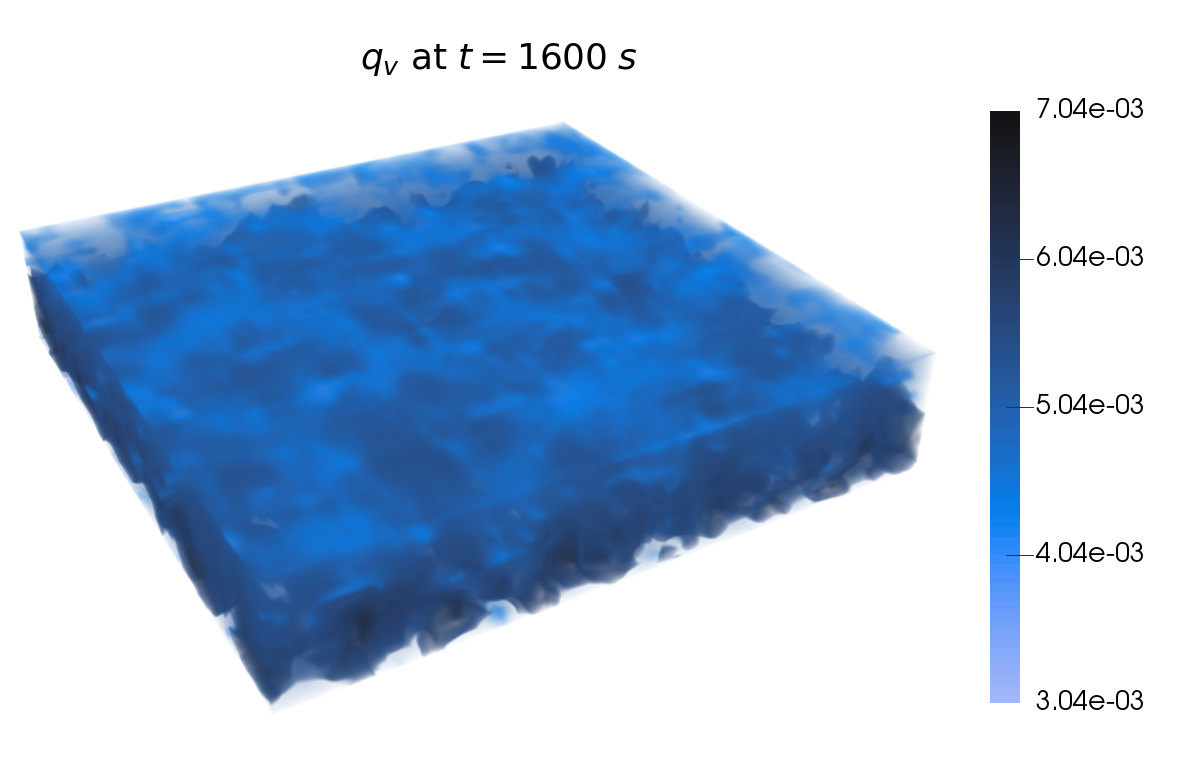}
\includegraphics[scale=0.18]{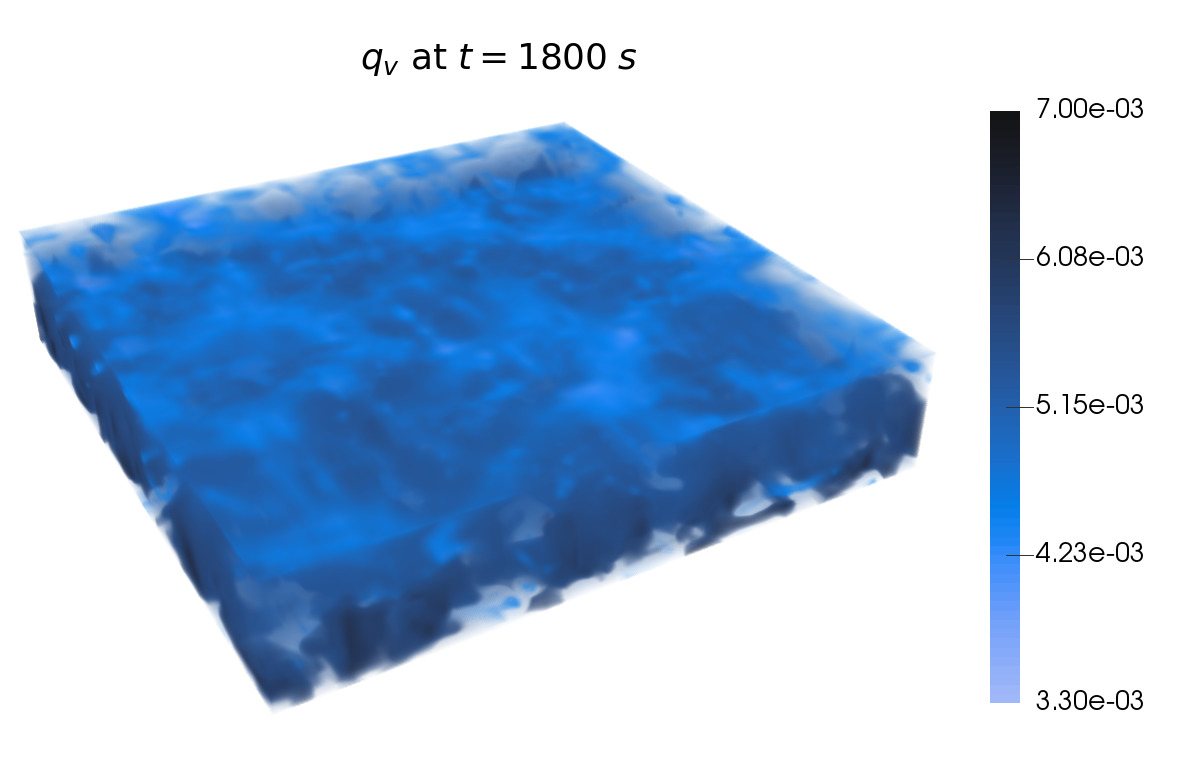}\hspace*{0.5cm}\includegraphics[scale=0.18]{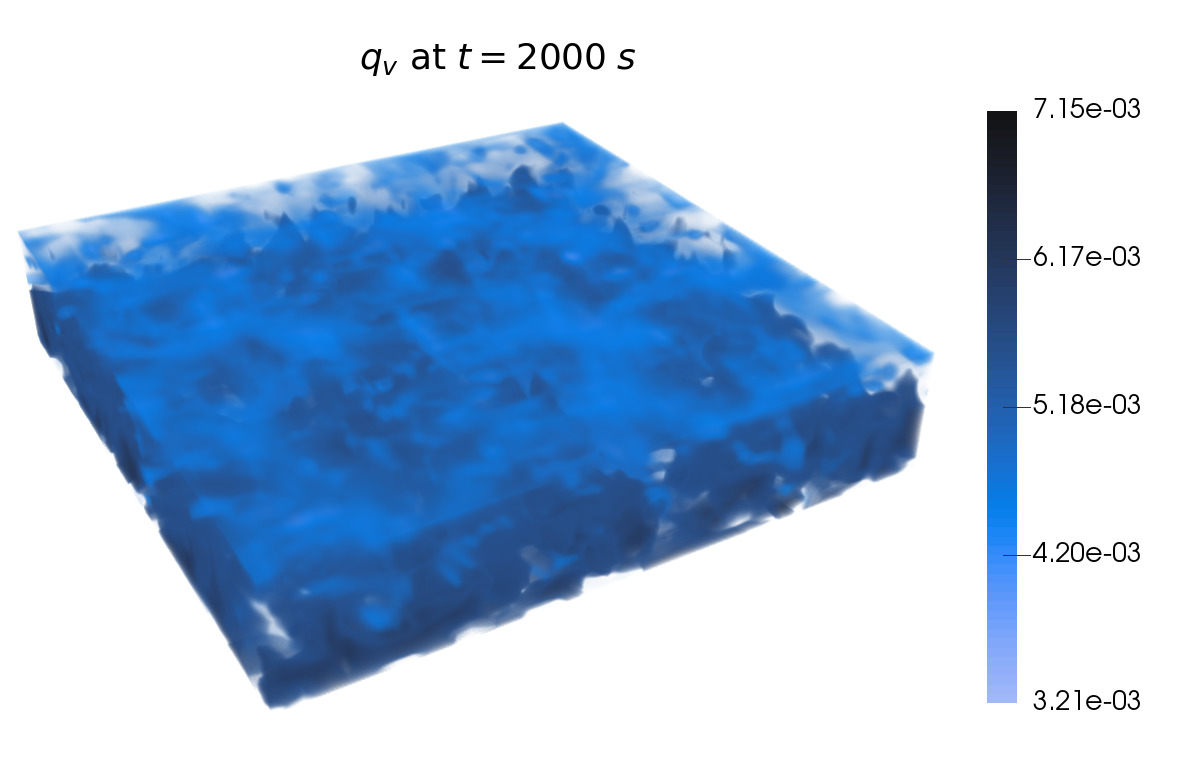}
\caption {Example 3: Time evolution of the water vapor concentration $q_v$.\label{fig6}}
\end{figure}

\begin{figure}[ht!]
\includegraphics[scale=0.18]{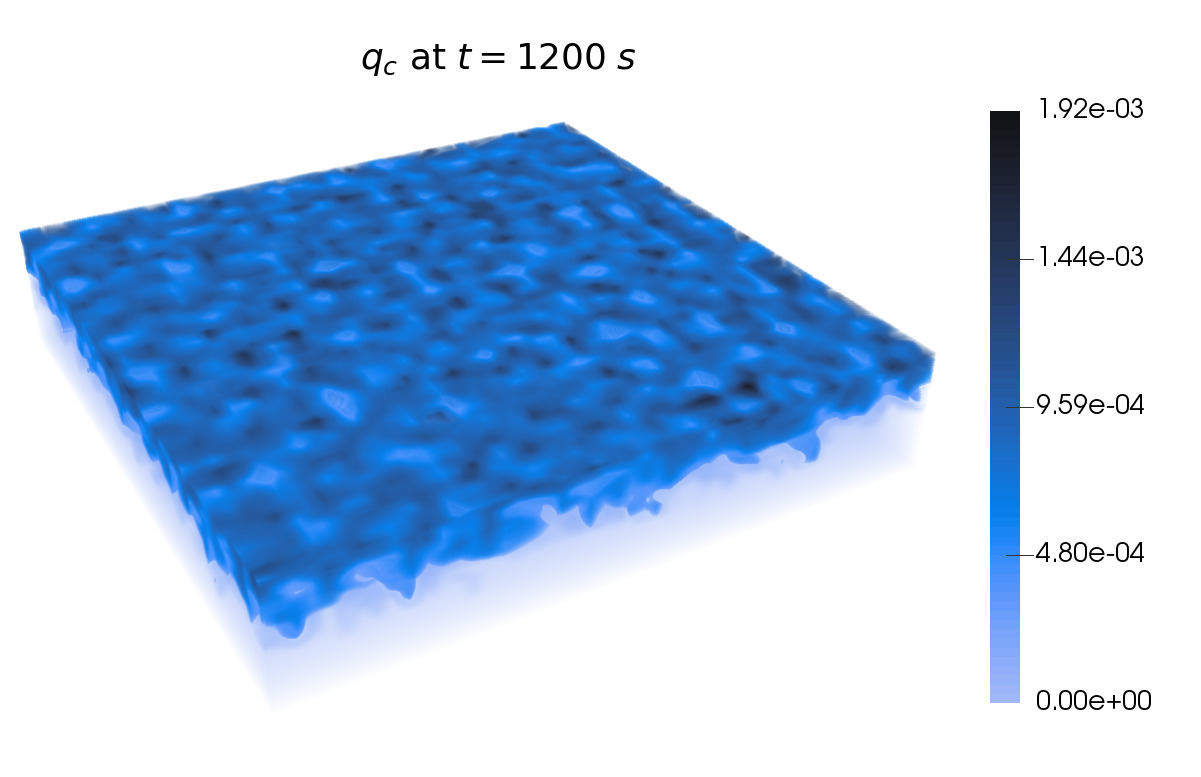}\hspace*{0.5cm}\includegraphics[scale=0.18]{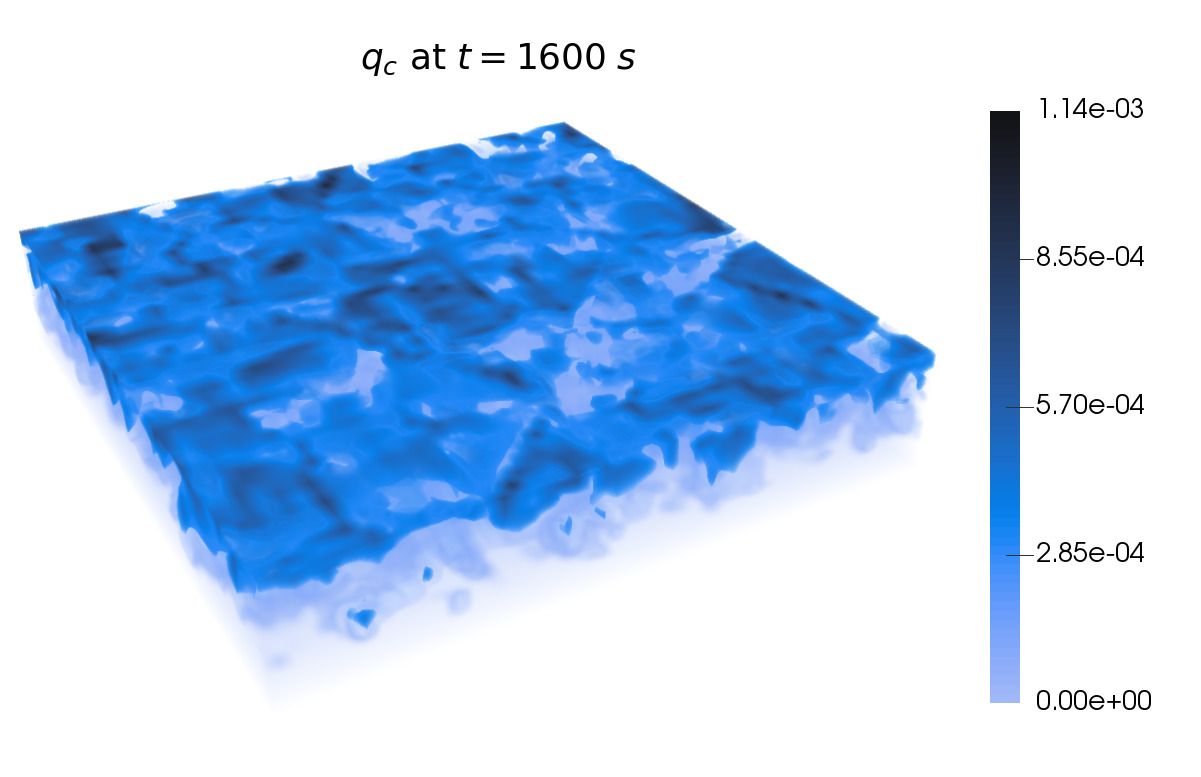}
\includegraphics[scale=0.18]{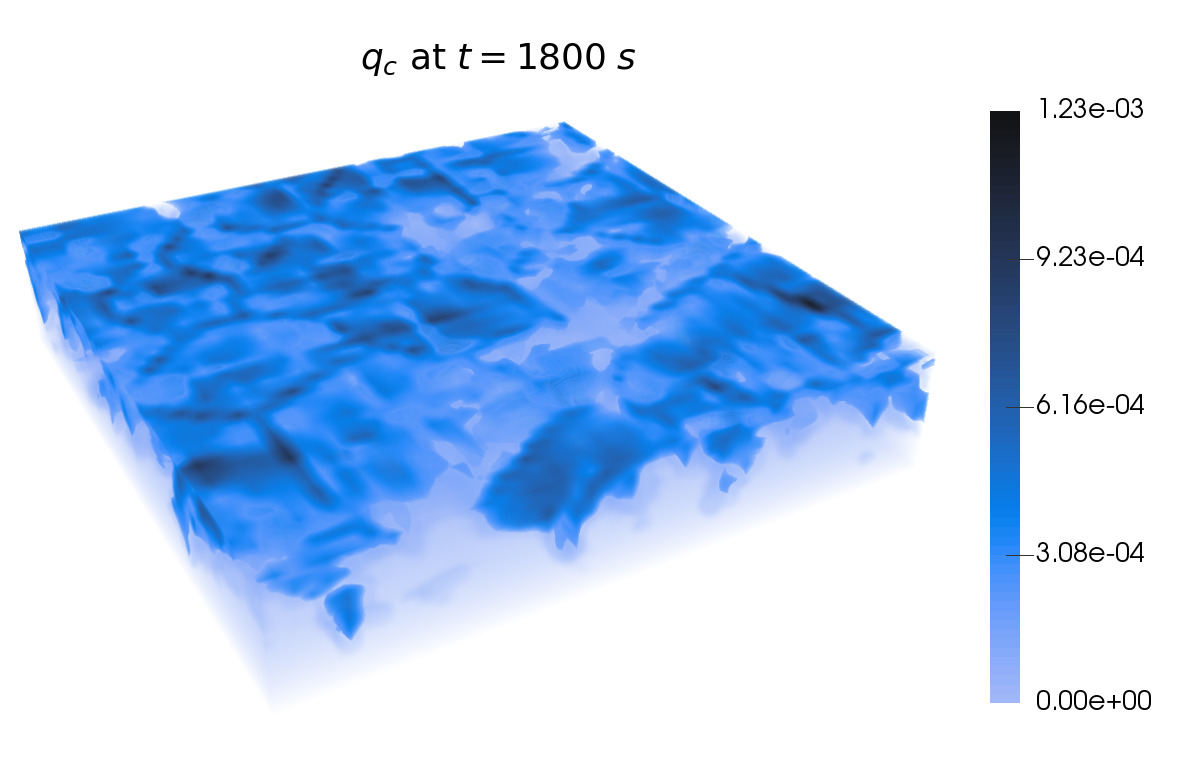}\hspace*{0.5cm}\includegraphics[scale=0.18]{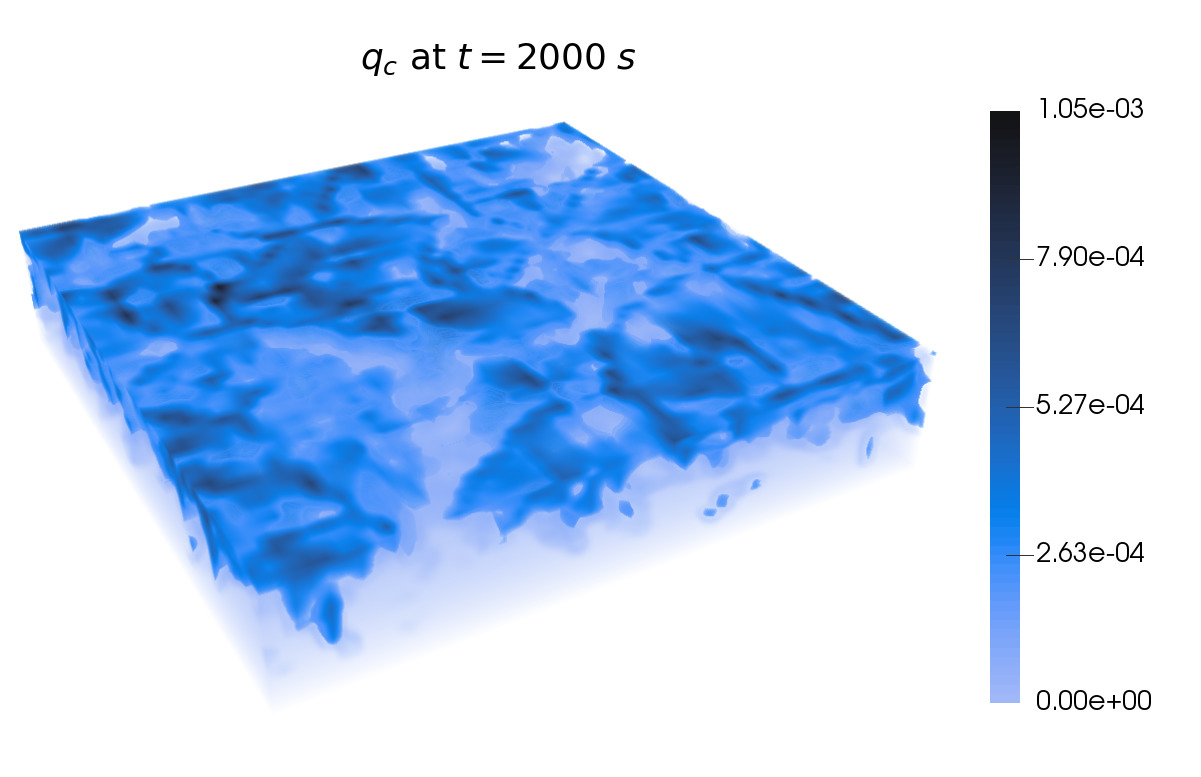}
\caption {Example 3: Time evolution of the cloud drops concentration $q_c$.\label{fig7}}
\end{figure}

\begin{figure}[ht!]
\includegraphics[scale=0.18]{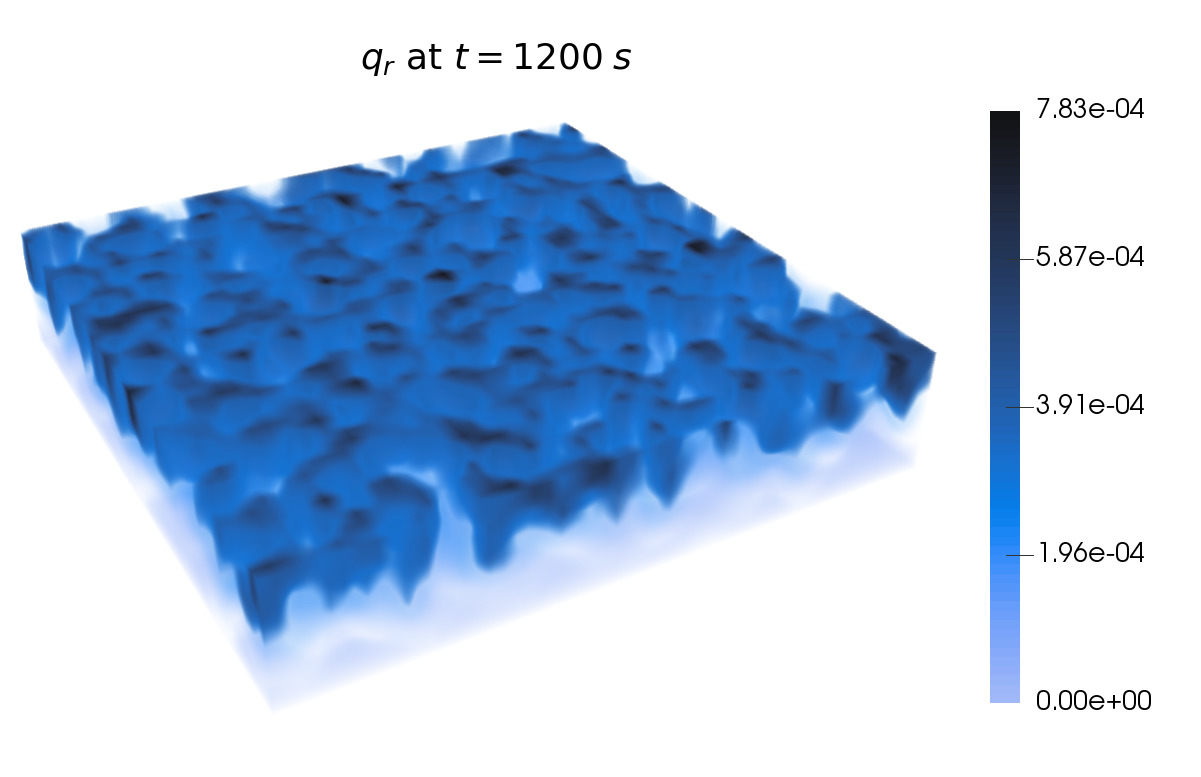}\hspace*{0.5cm}\includegraphics[scale=0.18]{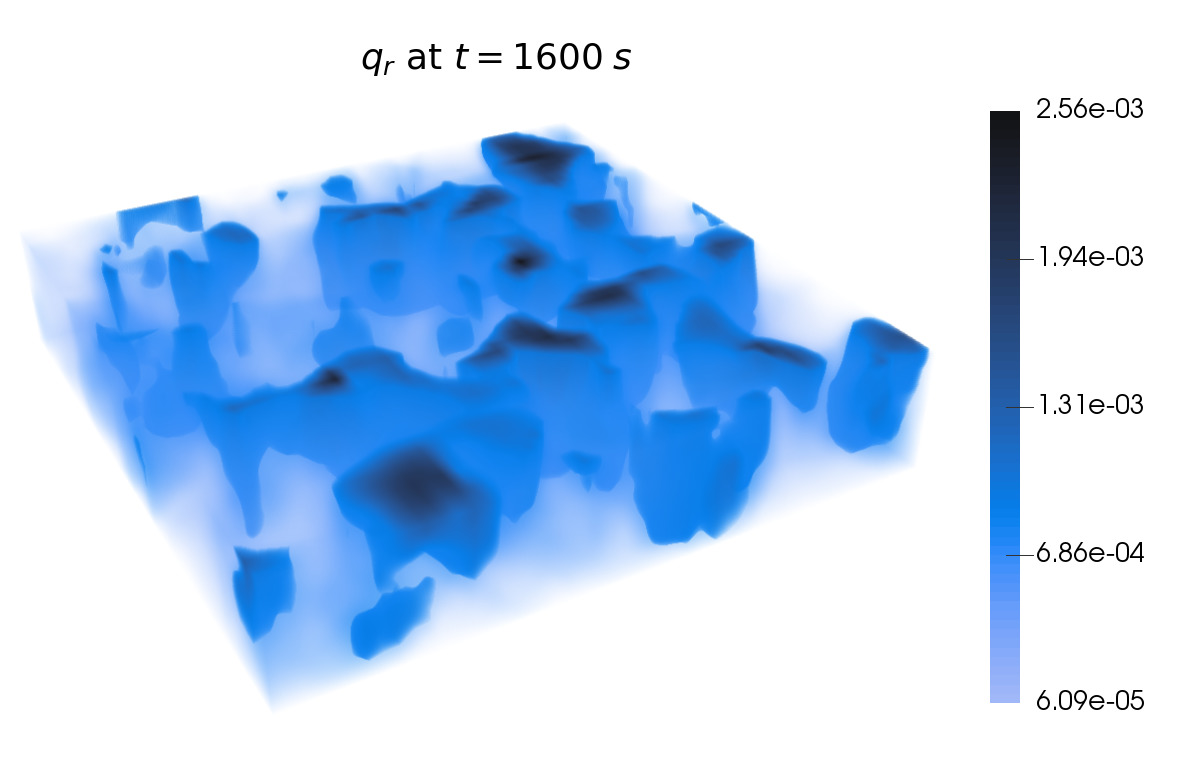}
\includegraphics[scale=0.18]{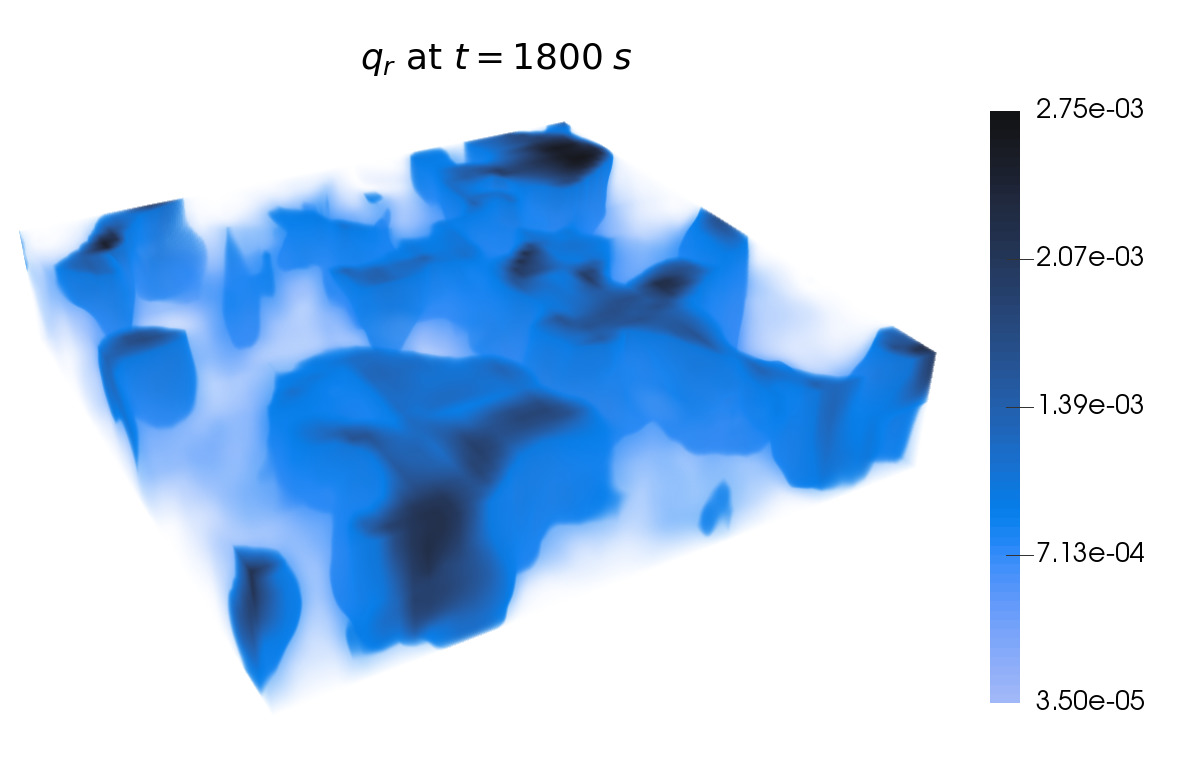}\hspace*{0.5cm}\includegraphics[scale=0.18]{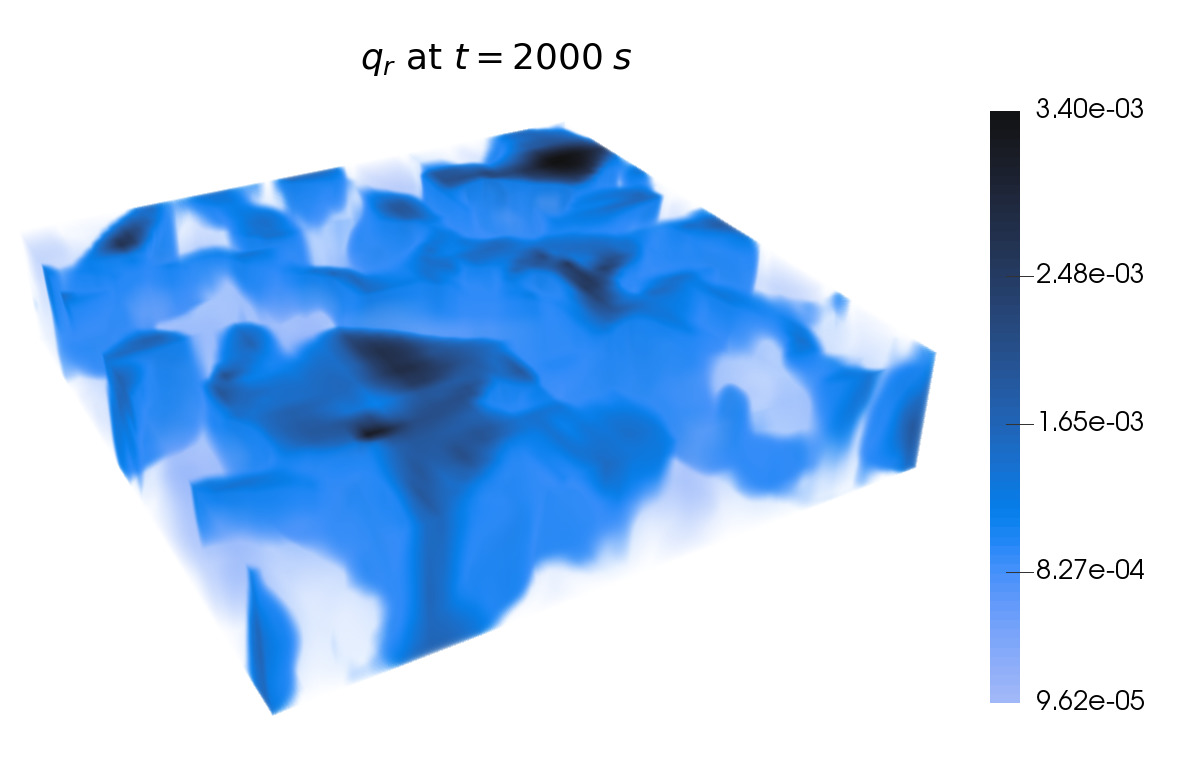}
\caption {Example 3: Time evolution of the rain concentration $q_r$.\label{fig8}}
\end{figure}

In this example, we compute the numerical solution in a domain $\Omega=[0,5000]\times[0,5000]\times[0,1000]\,\mathrm{m^3}$ that has been
discretized using $50\times50\times50$ mesh cells. {\cblue Figures \ref{fig5}--\ref{fig8} show $\theta$, $q_v$, $q_c$ and $q_r$ computed at
times $t=1200$, $1600$, $1800$ and $2000s$. In order to better visualize the computed structures, we have plotted the solution in a slightly
smaller domain $[0,5000]\times[0,5000]\times[0,980]\,\mathrm{m^3}$.}
{\cblue As in the 2-D case, different structures are formed in the different variables. At an earlier time $t=1200\,s$, small scale
structures can be seen at the top layer of the domain, especially in the potential temperature (Figure \ref{fig5}) and cloud water (Figure
\ref{fig7}). As time progresses, these structures aggregate and reorganize to quasi-hexagonal structures, which would be typical for
classical dry thermal convection---Rayleigh-B\'{e}nard convection. These structures in potential temperature seem to be quite robust as the
overall scales and pattern do not change from $t=1800\,s$ to $t=2000\,s$. As one can see in Figure \ref{fig7}, the structures in the cloud
water are very similar since the variables $\theta$ and $q_c$ are closely connected; the structures are mainly visible in horizontal planes.
For the rain distribution, the evolved structures are quite different since rain is formed at regions with high cloud water (that is, at the
top layers) and then transported by sedimentation leading to more pronounced pattern in the vertical direction, since sedimentation is the
dominant process after the rain has formed.}


{\cblue \begin{remark}
Let us note that the Rayleigh-B\'enard convection can be understood as a very simplified model for atmospheric convection in the
turbulent planetary boundary layer. In \cite{pauluis_schumacher_2010,WS}, numerical simulations for moist Rayleigh-B\'enard convection have
been realized using the Boussinesq approximation, a simplified equation of state, and the rigid-lid boundary conditions at the top and
bottom of the computational domain. Our mathematical model is more general and  takes weakly compressible effects into account. Numerical
experiments presented in Examples 2 and 3 are in good agreement with the results presented in the literature, but the focus of those studies
differs from ours.
\end{remark}
}

\section{Stochastic mathematical model}\label{sec3}
{\cred In meteorological applications, it is typical that initial or boundary data as well as parameters are uncertain. In order to take
this into account and analyze the influence of such uncertainties on the output of numerical simulations, we apply the stochastic Galerkin
method that will be described below. In this paper, we consider the case} where the uncertainty arises from the initial data or some
coefficients in the microphysical cloud parametrizations. In order to mathematically describe the uncertainty, we introduce a random
variable $\omega$. We assume that either the initial data or some well-chosen model parameters depend on $\omega$, that is,
\begin{equation*}
(\rho q_\ell)\big|_{t=0}=(\rho q_\ell)(\x,t=0,\omega)\quad\mbox{with}\quad\ell\in\{v,c,r\}
\end{equation*}
or
\begin{equation*}
k_1=k_1(\omega),\quad k_2=k_2(\omega),\quad\alpha=\alpha(\omega).
\end{equation*}
Consequently, the solution at later time will also depend on $\omega$, that is, $(\rho q_\ell)(\x,t,\omega)$ for $\ell\in\{v,c,r\}$, and the
system \eqref{cloud_equations} for cloud variables will read as
\begin{align}
\label{cloud_equations_random}
((\rho q_v)(\omega))_t+\nabla\cdot\((\rho q_v)(\omega)\bm{u}-\mu_q\rho\nabla q_v(\omega)\)&=\rho(-C(\omega)+E(\omega)),\notag\\
((\rho q_c)(\omega))_t+\nabla\cdot\((\rho q_c)(\omega)\bm{u}-\mu_q\rho\nabla q_c(\omega)\)&=\rho(C(\omega)-A_1(\omega)-A_2(\omega)),\\
((\rho q_r)(\omega))_t+\nabla\cdot\((\rho q_r)(\omega)(-v_q(\omega)\bm{e_3}+\bm{u})-\mu_q\rho\nabla q_r(\omega)\)&=
\rho(A_1(\omega)+A_2(\omega)-E(\omega)).\notag
\end{align}
From now on we will stress the dependence on $\omega$, but we will omit the dependence on $\x$ and $t$ to simplify the notation. We would
like to point out that the solution of the Navier-Stokes equations \eqref{NS_equations_pert} will also depend on $\omega$, because of the
source term $S_\theta$. In this paper, we will consider a simplified situation by replacing
\begin{equation*}
S_\theta(\omega)=\rho\frac{L\theta}{c_pT}\big\{C((\rho q_v)(\omega),(\rho q_c)(\omega))-E((\rho q_v)(\omega),(\rho q_r)(\omega))\big\}
\end{equation*}
in \eqref{NS_equations_pert} by $\bar{S}_\theta$ which only depends on the expected values of the cloud variables
\begin{equation*}
\bar{S}_\theta:=
\rho\frac{L\theta}{c_pT}\big\{C(\mathbb{E}[\rho q_v],\mathbb{E}[\rho q_c])-E(\mathbb{E}[\rho q_v],\mathbb{E}[\rho q_r])\big\}.
\end{equation*}
This ensures that all of the fluid variables, $\rhoprime$, $\rho\uvec$ and $\rhothetaprime$, remain deterministic.
	
\section{Numerical scheme for the stochastic model}\label{sec5}
In this section, we describe a generalized polynomial chaos stochastic Galerkin (gPC-SG) method for the system of cloud equations
\eqref{cloud_equations_random}. Such method belongs to the class of intrusive methods and the use of the Galerkin expansion leads to a
system of deterministic equations for the expansion coefficients. In the gPC-SG method, the solution is sought in the form of a polynomial
expansion
\begin{equation}
\rho q_\ell(\x,t,\omega)=\sum_{k=0}^M(\widehat{\rho q_\ell})_k(\x,t)\Phi_k(\omega)\quad\mbox{with}\quad\ell\in\{v,c,r\},~M\ge0,	
\label{expansion_rhoq}
\end{equation}
where $\Phi_k(\omega)$, $k=0,\ldots,M$, are  polynomials of $k$-th degree that are orthogonal with respect to the probability density
function $\mu(\omega)$. The choice of the orthogonal polynomials $\{\Phi_k(\omega)\}_{k=0}^M$ depends on the distribution of $\omega$. In
our case, we use a uniformly distributed $\omega\in\Gamma=(-1,1)$, which defines the Legendre polynomials, which satisfy
\begin{equation}
\int\limits_\Gamma\Phi_k(\omega)\Phi_{k'}(\omega)\mu(\omega)\,{\rm d}\omega = {\cred \frac{1}{2 k + 1}}\delta_{kk'}\quad\mbox{for}\quad0\le k,k'\le M,
\label{orthogonal_property}
\end{equation}
where $\delta_{kk'}$ is the Kronecker symbol and $\Gamma$ is the sample space.

We use the same expansion for the uncertain coefficients,
\begin{equation}
k_1(\omega)=\sum_{k=0}^M(\widehat{k_1})_k\Phi_k(\omega),\quad k_2(\omega)=\sum_{k=0}^M(\widehat{k_2})_k\Phi_k(\omega),\quad
\alpha(\omega)=\sum_{k=0}^M\widehat{\alpha}_k\Phi_k(\omega),
\label{expansion_coef}
\end{equation}
for the source terms on the RHS of \eqref{cloud_equations_random},
\begin{align}
\label{5.4}
\rho\(-C(\x,t,\omega)+E(\x,t,\omega)\)=:R_1(\x,t,\omega)&=\sum_{k=0}^M(\widehat{r_1})_k(\x,t)\Phi_k(\omega),\notag\\
\rho\(C(\x,t,\omega)-A_1(\x,t,\omega)-A_2(\x,t,\omega)\)=:R_2(\x,t,\omega)&=\sum_{k=0}^M(\widehat{r_2})_k(\x,t)\Phi_k(\omega),\\
\rho\(A_1(\x,t,\omega)+A_2(\x,t,\omega)-E(\x,t,\omega)\)=:R_3(\x,t,\omega)&=\sum_{k=0}^M(\widehat{r_3})_k(\x,t)\Phi_k(\omega),\notag
\end{align}
as well as for the raindrop fall velocity,
\begin{equation}
v_q(\x,t,\omega)=\sum_{k=0}^M(\widehat{v_q})_k(\x,t)\Phi_k(\omega).
\label{5.5}
\end{equation}
Since $\rho(\x,t)=\widehat{\rho}_0(\x,t)$, we also obtain
\begin{equation}
q_\ell(\x,t,\omega)=\sum_{k=0}^M(\widehat{q_\ell})_k(\x,t)\Phi_k(\omega)\quad\mbox{with}\quad(\widehat{q_\ell})_k(\x,t)=
\frac{(\widehat{\rho q_\ell})_k(\x,t)}{\rho(\x,t)}~~~\mbox{for}~\ell\in\{v,c,r\},~k=1,\ldots,M.
\label{expansion_q}
\end{equation}
We note that if $\rho(\x,t)$ is very small, the computation of the coefficients $(\widehat{q_\ell})_k(\x,t)$ should be desingularized;
see \cite[formulae (5.16)--(5.18)]{Kur_Acta}.

Applying the Galerkin projection to \eqref{cloud_equations_random} yields
\begin{align}
\label{cloud_equations_galerkin_projection}
\left\langle(\rho q_v)_t+\nabla\cdot\(\rho q_v\bm{u}-\mu_q\rho\nabla q_v\),\Phi_k\right\rangle&=
\left\langle\rho(-C+E),\Phi_k\right\rangle,\notag\\
\left\langle(\rho q_c)_t+\nabla\cdot\(\rho q_c\bm{u}-\mu_q\rho\nabla q_c\),\Phi_k\right\rangle&=
\left\langle\rho(C-A_1-A_2),\Phi_k\right\rangle,\\
\left\langle(\rho q_r)_t+\nabla\cdot\(\rho q_r(-v_q\bm{e_3}+\bm{u})-\mu_q\rho\nabla q_r\),\Phi_k\right\rangle&=
\left\langle\rho(A_1+A_2-E),\Phi_k\right\rangle,\notag
\end{align}
for $k=0,\ldots,M$, where $\langle\cdot,\cdot\rangle$ is the scalar product in our probability space which is given through
\begin{equation*}
\langle u,v\rangle=\int\limits_{-1}^1u(\omega)v(\omega)\mu(\omega)\,{\rm d}\omega.
\end{equation*}
We now substitute \eqref{expansion_rhoq}, \eqref{5.4}--\eqref{expansion_q} into \eqref{cloud_equations_galerkin_projection} and use the
orthogonality property \eqref{orthogonal_property} to obtain the following $3(M+1)\times3(M+1)$ deterministic system for the gPC
coefficients:
\begin{align}
\label{cloud_equations_coefficients}
\frac{\partial}{\partial t}(\widehat{\rho q_v})_k+\sum_{s=1}^d\frac{\partial}{\partial x_s}\((\widehat{\rho q_v})_ku_s\)-
\mu_q\sum_{s=1}^d\(\frac{\partial\rho}{\partial x_s}\frac{\partial}{\partial x_s}(\widehat{q_v})_k+
\rho\frac{\partial^2}{\partial x_s^2}(\widehat{q_v})_k\)&=(\widehat{r_1})_k,\notag\\
\frac{\partial}{\partial t}(\widehat{\rho q_c})_k+\sum_{s=1}^d\frac{\partial}{\partial x_s}\((\widehat{\rho q_c})_ku_s\)-
\mu_q\sum_{s=1}^d\(\frac{\partial\rho}{\partial x_s}\frac{\partial}{\partial x_s}(\widehat{q_c})_k+
\rho\frac{\partial^2}{\partial x_s^2}(\widehat{q_c})_k\)&=(\widehat{r_2})_k,\\
\frac{\partial}{\partial t}(\widehat{\rho q_r})_k-\frac{\partial}{\partial x_d}\widehat{\alpha}_k+
\sum_{s=1}^d\frac{\partial}{\partial x_s}\((\widehat{\rho q_r})_ku_s\)-
\mu_q\sum_{s=1}^d\(\frac{\partial\rho}{\partial x_s}\frac{\partial}{\partial x_s}(\widehat{q_r})_k+
\rho\frac{\partial^2}{\partial x_s^2}(\widehat{q_r})_k\)&=(\widehat{r_3})_k,\notag
\end{align}
for $k=0,\ldots,M$. Here, the coefficients $\{\widehat{\alpha}_k\}_{k=0}^M$ are obtained using the following expansion:
\begin{equation*}
v_q(\x,t,\omega)(\rho q_r)(\x,t,\omega)=
\sum_{j=0}^M(\widehat{v_q})_j(\x,t)\Phi_j(\omega)\sum_{m=0}^M(\widehat{\rho q_r})_m(\x, t)\Phi_m(\omega)=:
\sum_{k=0}^M\widehat{\alpha}_k(\x,t)\Phi_k(\omega).
\end{equation*}
The coefficients $\{(\widehat{r_1})_k,(\widehat{r_2})_k,(\widehat{r_3})_k\}_{k=0}^M$, as well as $\{\widehat\alpha_k\}_{k=0}^M$ are
calculated via discrete Legendre transform (DLT) and inverse discrete Legendre transform (IDLT), which can be briefly described as
follows.
\begin{itemize}
\item DLT: First, the Galerkin projection applied to the expansion $f(\x,t,\omega)=\sum_{k=0}^M\widehat f_k(\x,t)\Phi_k(\omega)$ yields
\begin{equation}
\widehat f_k(\x,t)=\frac{2k+1}{2}\int\limits_{-1}^1f(\x,t,\omega)\Phi_k(\omega)\,{\rm d}\omega\quad\mbox{for}\quad0\le k\le M.
\label{DLT}
\end{equation}
Then, approximating the integral in \eqref{DLT} using the Gauss-Legendre quadrature leads to
\begin{equation*}
\mbox{DLT}\left[\left\{f(\x,t,\omega_\ell)\right\}_{\ell=0}^M\right]=\left\{\widehat f_k(\x,t)\right\}_{k=0}^M=
\left\{\frac{2k+1}{2}\sum_{l=0}^M\beta_\ell f(\x,t,\omega_\ell)\Phi_k(\omega_\ell)\right\}_{k=0}^M,
\end{equation*}
where $\beta_\ell$ are the Gauss-Legendre quadrature weights and $\omega_\ell$ is the $\ell$-th root of $\Phi_{M+1}$.
\item IDLT: Given the coefficients $\{\widehat f_k\}_{k=0}^M$, we compute the function $f$ {\cred through} the gPC expansion
\begin{equation*}
\mbox{IDLT}\left[\left\{\widehat f_k(\x,t)\right\}_{k=0}^M\right]=\left\{f(\x,t,\omega_\ell)\right\}_{l=0}^M=
\left\{\sum_{k=0}^M\widehat f_k(\x,t)\Phi_k(\omega_\ell)\right\}_{\ell=0}^M.
\end{equation*}
\end{itemize}
Consequently, we obtain
\begin{equation*}
\{(\widehat{r_1})_k\}_{k=0}^M=\mbox{DLT}\left[R_1\(\mbox{IDLT}\left[\left\{(\widehat{\rho q_v})_k\right\}_{k=0}^M\right],
\mbox{IDLT}\left[\left\{(\widehat{\rho q_c})_k\right\}_{k=0}^M\right],\mbox{IDLT}\left[\left\{(\widehat{\rho q_r})_k\right\}_{k=0}^M
\right]\)\right],
\end{equation*}
and analogously for $\{(\widehat{r_2})_k\}_{k=0}^M$, $\{(\widehat{r_3})_k\}_{k=0}^M$ and $\{\widehat{\alpha}_k\}_{k=0}^M$.
\begin{remark}
We stress that since the values $\Phi_k(\omega_\ell)$, $0\le k,\ell\le M$, are needed every time either DLT or IDLT is applied, one can
{\color{blue} pre-compute them} for the code efficiency.
\end{remark}

For the spatial and temporal discretizations of the system \eqref{cloud_equations_coefficients}, we apply the same finite volume method as
described in Section \ref{discr_space} and the same large stability domain explicit time integration method mentioned in Section
\ref{discr_time}. As in the deterministic case, we implement the ODE solver DUMKA3, which we provide with the following time step stability
restriction for the forward Euler method:
\begin{equation*}
\max_{s=1,2}\,\max_{i=1,\ldots,N}(|(u_s)_i|,|(u_3)_i+v_q(\omega_l)|)\,\frac{\Delta t^n_{\rm cloud}}{h}<0.5,
\end{equation*}
which should be satisfied for all of the Legendre roots $\omega_\ell$, $\ell={\cblue 0},\ldots,M$.

\section{Stochastic numerical experiments}\label{sec7}
In this section, we conduct numerical experiments with the stochastic Galerkin method described in Section \ref{sec5} for the free
convection of a moist warm air bubble and the Rayleigh-B\'enard convection. We demonstrate the influence of uncertainty in initial data as
well as in cloud parameters on the solution of the coupled Navier-Stokes-cloud model \eqref{NS_equations_pert},
\eqref{cloud_equations_random}. In all of our numerical examples below, we take $M=3$. {\cred Our extensive tests, from which we present
here only a selected part, showed that this was sufficient. Indeed, as documented in Example 6, high-order stochastic coefficients typically
have very small influence on a solution (see Figure \ref{coeff}), and thus can be neglected. Similar behavior was observed in other
experiments.}

\subsection{Free convection of a smooth warm air bubble}
In this test, we modify Example 1 by randomly perturbing either the initial data or selected model parameters.

\subsubsection*{\cred Example 4: 2-D case with stochastic initial data}
{\cred We begin by considering the following experiment with a $10\%$ perturbation of the initial water vapor concentration:}
\begin{align*}
(\widehat{q_v})_0(\x,0)&=0.02\thetaprime(\x,0),~~~(\widehat{q_v})_1(\x,0)=0.1(\widehat{q_v})_0(\x,0),~~~(\widehat{q_v})_k(\x,0)=0\quad
\mbox{for}\quad2\le k\le M,\\
(\widehat{q_c})_k(\x,0)&=(\widehat{q_r})_k(\x,0)=0\quad\mbox{for}\quad0\le k\le M.
\end{align*}
We compute the solution using different meshes until the final time {\cred $t=10s$}.

\begin{figure}[ht!]
\centerline{\hspace*{0.5cm}\includegraphics[scale=0.70]{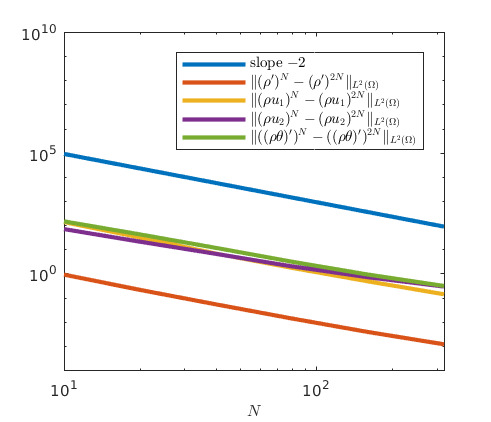}\hspace*{1cm}
\includegraphics[scale=0.70]{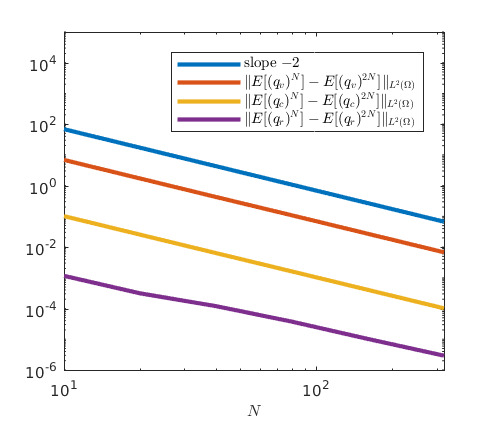}}
\caption{{\cred Example 4: Convergence study for the flow variables (left) and the expected values of $(\rho q_\ell)_k$, $\ell\in\{v,c,r\}$
(right) in space and time.\label{initdata_Ndt}}}
\end{figure}

\begin{figure}[ht!]
\centerline{\hspace*{0.5cm}\includegraphics[scale=1.05]{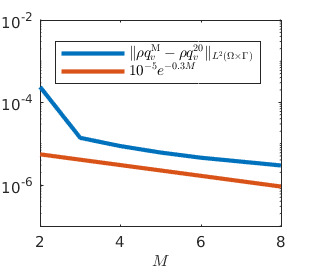}
\includegraphics[scale=1.05]{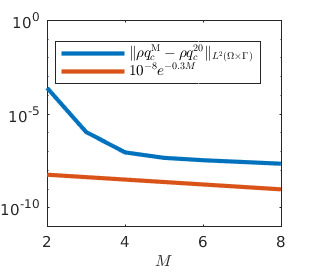}\includegraphics[scale=1.05]{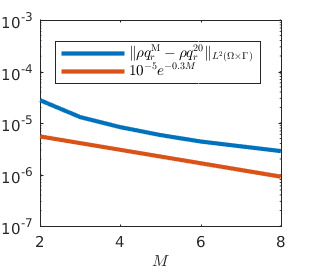}}
\caption{{\cred Example 4: Convergence study for the cloud variables $(\rho q_\ell)$, $\ell\in\{v,c,r\}$ in the stochastic space.
\label{fig100}}}
\end{figure}

{\cred The experimental convergence study for the cloud and flow variables is presented in Figure \ref{initdata_Ndt}}. Similarly to the
deterministic case, one can observe second-order convergence in space and time. In order to test the convergence in the stochastic space, we
obtain a reference solution computed by the stochastic {\cred Galerkin method with 20 stochastic modes.}
{\cred The convergence study is presented in Figure \ref{fig100}, where we plot the difference between the approximate and reference
solutions, both computed using a mesh with $160\times160$ cells and $\Delta t=0.01$ at time $t=10s$. One can clearly see a spectral
convergence with the rate $e^{-0.3M}$}.


\subsubsection*{\cred Example 5: 2-D case with stochastic parameters}
\begin{figure}[ht!]
\centerline{\hspace*{0.5cm}\includegraphics[scale=1.05]{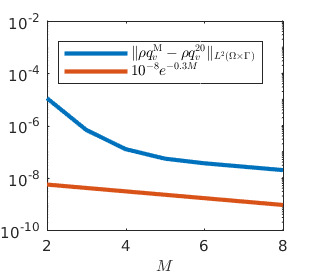}
\includegraphics[scale=1.05]{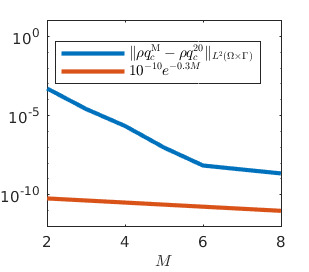}\includegraphics[scale=1.05]{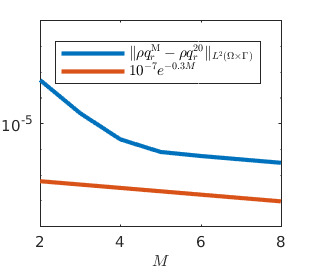}}
\caption{{\cred Example 5: Convergence study for the cloud variables $(\rho q_\ell)$, $\ell\in\{v,c,r\}$ in the stochastic space.
\label{fig10}}}
\end{figure}

{\cred In this experiment, we perturb the three selected model parameters ($k_1$, $k_2$ and $\alpha$) by $50\%$ each:
\begin{align*}
(\widehat{k_1})_0&=0.0041,\quad(\widehat{k_1})_1=0.5(\widehat{k_1})_0,\quad(\widehat{k_1})_k=0\quad\mbox{for}\quad2\le k\le M,\\
(\widehat{k_2})_0&=0.8,\quad(\widehat{k_2})_1=0.5(\widehat{k_2})_0,\quad(\widehat{k_2})_k=0\quad\mbox{for}\quad2\le k\le M,\\
\widehat\alpha_0&=190.3,\quad\widehat\alpha_1=0.5\widehat\alpha_0,\quad\widehat\alpha_k=0\quad\mbox{for}\quad2\le k\le M.
\end{align*}
These parameters were proposed in \cite{porz18} as the most sensitive model parameters. We study the convergence in the stochastic space. To this end, we plot in Figure~\ref{fig10} the difference between the approximate and
reference (obtained with 20 stochastic modes) solutions, both computed using a mesh with $160\times160$ cells and $\Delta t=0.01$ at time
$t=10s$. As in Example 4, one can observe a spectral convergence with the rate $e^{-0.3M}$.}

\subsubsection*{\cred Example 6: Comparison of the stochastic Galerkin and stochastic collocation methods in the 3-D case}
{\cred In this comparison test, we consider free convection of a moist smooth warm air bubble with perturbed initial data as in Example 4.
Numerical solutions are computed on $50\times50\times50$ mesh at time $t=1s$. We compare the performance of the stochastic Galerkin and
stochastic collocation methods.}
{\cblue The collocation method (see, e.g., \cite{collocation}) is an interpolation method in the stochastic space, which uses a
deterministic model with the values of the stochastic variable taken at collocation points suitably chosen on the interval $(-1,1)$; here,
we use the Gauss-Legendre points.}

{\cred
In Figure \ref{coeff}, the $L^2$-norms of the stochastic coefficients of the stochastic Galerkin
($\|(\widehat{\rho q_\ell})_m\|_{L^2(\Omega)}$) and stochastic collocation ($\|(\widetilde{\rho q_\ell})_m\|_{L^2(\Omega)}$) methods for
$\ell\in\{v,c,r\}$ and $m=0,\ldots,19$ are shown. One can observe an exponential decay with respect to $m$ and  good agreement between both
methods that demonstrates the reliability of the stochastic Galerkin method. We note, however, that the norms of the solutions computed by
these two different methods are not equal since the stochastic Galerkin method uses the expected values of the cloud variables in the
Navier-Stokes equations.

\begin{figure}[ht!]
\centerline{\hspace*{0.5cm}\includegraphics[scale=1.05]{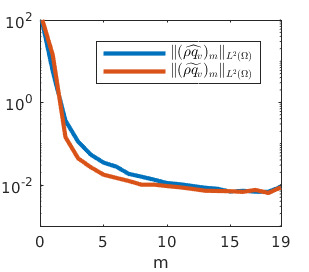}\includegraphics[scale=1.05]{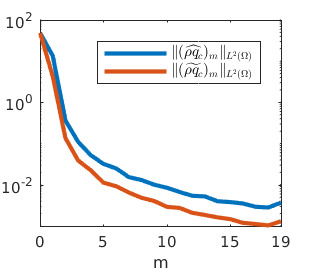}
\includegraphics[scale=1.05]{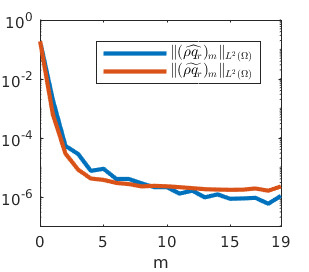}}
\caption{{\cred Example 6: $L^2(\Omega)$-norms of the stochastic coefficients computed by the stochastic Galerkin
($(\widehat{\rho q_\ell})_m$) and stochastic collocation ($(\widetilde{\rho q_\ell})_m$) methods, $\ell\in\{v,c,r\}$.\label{coeff}}}
\end{figure}

In Figure \ref{CPU}, we compare the CPU times consumed by the stochastic Galerkin and stochastic collocation methods with the same number
$M$ of stochastic modes/collocation points. Since the stochastic Galerkin method solves the Navier-Stokes equations just once instead of $M$
times, as needed by the stochastic collocation method does, it is expected to outperform the stochastic collocation method. This has been
confirmed by our simulations.}

\begin{figure}[ht!]
\centerline{\includegraphics[scale=0.8]{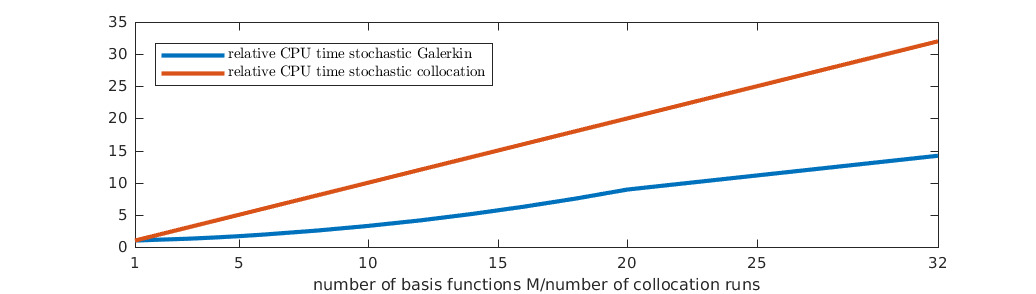}}
\caption{{\cred Example 6: Relative CPU times consumed by the stochastic Galerkin and stochastic collocation methods.\label{CPU}}}
\end{figure}

\subsection{Rayleigh-B\'{e}nard convection}
In this section, we present results of uncertainty study for the Rayleigh-B\'enard convection in both 2-D and 3-D. We investigate
uncertainty propagation, which is triggered either by the initial data or cloud parameters.

\subsubsection*{Example 7: 2-D case with stochastic initial data}
{\cred In this experiment, we choose the same initial data for the flow variables as in Section \ref{RB_det} and perturb the initial data in
$q_v$ by $5\%$, $10\%$, $20\%$ and $50\%$:
\begin{equation}
\begin{aligned}
(\widehat{q_v})_0(\x,0)&=2\cdot10^{-5}\thetabar,\,(\widehat{q_v})_1(\x,0)=\nu(\widehat{q_v})_0(\x,0),\,(\widehat{q_v})_k=0\quad\mbox{for}
\quad2\le k\le M,\\
(\widehat{q_c})_k(\x,0)&=(\widehat{q_r})_k(\x,0)=0\quad\mbox{for}\quad0\le k\le M,
\end{aligned}
\label{ic_rb}
\end{equation}
where $\nu=0.05$, $0.1$, $0.2$ and $0.5$, respectively. {\cred It should be observed that uniform perturbations in the initial conditions
for $q_v$ may lead to either reduced or enhanced water vapor concentrations as compared to the deterministic simulations. Since the
temperature gradient is quite small, a change in $q_v$ translates to an (almost) linear change in saturation ratio, which directly controls
cloud formation. Thus, in the case of positive perturbations, a higher water vapor concentration leads to earlier cloud formation and, in
addition, a higher potential temperature change since more water is available in the system. On the other hand, lower values of $q_v$ lead
to a time delay in the formation of clouds, even if small convective cells are driven by the dry unstable situation. In a feedback cycle, a
reduced or even delayed formation of cloud water propagates further to a weaker rain formation. Finally, the evaporation of rain water
leads to a strong cooling effect of the lower layers of the domain, which also crucially depends on the amount of sedimenting rain water.
These effects have to be taken into account for the evaluation of the different perturbation scenarios.}

Numerical solutions for both the potential temperature $\theta$ (Figure \ref{fig12}) and the expected value of the cloud drops concentration
$q_c$ (Figure \ref{fig13}) are computed at time $t=1400s$ using $320\times320$ mesh cells and presented for $\nu=0$, $0.05$, $0.1$, $0.2$
and $0.5$. For a better comparison, we have used the same range of values for different perturbations in all of the plots.

\begin{figure}[ht!]
\centerline{\includegraphics[scale=1]{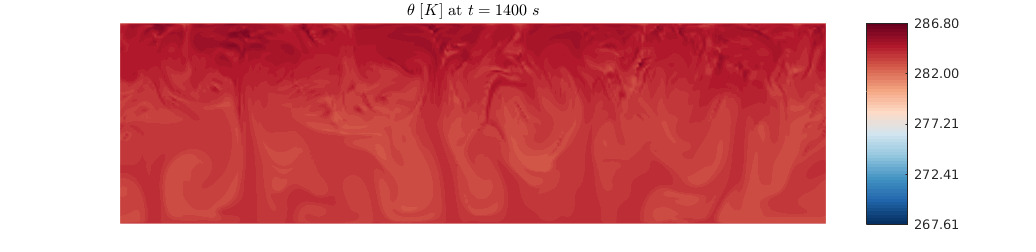}}
\centerline{\includegraphics[scale=1]{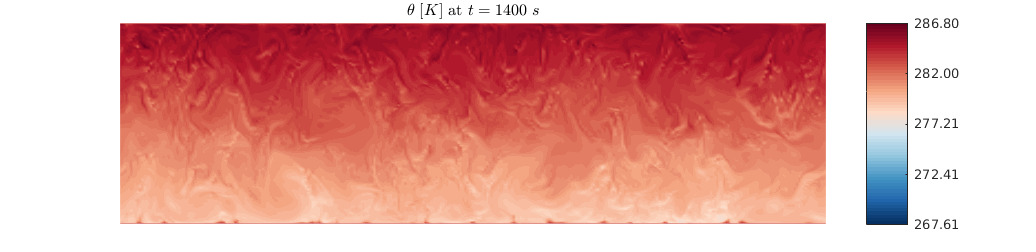}}
\centerline{\includegraphics[scale=1]{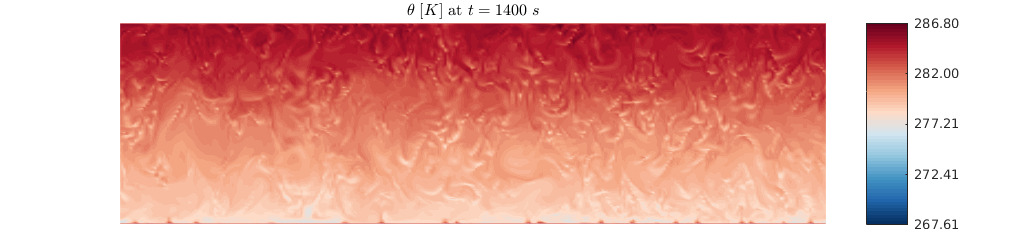}}
\centerline{\includegraphics[scale=1]{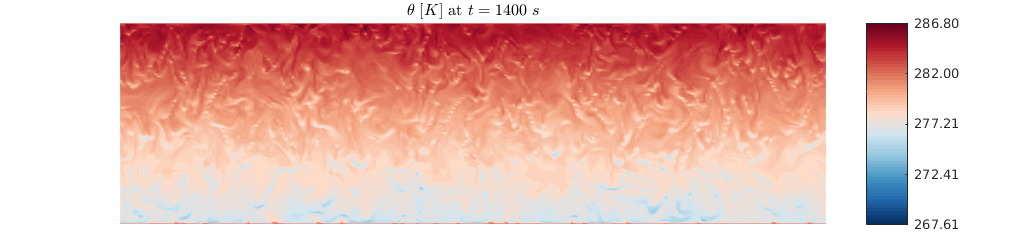}}
\centerline{\includegraphics[scale=1]{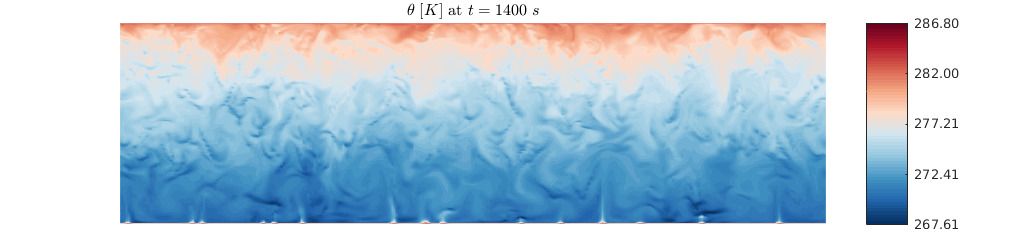}}
\caption{Example 7: Potential temperature $\theta$ for $\nu=0$, $0.05$, $0.1$, $0.2$ and $0.5$ (from top to down).\label{fig12}}
\end{figure}
\begin{figure}[ht!]
\centerline{\includegraphics[scale=1]{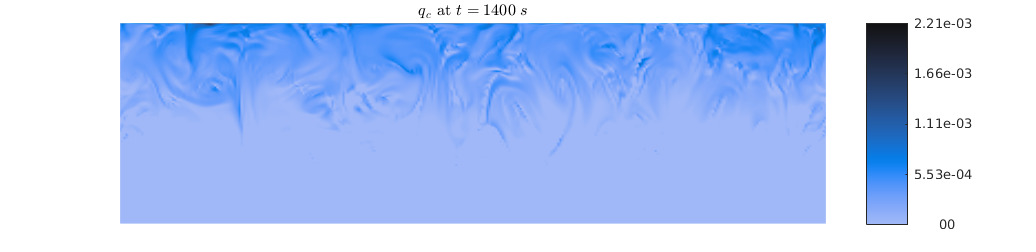}}
\centerline{\includegraphics[scale=1]{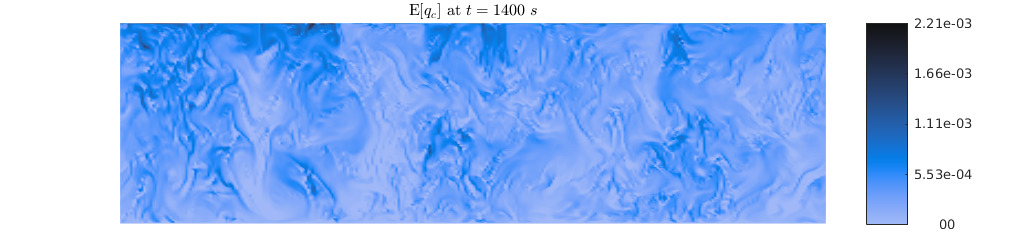}}
\centerline{\includegraphics[scale=1]{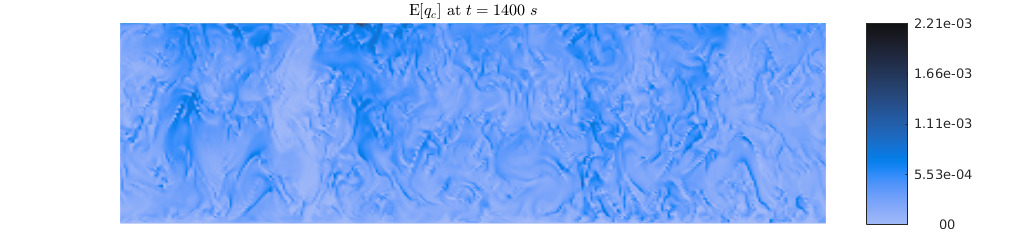}}
\centerline{\includegraphics[scale=1]{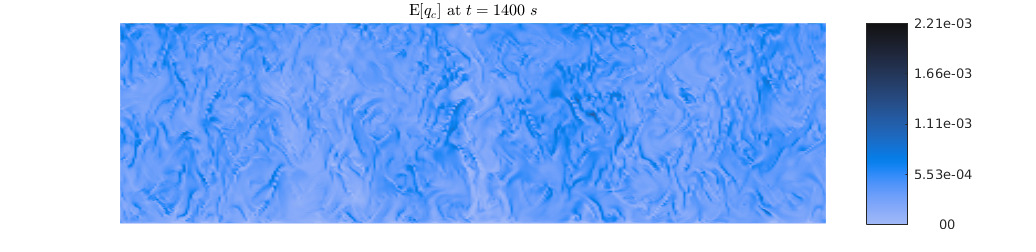}}
\centerline{\includegraphics[scale=1]{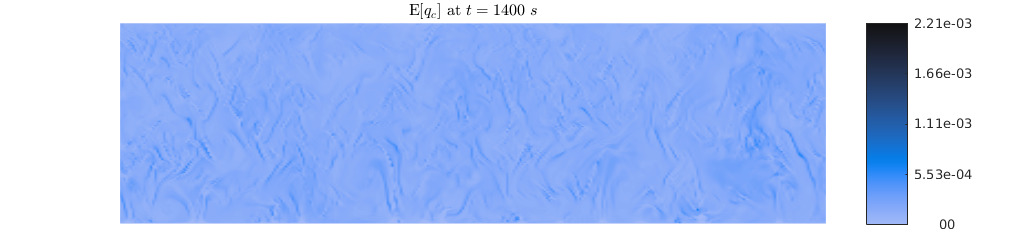}}
\caption{Example 7: Cloud drops concentration $q_c$ for $\nu=0$, $0.05$, $0.1$, $0.2$ and $0.5$ (from top to down).\label{fig13}}
\end{figure}

{\cred From the results shown in Figure \ref{fig12} for the potential temperature $\theta$, one can observe two major features that
increase in strength with increasing strength of perturbations. First, larger perturbations lead to more concise fine structures, that is,
for the deterministic simulation ($\nu=0$), the variable is quite smooth, whereas for large perturbations, large gradients on a very small
scale appear. This feature can also be recognized for the cloud water, that is, for the expected values $\mathbb{E}[q_c]$ of cloud water
concentrations in Figure \ref{fig13}. This is probably due to the fact that even small variations in water vapor have a strong impact on
cloud formation, since the activation of cloud droplet is basically a threshold process. For values closer to saturation, even small
variations in vertical upward motions can trigger cloud formation. Thus, the small scale variations are more prominent in the massively
perturbed scenarios. The second feature is even more striking. Increasing the perturbation in the initial water vapor distribution leads to
a stronger vertical gradient in potential temperature, that is, at low levels the temperature is much smaller than in the deterministic
case. This feature is due to cooling of sedimenting rain water. For simulations with a high water vapor loading, more cloud and thus more
rain is formed, which is subsequently falling down into lower levels and cools the environment by evaporation. Since this process is very
effective, the temperature can be reduced drastically. Note that this feature is well-known for the real atmosphere: Falling rain can cool
lower levels efficiently, so that a transition from melting rain droplets to snow can be possible for winter seasons. The efficient
formation of rain also leads to a strong reduction in cloud water, since accretion can eat up cloud droplets in the lower level; see
Figure~\ref{fig13}. For reduced values of initial water vapor, processes of cloud and precipitation formation is strongly reduced. However,
on average the positive perturbations dominate the statistical picture.}

In Figure \ref{fig14}, we plot the standard deviation of the cloud drops concentration. Here, we compute the standard deviation for a
function $f(\x,t,\omega)=\sum_{k=0}^M \widehat{f_k}(\x,t)\Phi_k(\omega)$ as
\begin{equation}
\sigma(f(\x,t,\omega))=\sqrt{\sum_{k=1}^M \frac{1}{2k+1}\widehat{f_k}(\x,t)^2},
\end{equation}
which follows from the orthogonal property \eqref{orthogonal_property}.

\begin{figure}[ht!]
\centerline{\includegraphics[scale=1]{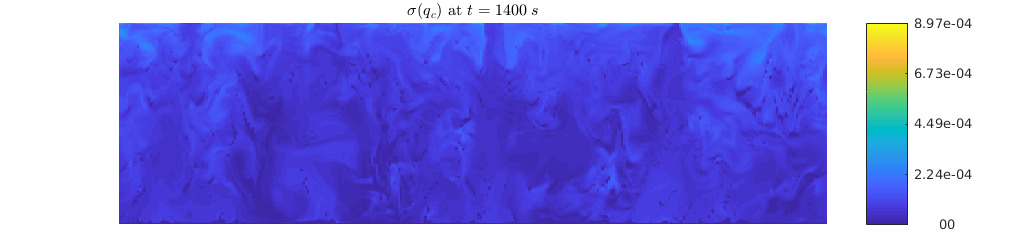}}
\centerline{\includegraphics[scale=1]{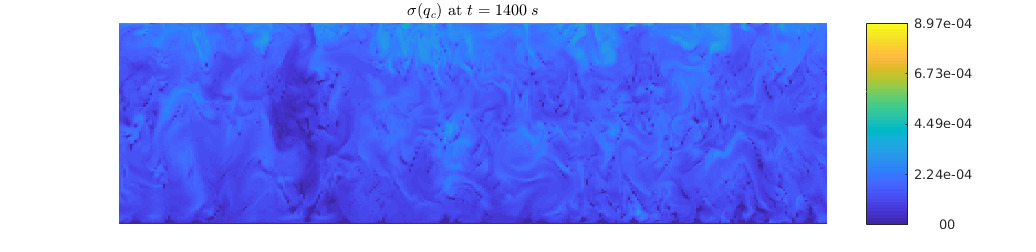}}
\centerline{\includegraphics[scale=1]{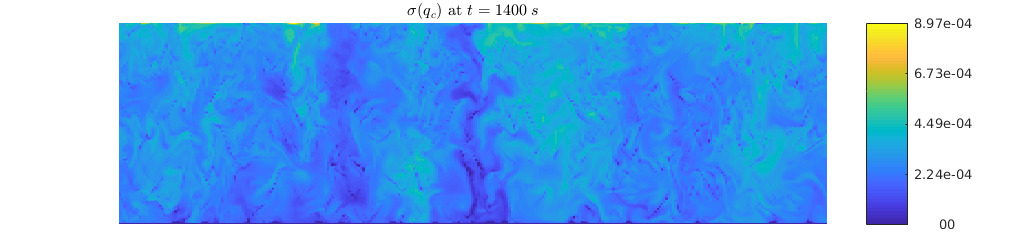}}
\centerline{\includegraphics[scale=1]{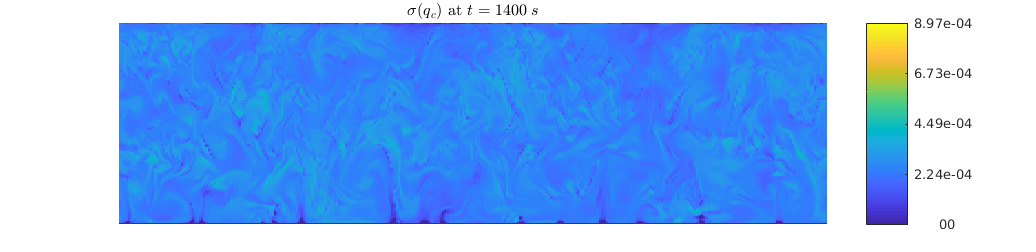}}
\caption{Example 7: Standard deviation of the cloud drops concentration $q_c$ for $\nu=0.05$, $0.1$, $0.2$ and $0.5$ (from top to down).
\label{fig14}}
\end{figure}

In order to investigate the influence of perturbations, the time evolution of the mean expected value per $m^2$ as well as the mean
standard deviation per $m^2$ for the cloud variables are presented in Figure \ref{fig15}. In $d$-space dimensions these quantities can be
computed in the following way:
$$
\begin{aligned}
\mathbb{E}\left[\frac{h^d}{|\Omega|}\sum_{i=1}^N(q_\ell)_i\right]&=\frac{h^d}{|\Omega|}\sum_{i=1}^N\mathbb{E}\left[(q_\ell)_i\right]=
\frac{h^d}{|\Omega|}\sum_{i=1}^N\widehat{((q_\ell)_i)_0},\\
\sigma\left(\frac{h^d}{|\Omega|}\sum_{i=1}^N(q_\ell)_i\right)&=
\frac{h^d}{|\Omega|}\sqrt{\sum_{k=1}^M \left(\sum_{i=1}^N \widehat{((q_\ell)_i)_k}\right)^2 \frac{1}{2k+1}},
\end{aligned}
$$
where $N$ is the number of mesh cells and $\ell\in\{v,c,r\}$.}

\begin{figure}[ht!]
\centerline{\hspace*{0.25cm}\includegraphics[scale=0.80]{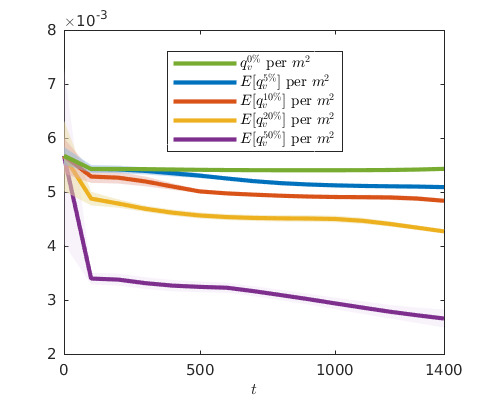}\hspace*{0.5cm}
\includegraphics[scale=0.80]{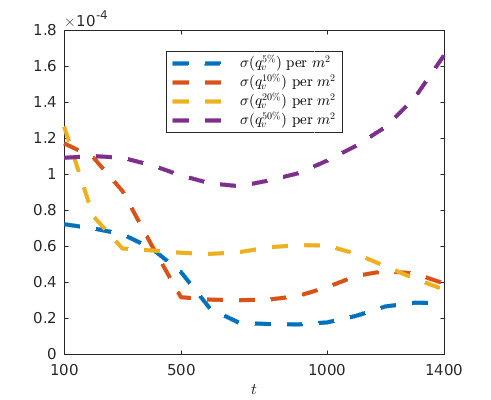}}

\vspace*{0.5cm}
\centerline{\hspace*{0.25cm}\includegraphics[scale=0.80]{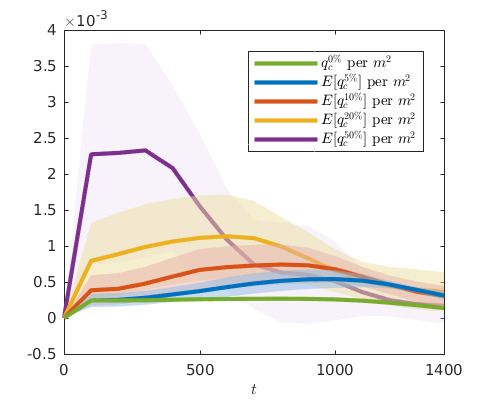}\hspace*{0.5cm}
\includegraphics[scale=0.80]{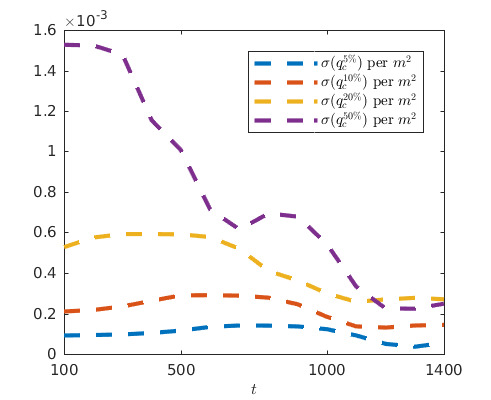}}

\vspace*{0.5cm}
\centerline{\hspace*{0.25cm}\includegraphics[scale=0.80]{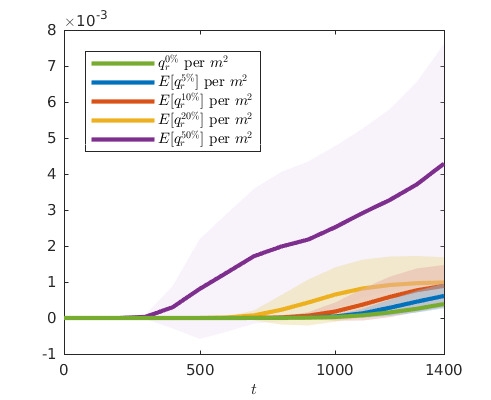}\hspace*{0.5cm}
\includegraphics[scale=0.80]{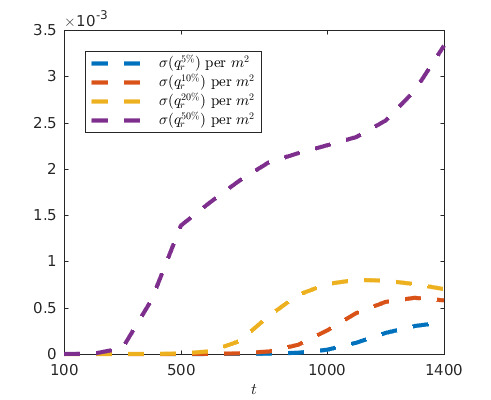}}
\caption{{Example 7: Time evolution of the expected values with their standard deviations for the cloud variables per $m^2$ (shaded region,
left column) and standard deviation (right column) for $\nu=0.05$, $0.1$, $0.2$ and $0.5$ (from top to down).\label{fig15}}}
\end{figure}

{\cred These averaged quantities show the qualitative difference in the different perturbations scenarios, as already described above.
For low perturbations, the difference of the expectation value is quite small and also the standard deviation remains small. For increasing
perturbations, the spread is increased. As noted before, the averaged quantities are dominated by the positive perturbations, leading to
(i) earlier cloud formation, (ii) thicker clouds due to more available water vapor, and (iii) to enhanced rain formation. These three
features can be seen very nicely in the strongest perturbations ($\nu=0.5$), with a large drop in water vapor concentration accompanied by a
strong increase in cloud water and earlier onset of precipitation. We would also like to note that the spread is only given by the standard
deviation, whereas the actual minima (for instance, almost no cloud formation) cannot be seen directly, although these scenarios are
possible.}

\subsubsection*{Example 8: 2-D case with stochastic parameters} 
\begin{figure}[ht!]
\centerline{\includegraphics[scale=1]{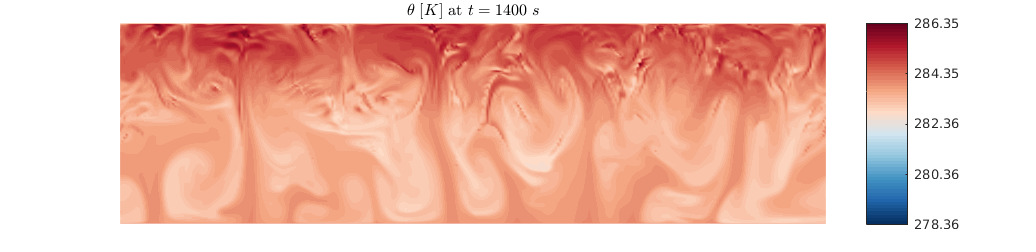}}
\centerline{\includegraphics[scale=1]{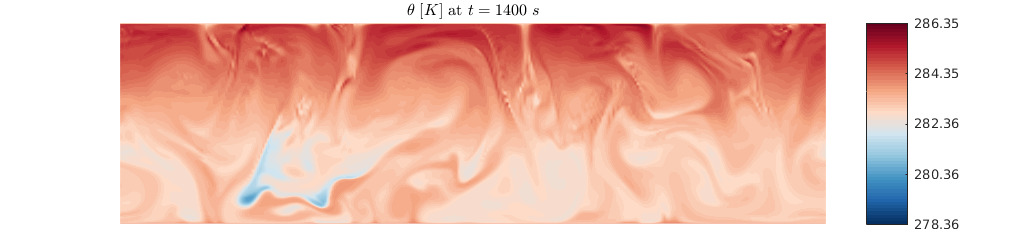}}
\centerline{\includegraphics[scale=1]{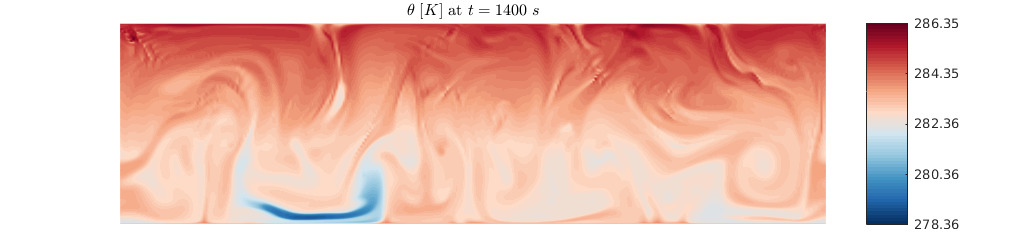}}
\centerline{\includegraphics[scale=1]{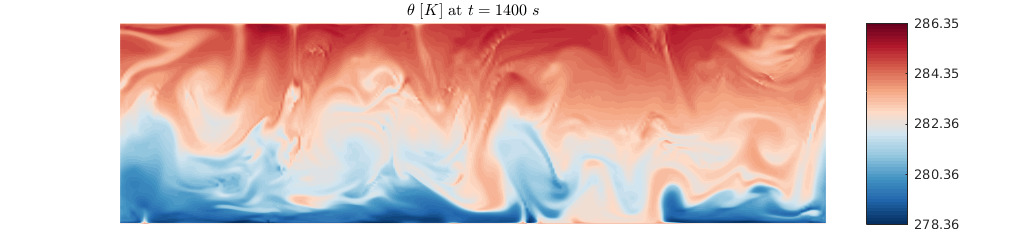}}
\caption{{\cblue Example 8: Potential temperature $\theta$ for $\nu=0$, $0.1$, $0.2$ and $0.5$ (from top to down).\label{fig26}}}
\end{figure}

\begin{figure}[ht!]
\centerline{\includegraphics[scale=1]{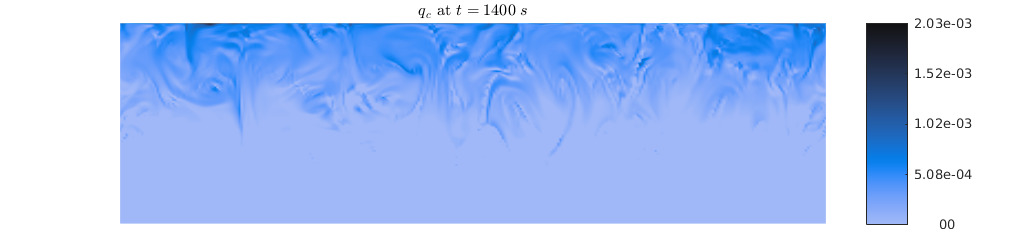}}
\centerline{\includegraphics[scale=1]{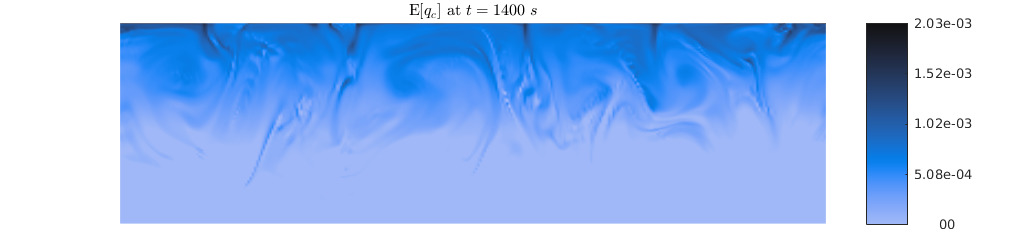}}
\centerline{\includegraphics[scale=1]{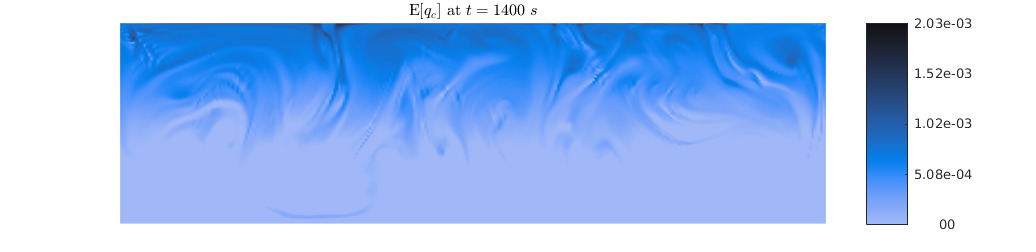}}
\centerline{\includegraphics[scale=1]{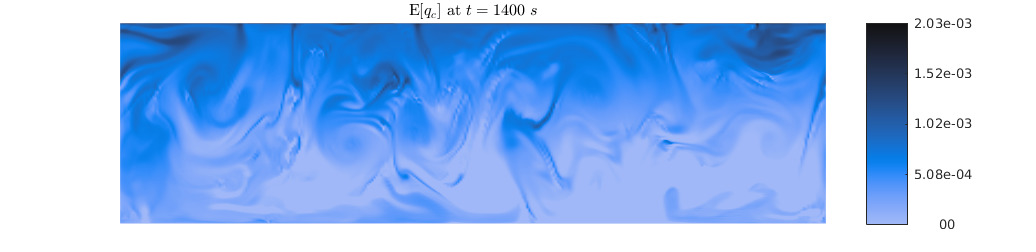}}
\caption{{\cblue Example 8: Cloud drops concentration $q_c$ for $\nu=0$, $0.1$, $0.2$ and $0.5$ (from top to down).\label{fig27}}}
\end{figure}

In the following experiment we study uncertainty propagation due to
incomplete information about the model parameters which is a very typical problem arising in atmospheric science. We chose the same initial
data for the flow and cloud variables as in Section \ref{RB_det}. More precisely, we take the following initial cloud variables:
$$
(\widehat{q_v})_0(\x,0)=0.02\thetaprime(\x,0),\quad(\widehat{q_v})_k(\x,0)=(\widehat{q_c})_k(\x,0)=(\widehat{q_r})_k(\x,0)=0\quad\mbox{for}
\quad1\le k\le M.
$$
Consequently, in order to investigate uncertainty propagation in the numerical solution we choose $10\%$, $20\%$ and $50\%$ perturbation of
these coefficients, namely, we take
\begin{align*}
(\widehat{k_1})_0&=0.0041,\quad(\widehat{k_1})_1=\nu(\widehat{k_1})_0,\quad(\widehat{k_1})_k=0\quad\mbox{for}\quad2\le k\le M,\\
(\widehat{k_2})_0&=0.8,\quad(\widehat{k_2})_1=\nu(\widehat{k_2})_0,\quad(\widehat{k_2})_k=0\quad\mbox{for}\quad2\le k\le M,\\
\widehat\alpha_0&=190.3,\quad\widehat\alpha_1=\nu\widehat\alpha_0,\quad\widehat\alpha_k=0\quad\mbox{for}\quad2\le k\le M,
\end{align*}
where $\nu=0.1$, 0.2 and 0.5, respectively. The numerical solution is computed at time {\cblue $t=1400s$ on a $320\times 320$ mesh. Figures
\ref{fig26} and \ref{fig27} present the potential temperature $\theta$ and the expected values of the cloud drops concentration $q_c$,
respectively, for $\nu=0$, $0.1$, $0.2$ and $0.5$.}


\begin{figure}[ht!]
\centerline{\includegraphics[scale=1]{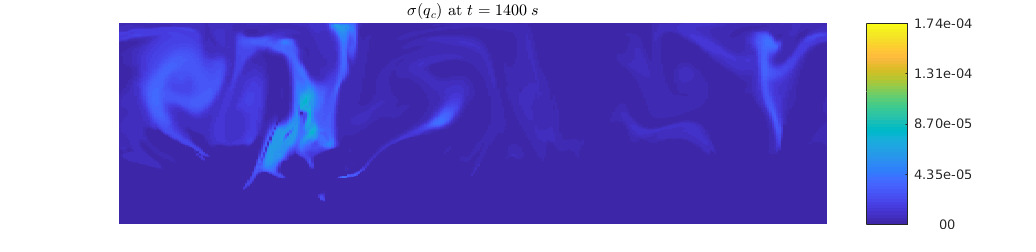}}
\centerline{\includegraphics[scale=1]{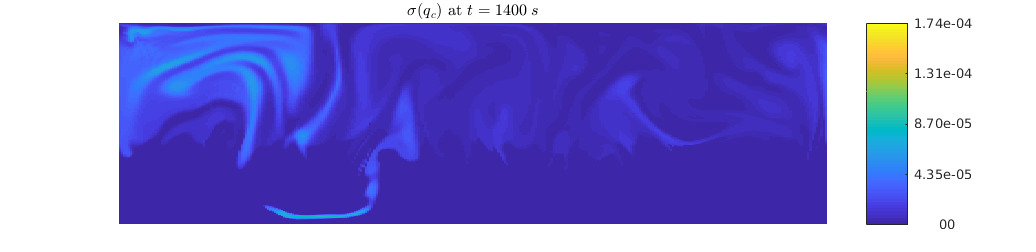}}
\centerline{\includegraphics[scale=1]{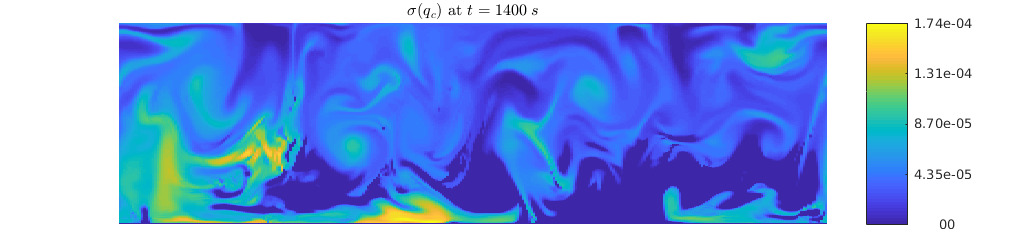}}
\caption{{\cblue Example 8: Standard deviation of the cloud drops concentration $q_c$ for $\nu=0.1$, $0.2$ and $0.5$ (from top to down).
\label{fig28}}}
\end{figure}

{\cblue In these scenarios with perturbed cloud model parameters, the overall structures are more stable and changes are less pronounced
than in Example 7. Nevertheless, the main features of the variations are obviously driven by precipitation processes since the only
perturbed parameters are those determining rain processes. Again, one key feature of the perturbed scenarios is the signature of
evaporating rain in lower levels of the 2-D domain. For positive perturbations, that is, larger parameters $k_1$, $k_2$ and $\alpha$, rain
formation is enhanced (more rain is formed from cloud water, due to larger $k_1$ and $k_2$) and sedimentation is enhanced (more rain
is falling downwards due to larger $\alpha$). Thus, more rain water is transported downwards into subsaturated regions, which is then
evaporated inducing cooling due to latent heat consumption. These positive variations again dominate the potential temperature field; see
Figure \ref{fig26}. For the cloud water field (Figure \ref{fig27}), one can observe higher expected values for stronger perturbations. This
is probably due to the fact that more rain is evaporated in lower levels, thus more water vapor is then available for cloud formation in
upward motions of the convective cells. This redistribution of water vapor as well as the reduction of rain for negative perturbation lead
to larger variations of cloud water in the 2-D domain, as can be seen in the standard deviation for the cloud drops concentration $q_c$, as
depicted in Figure \ref{fig28}.}

\begin{figure}[ht!]
\centerline{\hspace*{0.25cm}\includegraphics[scale=0.80]{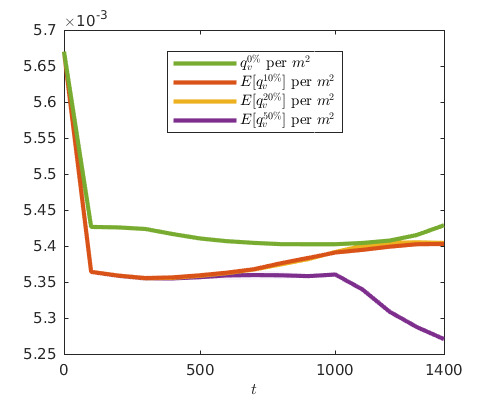}\hspace*{0.5cm}
\includegraphics[scale=0.80]{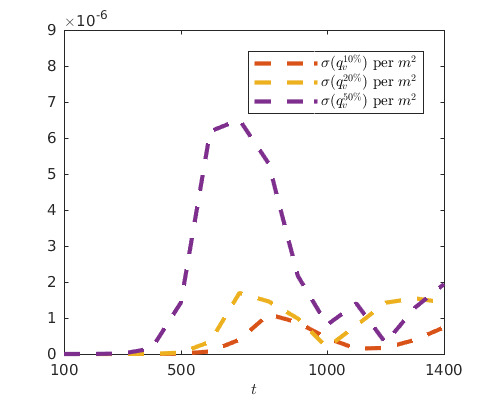}}

\vspace*{0.5cm}
\centerline{\hspace*{0.25cm}\includegraphics[scale=0.80]{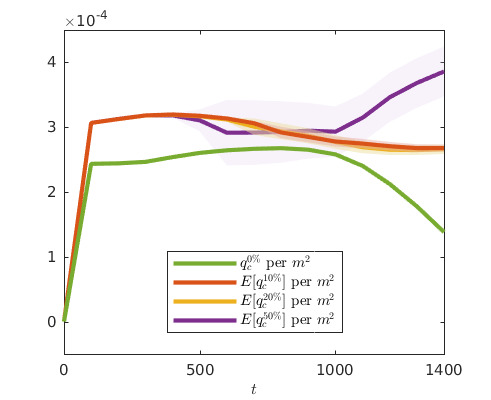}\hspace*{0.5cm}
\includegraphics[scale=0.80]{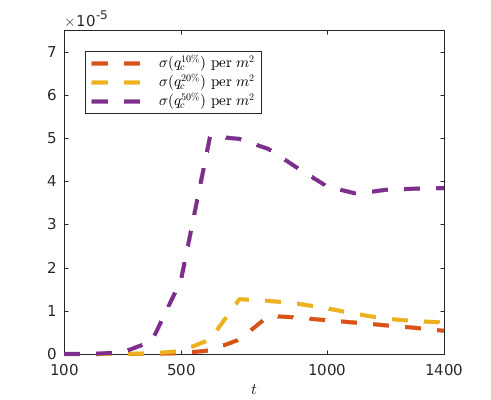}}

\vspace*{0.5cm}
\centerline{\hspace*{0.25cm}\includegraphics[scale=0.80]{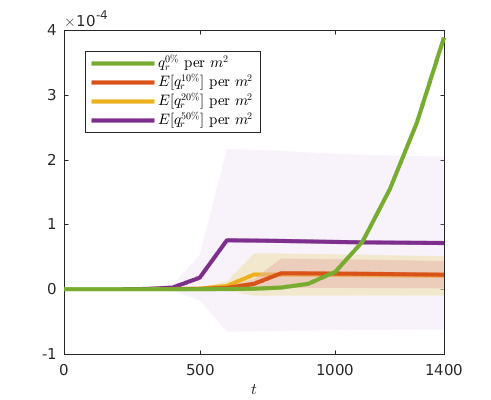}\hspace*{0.5cm}
\includegraphics[scale=0.80]{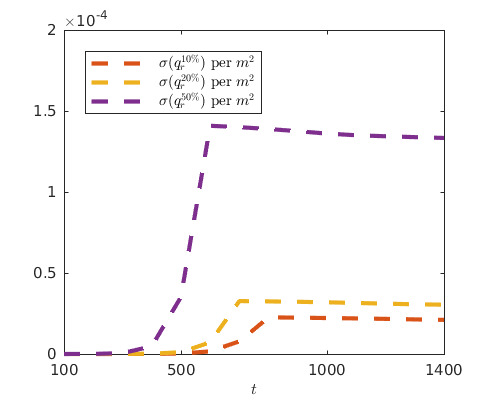}}
\caption{Example 8: Time evolution of the expected values with their standard deviations for the cloud variables per $m^2$ (shaded region,
left column) and standard deviation (right column) for $\nu=0.1$, $0.2$ and $0.5$ (from top to down).\label{fig29}}
\end{figure}

{\cblue The expected values as well as the standard  deviation for the cloud variables vary considerably with respect to the perturbation of
model parameters; see Figure \ref{fig29}. Our numerical experiments indicate that the standard deviation increases in time and also depends
on the size of the parameter perturbation. Indeed, the larger the parameter perturbation, the higher is the standard deviation of the cloud
variables. The size of the corresponding standard deviations is depicted in the {\cblue right} column of Figure~\ref{fig29}}. {\cblue As
expected, rain formation sets on earlier for large perturbations, leading also to a strong decrease in the overall water
vapor in comparison to less perturbed scenarios. However, the variations in all water variables increase a lot from the onset of
precipitation to later times, and even the spread in cloud water increases in contrast to the time evolution in Example 7. Generally, the spread in the mean
    values is smaller than in Example~7; changes in initial data
    produce a larger variation, that means the model physics is quite stable
    with respect to  perturbations in rain process formulations.}

\subsubsection*{Example 9: 3-D case with stochastic initial data}
\begin{figure}[ht!]
\includegraphics[scale=0.18]{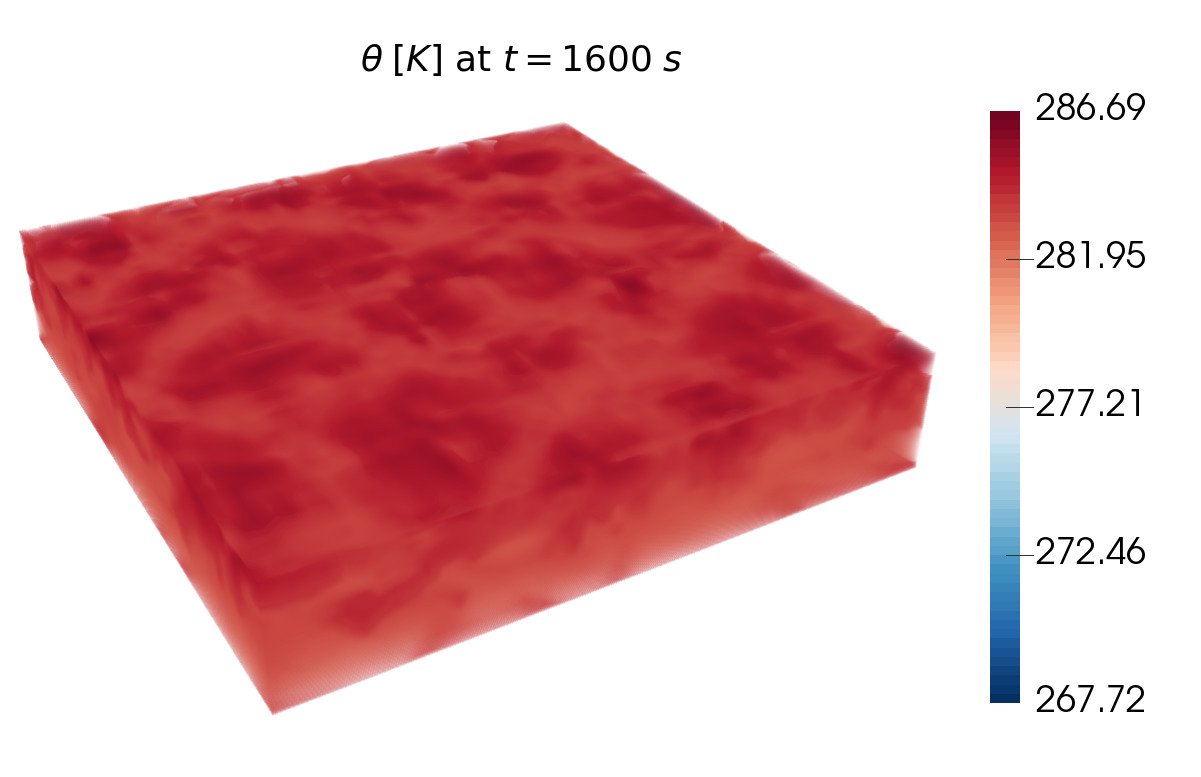}\hspace*{0.5cm}\includegraphics[scale=0.18]{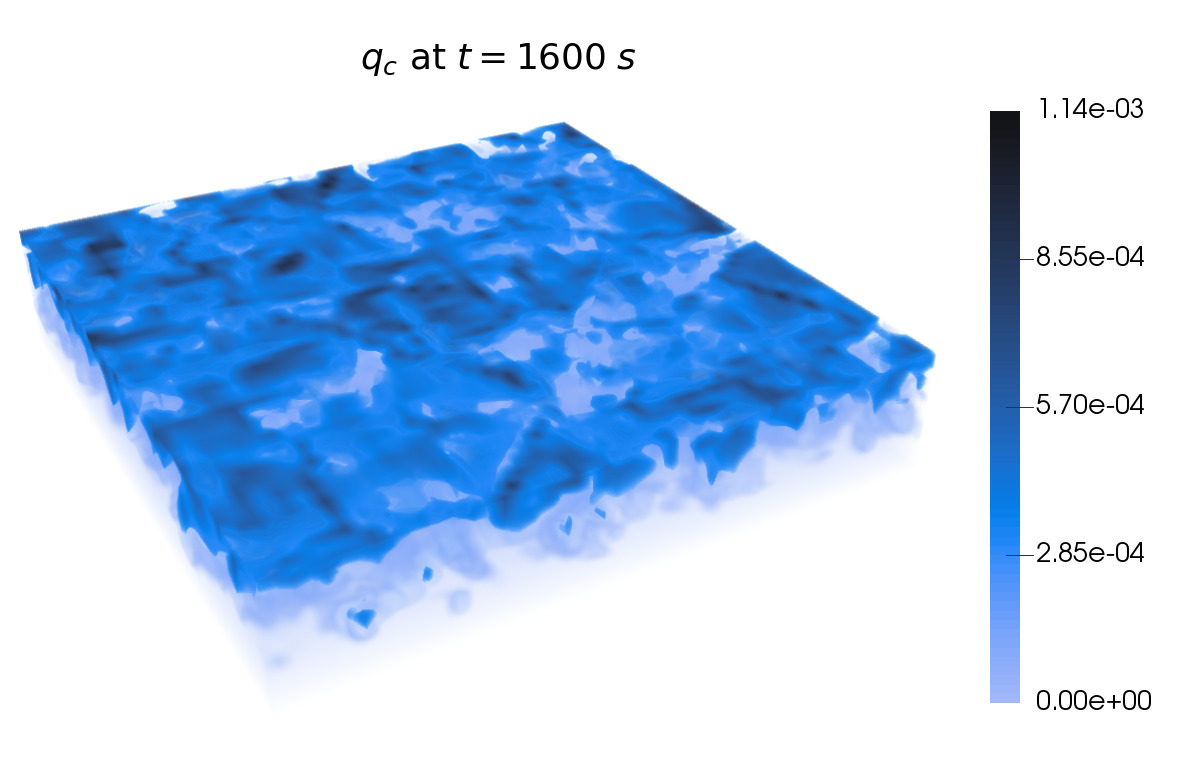}
\includegraphics[scale=0.18]{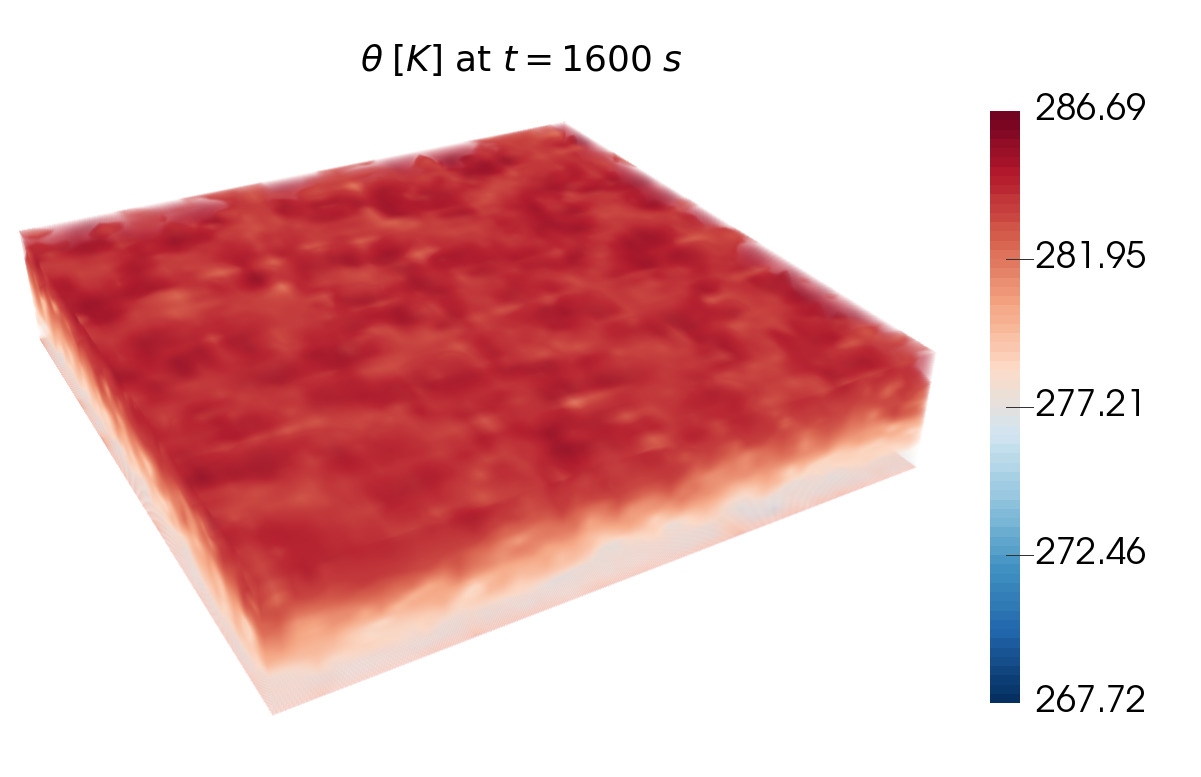}\hspace*{0.5cm}\includegraphics[scale=0.18]{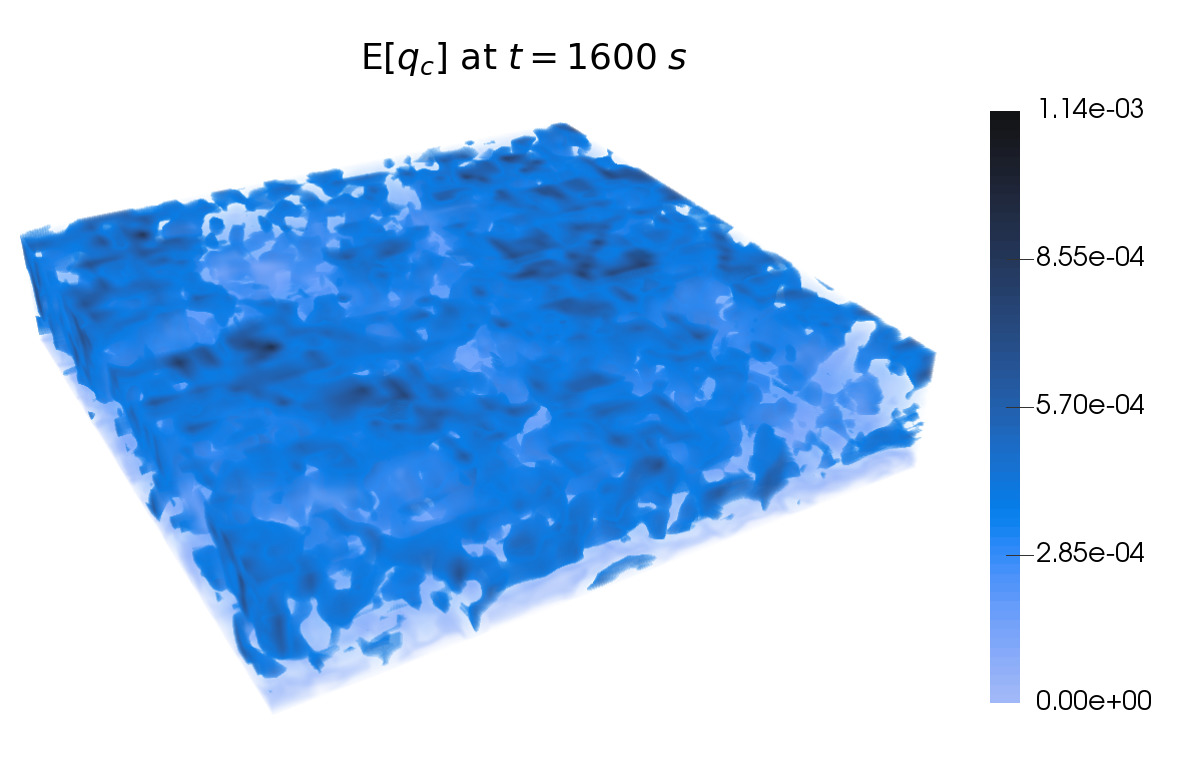}
\includegraphics[scale=0.18]{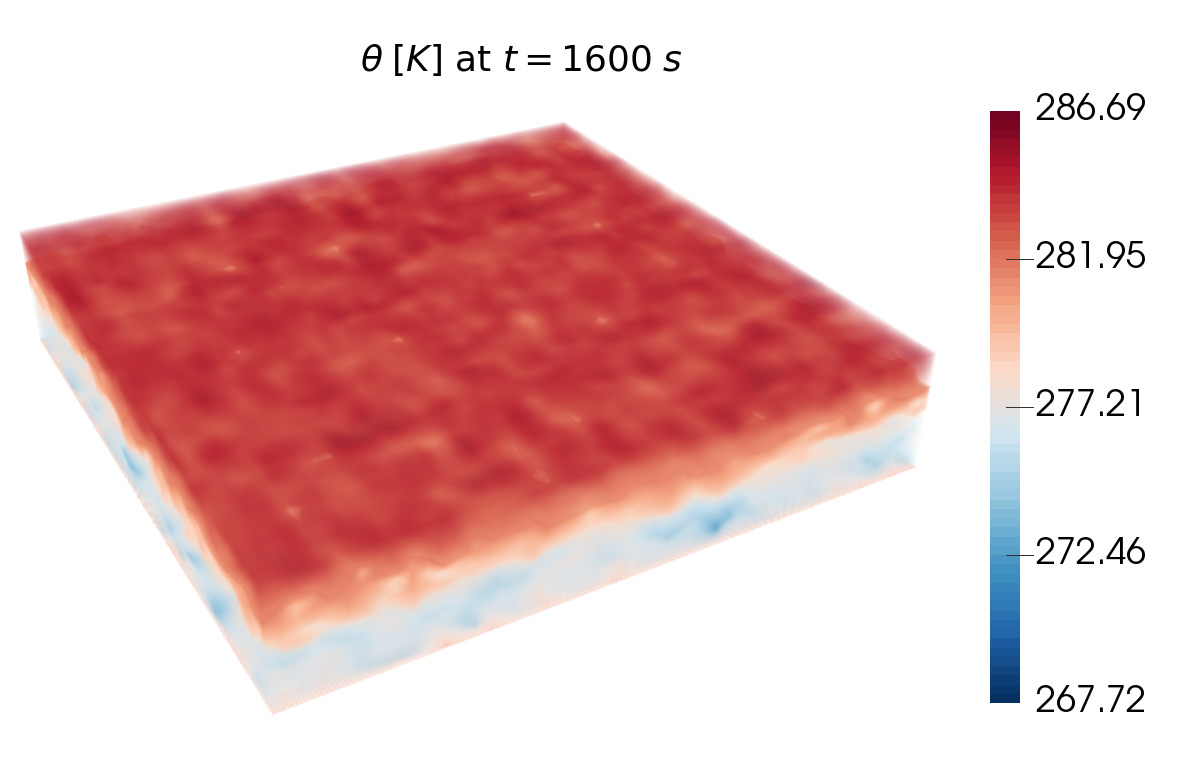}\hspace*{0.5cm}\includegraphics[scale=0.18]{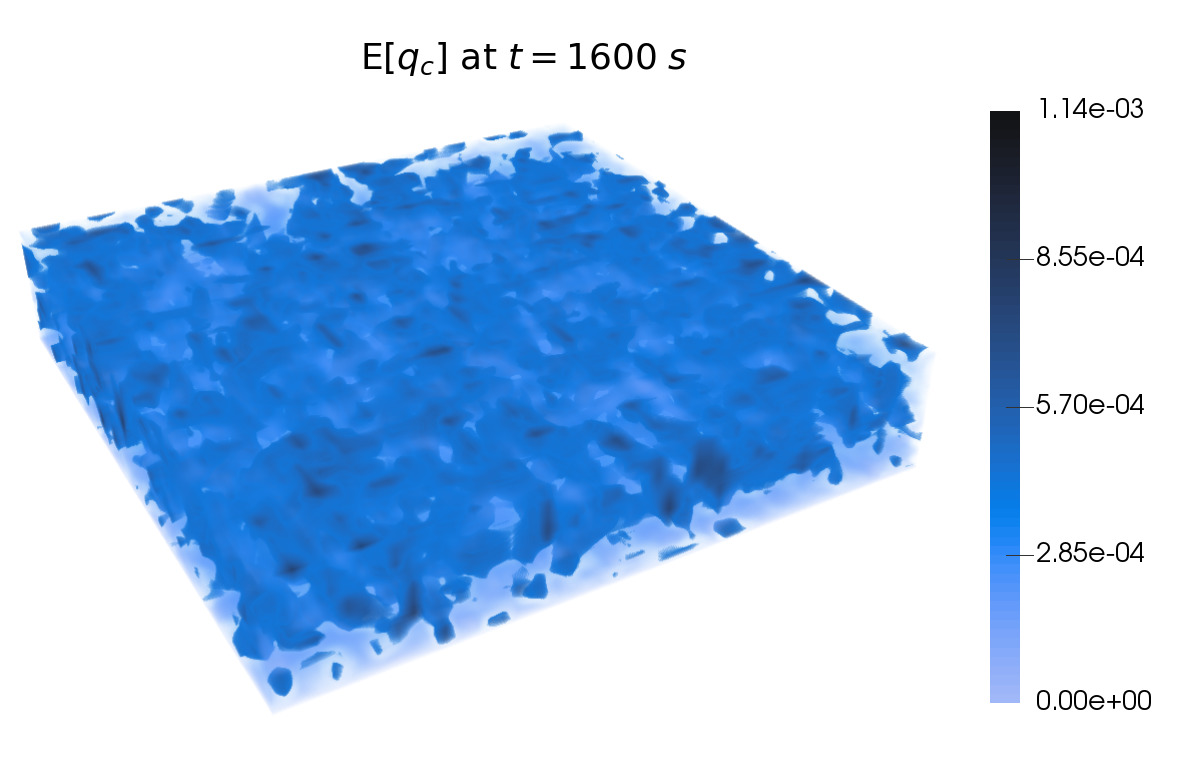}
\includegraphics[scale=0.18]{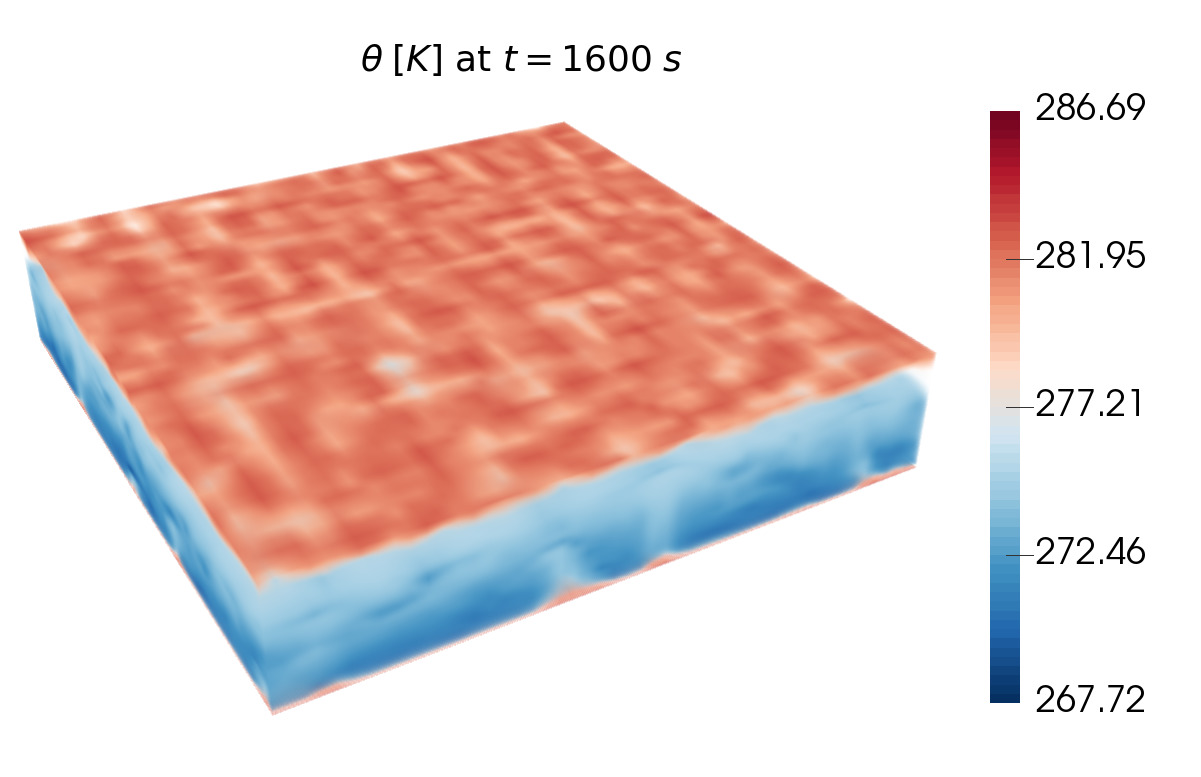}\hspace*{0.5cm}\includegraphics[scale=0.18]{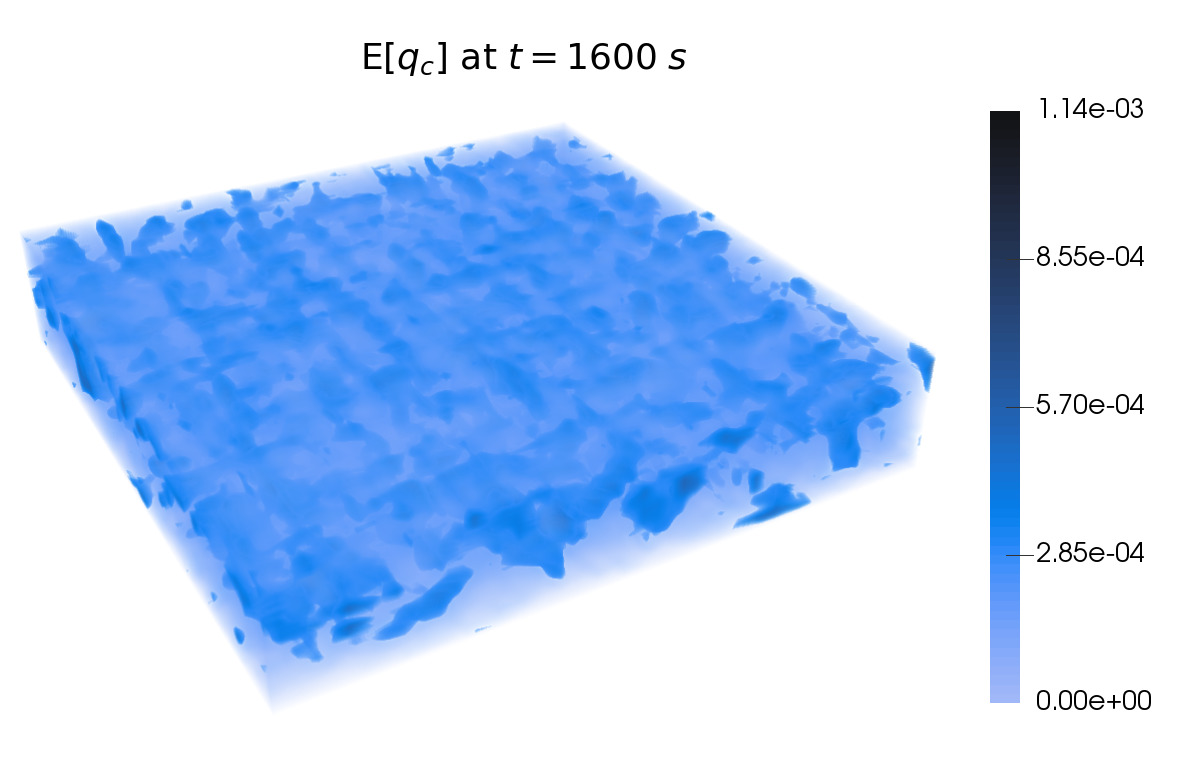}
\caption {{\cblue Example 9: Potential temperature $\theta$ (left column) and cloud drops concentration $q_c$ (right column)
using $0\%$, $10\%$, $20\%$ and $50\%$ (from top to down) perturbation of the initial data in $q_v$.\label{fig33}}}
\end{figure}

\begin{figure}[ht!]
\centerline{\hspace*{0.25cm}\includegraphics[scale=0.80]{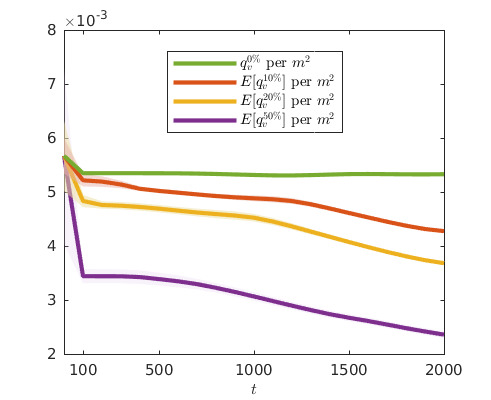}\hspace*{0.5cm}
\includegraphics[scale=0.80]{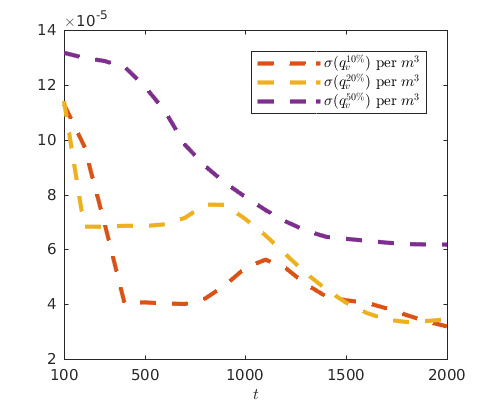}}

\vspace*{0.5cm}
\centerline{\hspace*{0.25cm}\includegraphics[scale=0.80]{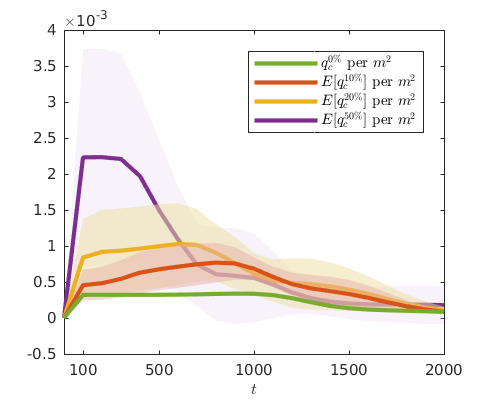}\hspace*{0.5cm}
\includegraphics[scale=0.80]{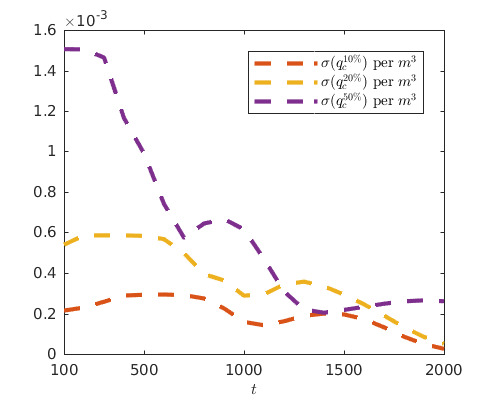}}

\vspace*{0.5cm}
\centerline{\hspace*{0.25cm}\includegraphics[scale=0.80]{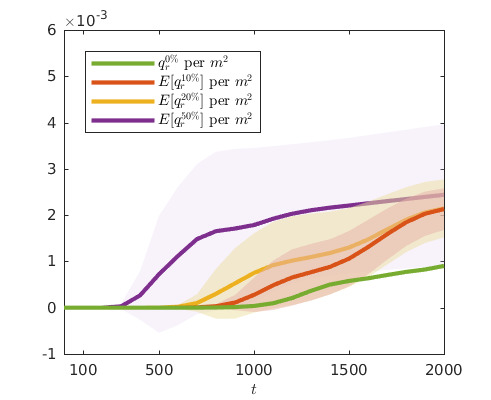}\hspace*{0.5cm}
\includegraphics[scale=0.80]{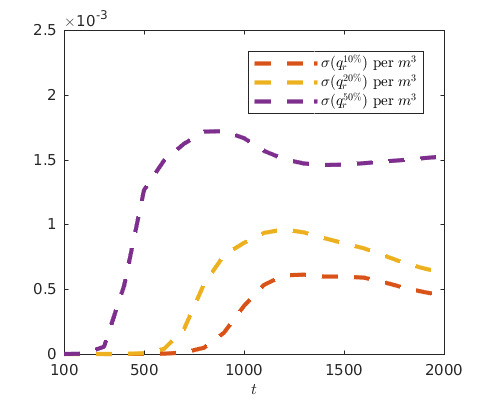}}
\caption{{\cblue Example 9: Time evolution of the expected values with their standard deviations for the cloud variables per $m^3$ (shaded
region, left column) and standard deviation (right column) using $0\%$, $10\%$, $20\%$ and $50\%$ (from top to down) perturbation of the
initial data in $q_v$.\label{fig35}}}
\end{figure}

Similarly to Example 7, we now investigate the uncertainty quantification in the 3-D Rayleigh-B\'enard convection for stochastically
perturbed initial data of the cloud variables given by \eqref{ic_rb}. The numerical solution is computed in a domain
$\Omega=[0,5000]\times[0,5000]\times[0,1000]\,\mathrm{m^3}$, which is discretized using $50\times50\times50$ mesh cells. {\cblue In Figure
\ref{fig33}, we show the influence of $0\%$, $10\%$, $20\%$ and $50\%$ perturbation on the potential temperature $\theta$ as well as the
expected value of the cloud drop concentration $q_c$. {\cblue For the perturbation scenarios in the case of the 3-D moist
Rayleigh-B\'{e}nard convection, one can see overall a qualitatively similar picture as for the corresponding 2-D simulations. However, due
to an additional spatial direction, the 3-D structures can further change their patterns. For positive perturbations, clouds can be formed
even at low vertical upward motions. The latent heat release increases the vertical motions in the convective cells, which leads to
additional feedback, such as stronger cloud formation, which in turn leads to formation of larger amount of rain water. As in the 2-D case,
the potential temperature distribution changes considerably at lower levels of the 3-D domain since evaporative cooling of precipitation is a
dominant process. Similarly, the mean cloud water distribution is crucially changed. For positive perturbations, which are dominant on average in
the beginning, more cloud water is formed and is later removed. Thus, less cloud water is available in the domain at later times. This can
also be seen in the time evolution of spatially averaged variables, as shown in Figure~\ref{fig35}, where the time evolution of the mean
expected value per $m^3$ and the mean standard deviation per $m^3$ of the cloud variables are plotted. The time evolution of these variables
is very similar to the 2-D case, but the 3-D scenarios show a very interesting feature, which is not available in 2-D. For larger
perturbations, the pattern of the convective cells is changed. While convective cells have a quasi hexagonal shape for the deterministic
simulation, they are ordered in a different way for larger perturbations. At upper levels, a more roll-like structure or even rectangular
shape is formed. Thus, there is a transition of structures due to perturbations in the initial conditions. This feature has not been
documented until now. A more detailed analysis of these structures is left for future studies.}}

\section{Conclusion}
In the present paper, we have studied uncertainty propagation in an atmospheric model that combines the Navier-Stokes equations for weakly
compressible fluids \eqref{NS_equations} with the cloud equations \eqref{cloud_equations}. The latter has been recently proposed in
\cite{porz18} and is based on the so-called single moment approach considering the evolution equations for the mass concentrations of the
water vapor, cloud drops and rain. Our numerical strategy is based on the stochastic Galerkin method that combines a finite-volume method
for space-time discretization with a spectral approximation in the stochastic space. We point out that atmospheric flows are weakly
compressible which leads to the low Mach number problem. One therefore needs to use a finite-volume method, which is accurate and efficient
in the low Mach number regime; see \cite{bispen1,diss}. To this end, we have chosen a suitable linear-nonlinear splitting between the fast
and slow flow variables and the second-order IMEX discretization in time (the ARS (2,2,2) scheme) as described in Section \ref{sec4}.
Coupling between the cloud model \eqref{cloud_equations} and the Navier-Stokes system \eqref{NS_equations_pert} is realized numerically by
the second-order Strang splitting. The cloud equations are approximated in space by the finite-volume method and in time using the explicit
third-order Runge-Kutta method with an enlarged stability region as explained in Section \ref{sec4}. Note that microscopic cloud dynamics
requires a smaller time step than the flow dynamics and thus several microscopic cloud subiterations are realized within one macroscopic
splitting time step, whose size is dictated by the flow dynamics. To the best of our knowledge, this is the first contribution that
combines an accurate and efficient method for the weakly compressible Navier-Stokes equations with the stochastic Galerkin method for the
uncertainty quantification of time evolution of the mass densities of water vapor, cloud drops and rain.

We have conducted extensive numerical benchmarking for both the deterministic and stochastic models and present the obtained numerical
results in Sections \ref{sec6} and \ref{sec7}. In the latter, we took into account the uncertainties in both initial data and cloud model
parameters. Our numerical study clearly demonstrates applicability of the stochastic Galerkin method for the uncertainty quantification in
complex atmospheric models. {\cblue We have obtained interesting results illustrating the behavior of clouds in different perturbed
scenarios and demonstrated that perturbations in the initial conditions can crucially change the time evolution of the moist
Rayleigh-B\'{e}nard convection. In particular, it has been shown that for larger perturbation, on average the positive perturbations
dominate the expectation values, although the standard deviation can be quite large. The main feature is a strong evaporative cooling in
lower levels of the (2-D/3-D) domain due to enhance rain formation and sedimentation into low humidity levels. In the 3-D case, a change in
the formed pattern can be seen, changing from hexagons to rolls/rectangles, which is quite surprising. For perturbations in parameters for
rain processes, the results are also dominated by the positive part of the perturbations. Since rain processes are affected, the spread in
the time evolution is increasing since the changes depend crucially on the interaction of formed rain drops with other variables. Overall,
it seems that for cloud physics, the expectation values are dominated by the positive perturbations leading to a change in the
distributions. This is an interesting topic for detailed studies in future.} Our further goal is to extend the developed numerical method to the
fully random Navier-Stokes-cloud system by considering random weakly compressible Navier-Stokes equations. {\cred We are also
interested in considering different random effects, such as initial data, boundary data and model parameters simultaneously, which would
require a multivariate stochastic Galerkin method.} This will allow one to quantify more precisely the propagation of small scale stochastic
errors initiated at cloud scales to macroscopic scales of flow dynamics.

\bigskip
\textbf{Acknowledgement}\\
The research leading to	these results has been done within the sub-project A2 of the Transregional Collaborative Research Center
SFB/TRR 165 ``Waves to Weather'' funded by the German Science Foundation (DFG). The work of A. Chertock was supported in part by NSF grants
DMS-1521051 and DMS-1818684. The work of A. Kurganov was supported in part by NSFC grant 11771201 and NSF grants DMS-1521009 and
DMS-1818666.

The authors gratefully acknowledge the support of the Data Center (ZDV) in Mainz for providing computation time on MOGON and MOGON~II
clusters and of NSF RNMS Grant DMS-1107444 (KI-Net). M. Luk\'a\v{c}ov\'a-Medvi{\softd}ov\'a and B. Wiebe would like to thank L. Yelash
(University of Mainz) for fruitful discussions. P. Spichtinger would like to thank P. Reutter (University of Mainz) for fruitful
discussions on the results of the 2-D/3-D simulations.

{\cblue
\section{Appendix A: Closure for single moment schemes}
The number concentration of rain drops can be approximated by a function of the respective mass concentration $n_r=f(q_r,c_r)$. Since we
implicitly assume that rain drops are distributed according to a size distribution, this approach should be used for mimicking the shape of
the distribution in a proper way. If we use a constant mean mass of rain drops, the function will be a simple linear relation
$n_r=\frac{1}{\overline{m}_r}q_r$. We extend this approach and propose the following nonlinear relation:
\begin{equation*}
n_r=c_r\cdot q_r^\gamma,\quad0<\gamma\le1.
\end{equation*}
Using this approach, one can replace the quantity $n_r$ in the processes related to rain drop number concentration. For the simple case of a
constant mean mass $\overline{m}_r$, we can determine the constants as $c_r=\overline{m}_r^{-1}$ and $\gamma=1$. This approach would be
meaningful for the case of a symmetric size distribution of rain droplets centered around the mean mass. However, it is well-known that size
distributions of rain are usually skew to larger sizes and thus a linear relation is inappropriate. For sizes of rain drops, an exponential
distribution is often assumed (see \cite{marshall_palmer1948}), namely:
\begin{equation*}
f(r)=B_re^{-\lambda r}
\end{equation*}
with a constant parameter $B_r=2\cdot 10^7\,\mathrm{m}^{-4}$ and the drop radius $r$. Using the general moments of the distribution,
\begin{equation*}
\mu_k[r]=\frac{\Gamma(k+1)}{\lambda^{k+1}}B_r
\end{equation*}
with the gamma function $\Gamma(x):=\int_0^\infty t^{x-1}\exp(-t)\,dt$, we obtain
\begin{equation*}
\rho n_r=\mu_0[r]=\frac{B_r}{\lambda}\quad\mbox{and}\quad\rho q_r=\frac{4}{3}\pi\rho_\ell\frac{\Gamma(4)}{\lambda^4}B_r.
\end{equation*}
Using these relations, one can derive the following function for the number concentration $n_r$:
\begin{equation*}
n_r=\underbrace{B_r^\frac{3}{4}\rho^{-\frac{3}{4}}\left(8\pi\rho_\ell\right)^{-\frac{1}{4}}}_{=c_r}q_r^\frac{1}{4}=c_rq_r^\gamma,\quad
\gamma=\frac{1}{4}.
\end{equation*}
We stress that $c_r$ is, in fact, a function of the air density $\rho$, that is, $c_r=c_{r0}\cdot \rho^{-\frac{3}{4}}$.
}

{\cred
\section{Appendix B: Explicit formulation of the cloud equations}
We present the equations of microphysical processes in an explicit way as they are used in our numerical experiments. In Tables
\ref{tab:appendix_physical_constants} and \ref{tab:appendix_model_parameters}, we present physical constants and model parameters with their
values used in our numerical simulations.

$$
\begin{aligned}
n_c&=q_c\frac{8\cdot10^8}{q_c+4.1888\cdot10^{-7}}\coth\left(\frac{q_c}{5.236\cdot10^{-13}}\right),\quad
C_{{\rm act}}=6.2832\cdot10^{-3}D_vG\rho\(q_v-q_*\)_+,\\
C_1&=0.7796D_vG\(q_v-q_*\)\left(\frac{8\cdot10^8}{q_c+4.1888\cdot10^{-7}}\coth\left(\frac{q_c}{5.236\cdot10^{-13}}\right)\right)^\frac{2}{3}
\rho q_c,\\
p_s(T)&=\exp\Big\{54.842763-\nicefrac{6763.22}{T}-4.21\ln T+0.000367T+\tanh(0.0415(T-218.8))\\
      &\quad\cdot(53.878-\nicefrac{1331.22}{T}-9.44523\ln T+0.014025T)\Big\},\\
n_r&=23752.6753\rho^{-\frac{3}{4}}q_r^{\frac{1}{4}},\quad
r=\left(\frac{1.21\cdot10^{-5}}{q_r+0.2874\rho^{-\frac{3}{4}}q_r^{\frac{1}{4}}}\right)^{\frac{4}{15}},\\
E&=-0.7796D_vG\(q_*-q_v\)_+\left(644.5198\sqrt{\rho q_r}+17.5904\mu^{-\frac{1}{6}}D_v^{-\frac{1}{3}}\sqrt{\alpha r}\rho^{\frac{13}{24}}
q_r^{\frac{91}{120}}\right),\\
A_1&=10^{-3}k_1\rho q_c^2,\quad A_2=0.3846\alpha k_2\rho^{\frac{1}{4}}q_crq_r^{\frac{61}{60}},\quad
v_q=1.1068\alpha q_r^{\frac{4}{15}}r\rho^{-\frac{1}{2}}.
\end{aligned}
$$
}

\begin{table}[ht!]
\centering
\begin{tabular}{ll}
Constant & Description \\
\hline
$p_*=\SI{101325}{\pascal}$ & reference pressure\\
$T_*=\SI{288}{\kelvin}$ & reference temperature\\
$T_0=\SI{273.15}{\kelvin}$ & melting temperature \\
$\rho_*=\SI{1.225}{\kilogram\per\cubic\meter}$ & reference air density\\
$\rho_l=\SI{1000}{\kilogram\per\cubic\meter}$ & density of liquid water\\
$R_v=\SI{461.52}{\joule\per\kilogram\per\kelvin}$ & specific gas constant, water vapor\\
$c_p=\SI{1005}{\joule\per\kilogram\per\kelvin}$ & specific heat capacity, dry air\\
$g=\SI{9.81}{\meter\per\square\second}$ & acceleration due to gravity\\
$L=\SI{2.53e6}{\joule\per\kilogram}$ & latent heat of vaporization\\
$\varepsilon=\frac{M_{\text{mol},v}}{M_{\text{mol},a}}=0.622$ & ratio of molar masses of water and dry air\\
$D_0=\SI{2.11e-5}{\square\meter\per\second}$ & diffusivity constant\\
\end{tabular}
\caption{Physical constants and reference quantities, \cite{porz18}.\label{tab:appendix_physical_constants}}
\end{table}

\begin{table}[ht!]
\centering
\begin{tabular}{ll}
Parameter & Description \\
\hline
$\alpha=190.3\pm0.5\cdot\SI{190.3}{\meter\per\second\raiseto{-\beta}\kilogram}$ & parameter for terminal velocity\\
$k_1=0.0041\pm0.5\cdot\SI{0.0041}{\kilogram\per\second}$ & parameter for autoconversion\\
$k_2=0.8\pm0.5\cdot\SI{0.8}{\kilogram}$ & parameter for accretion\\
$\beta=\frac{4}{15}$ & parameter for terminal velocity\\
$m_t=\SI{1.21e-5}{\kilogram}$ & parameter for terminal velocity\\
$N_0=\SI{1000}{\per\cubic\meter}$ & parameter for activation\\
$N_\infty=\SI{8e+8}{\per\kilogram}$ & parameter for activation\\
$m_0=\SI{5.236e-16}{\kilogram}$ & parameter for activation\\
$a_E=0.78$ & parameter for evaporation\\
$a_v=0.78$ & parameter for ventilation\\
$b_v=0.308$ & parameter for ventilation\\
\end{tabular}
\caption{Model parameters, \cite{porz18}.\label{tab:appendix_model_parameters}}
\end{table}


\bibliographystyle{plain}
\bibliography{biblio}	
	
\end{document}